\documentclass[preprint,12pt,review,numbers,sort&compress]{elsarticle}

\usepackage{amssymb}
\usepackage{subfigure}
\usepackage{caption}
\usepackage{multirow}
\usepackage{tgtermes}
\usepackage{amsmath}
\usepackage{placeins}
\usepackage{natbib}
\usepackage{array}
\usepackage{graphicx}
\usepackage{gensymb}
\usepackage{textcomp}
\usepackage{bm}
\usepackage{enumerate}
\usepackage[margin=2.5cm]{geometry}
\usepackage{color, soul}
\usepackage[justification=centering]{caption}
\usepackage{booktabs}
\usepackage{mathtools}

\allowdisplaybreaks[4]

\journal{ }

\begin{document}

\begin{frontmatter}
\sethlcolor{yellow}

\title{A new meshless fragile points method (FPM) with minimum unknowns at each point, for flexoelectric analysis under two theories with crack propagation. Part I: Theory and implementation}

\author[add1]{Yue Guan\corref{cor1}}
\ead{yuguan@ttu.edu}
\cortext[cor1]{Corresponding author.}

\author[add2]{Leiting Dong}
\author[add1]{Satya N. Atluri}

\address[add1]{Department of Mechanical Engineering, Texas Tech University, Lubbock, TX 79415, United States}
\address[add2]{School of Aeronautic Science and Engineering, Beihang University, Beijing 100191, China}

\begin{abstract}
``Flexoelectricity'' refers to a phenomenon which involves a coupling of the mechanical strain gradient and electric polarization. In this study, a meshless Fragile Points Method (FPM), is presented for analyzing flexoelectric effects in dielectric solids. Local, simple, polynomial and discontinuous trial and test functions are generated with the help of a local meshless Differential Quadrature approximation of derivatives. Both primal and mixed FPM are developed, based on two alternate flexoelectric theories, with or without the electric gradient effect and Maxwell stress. The first theory is fully nonlinear and is recommended at the nano-scale, while the second theory is linear and is sufficient at the micro-scale. In the present primal as well as mixed FPM, only the displacements and electric potential are retained as explicit unknown variables at each internal Fragile Point in the final algebraic equations. Thus the number of unknowns in the final system of algebraic equations is kept to be absolutely minimal. An algorithm for simulating crack initiation and propagation using the present FPM is presented, with classic stress-based criterion as well as a Bonding-Energy-Rate(BER)-based criterion for crack development. The present primal and mixed FPM approaches represent clear advantages as compared to the current methods for computational flexoelectric analyses, using primal as well as mixed Finite Element Methods, Element Free Galerkin (EFG) Methods, Meshless Local Petrov Galerkin (MLPG) Methods, and Isogeometric Analysis (IGA) Methods, because of the following new features: they are simpler Galerkin meshless methods using polynomial trial and test functions; minimal DoFs per Point make it very user-friendly;  arbitrary polygonal subdomains make it flexible for modeling complex geometries; the numerical integration of the primal as well as mixed FPM weak forms is trivially simple; and FPM can be easily employed in crack development simulations without remeshing or trial function enhancement. In this first part of the two-part series, we focus on the theoretical formulation and implementation of the proposed primal as well as mixed FPM. Numerical results and validation are then presented in Part II of the present paper.
\end{abstract}


\end{frontmatter}

\section{Introduction}

``Flexoelectricity'' describes an electromechanical coupling effect between the electric polarization and mechanical strain gradient \cite{Wang2019, Zhuang2020, Maranganti2006}. As a result of the recent trend of miniaturization of electromechanical systems, which necessitates the consideration of the significant effects of strain gradients, the flexoelectric effect has gained much attention in recent years. Resulting from the size effect, the flexoelectric response can become significant and even dominant in micro/nano-electromechanical systems \cite{Huang2012, Choi_2017, Liu2019}. Therefore, there has been an increased demand for reliable and accurate theories and numerical methods for flexoelectric analyses. In 2006, \citet{Maranganti2006} proposed the first continuum theory for flexoelectricity. After that, a number of advanced continuum theories have been developed, considering the Maxwell stress and surface effect \cite{Shen2010, Hu2010} and incorporating  kinetic energy \cite{PhysRevB.77.125424}, etc.

These continuum theories of flexoelectricity have been applied in analyzing mechanical and electrical responses of various systems. For example, \citet{Mao2014} derived the governing equations for a flexoelectric solid under small deformation and presented the analytical solutions for a 2D axisymmetric boundary value problem. The electric field gradient effect, as well as the converse flexoelectricity (a linear coupling effect between the electric field gradient and mechanical strain) are omitted in the governing equations. In the present paper, we denote these equations in \cite{Mao2014} as a ``reduced'' theory. On the contrary, the generalized form considering both the electric gradient effect and the electroelastic stress is signified as a ``full'' theory \cite{Maranganti2006, Zhuang2020}. Whereas the reduced theory is linear and much easier for numerical implementation, the influence of the electric gradient effect and the electroelastic stress can be significant at the nano-scale and thus the reduced theory may lead to inaccurate results \cite{Hu2010}. In this paper, both the nonlinear and linear theories are numerically implemented. The appropriate theory should be chosen depending on the length scale of the problem under consideration.

In practice, the governing equations and boundary value problems describing the flexoelectric effect can be formulated from two different approaches. First, as presented by \citet{Maranganti2006}, the internal energy density can be given as a function of strain, strain gradient, electric polarization and polarization gradient. Since the electric polarization ($P$) is used as an independent variable, this is denoted as a `P-formulation'. Alternatively, we can employ the electric Gibbs free energy to formulate the boundary value problem. Thus instead of the electric polarization and polarization gradient, the independent variables become electric field ($E$) and electric field gradient, and the resulting equations are known as an ‘E-formulation’. Both the formulations have been considered for numerical implementations. For instance, \citet{Yvonnet2017} and \citet{Mao2016} have proposed primal and mixed finite element (FEM) frameworks based on the P- formulation. Whereas the E-formulation can also be incorporated with multiple numerical discretization methods, including primal and mixed FEM \cite{Sladek2018, Deng2018}, Element-Free-Galerkin (EFG) method \cite{He2019} and Meshless Local Petrov-Galerkin (MLPG) method \cite{Sladek2020}, etc. The two formulations are inherently equivalent and can be unified by a linear relation between the electric polarization and electric displacement (the conjugate component of the electric field). Whereas the P-formulation employs the electric polarization induced by strain gradient intuitively, when the full theory is taken into consideration, the E-formulation is more favorable due to mathematical simplicity \cite{Zhuang2020}.. In the current work, we also concentrate on the E-formulation.

As have been stated, a number of numerical methods have been developed in analyzing flexoelectric and / or strain gradient effects in dielectric solids. These methods can be divided into three groups: the traditional Finite Element Methods (FEM), meshless methods, and isogeometric analysis (IGA) methods. Generally, the flexoelectric behavior is governed by a fourth-order partial differential equation (PDE). As a result, the main difficulty in modelling flexoelectricity for classic element-based methods lies in the $C^1$ continuity requirement. On one hand, the $C^1$ requirement can be satisfied by using Argyris triangular element \cite{Yvonnet2017} or subparametric quadrilateral element \cite{Beheshti2017, Sladek2018a}. However, even for two dimensional problems, the $C^1$ continuous elements result in large numbers of degrees of freedom (DoFs); for both of the approaches \cite{Yvonnet2017, Beheshti2017}, 9 DoFs are required at each node, including 6 mechanical quantities and 3 electrical quantities, thus making the formulation rather unfriendly for users. Besides, generating $C^1$ continuous elements in 3D can be extremely intractable and has not been pursued in any prior literature. Alternatively, mixed FEM can be developed using displacement gradients and Lagrangian multipliers as additional independent variables and requiring only $C^0$ continuity of all the variables. These mixed FEM have been applied in analyzing strain gradient \cite{Shu1999, Amanatidou2002, Bishay2012}, flexoelectric \cite{Mao2016, Deng2017} and converse flexoelectric effects \cite{Mao2016a}, and have also been extended to 3D analysis \cite{Deng2018} and topology optimization \cite{Nanthakumar2017}. Yet they have even more DoFs in each element (e.g., 87 DoFs in a 9-node element \cite{Mao2016}, 36 DoFs in a 5-node element \cite{Mao2016a}, or 37 DoFs in a 7-node element \cite{Deng2017}), making it prohibitively complex for practical use.

Another category of methods, known as ``meshless methods'', is partly or completely free of mesh discretization and have also shown their capability in analyzing flexoelectric or strain-gradient effects. \citet{Abdollahi2014} proposed a meshfree method with a $C^\infty$ continuous basis function using Local Maximum Entropy (LME) approximants. After that, the Element-Free Galerkin (EFG) Method based on the Moving Least Squares (MLS) approximation which can generate $C^1$ continuos functions through an appropriate choice of the weight functions, was applied to composite beam analysis with flexoelectricity \cite{Ray2017a, Ray2017, Basutkar2019}, as well as 2D strain gradient \cite{Sidhardh2018} and flexoelectric analysis \cite{He2019}. The Meshless Local Petrov-Galerkin (MLPG) Method is also established based on the MLS approximation for trial functions, whereas multiple test functions (e.g., weight function, local fundamental solution, Heaviside function, etc.) can be applied and thus lead to an asymmetric formulation. The MLPG method is also easily applicable in analyzing strain gradient \cite{Tang2003} and flexoelectric \cite{Sladek2020, Sladek2012} problems. The main advantage of the MLS approximation used in the EFG and MLPG methods is with an appropriate weight function, we can generate smooth trial functions with $C^1$ or higher order continuity. However, the trial functions are very complicated rational polynomials, and its higher order derivatives can be rational, and even more complicated. In order to avoid that problem, a mixed MLPG method, or Meshless Finite Volume Method (FVM) is proposed by \citet{Atluri2004}, with MLS trial functions and Heaviside test functions. Several other kinds of mixed MLPG methods were also developed for solving fourth order ODEs or PDEs arising in formulations involving  strain gradient effects \cite{Atluri2005, Jarak2020}. In these mixed MLPG methods, even though the mechanical strain and / or strain gradients are used as independent variables, these additional variables are eliminated at the node level and only the displacement DoFs are retained at the nodes  in the final formulations. Nevertheless, even though higher order derivatives are avoided, the trial function based on MLS approximation is still complicated in the mixed MLPG method, making it difficult for high-precision numerical integration. Similar argument can be made when Compactly Supported Radial Basis Function (CSRBF)  is used as the trial function when solving nonlinear flexoelectric problems \cite{Zhuang2019}.

All the meshless methods mentioned above have certain drawbacks: First, the shape functions based on MLS, LME lack the Kronecker delta property. Therefore, additional treatments are required in imposing the essential boundary conditions, e.g., a modified collocation method based on interior penalty functions \cite{Zhu1998} or Lagrange multipliers \cite{Basutkar2019}. The imposition of higher order essential boundary conditions in flexoelectric analysis has not been mentioned in any meshless approaches. Second, as a result of the complicated trial functions, the numerical integration of the weak form in either EFG or MLPG methods can be tedious. Actually, this difficulty in domain integration is a challenge for most Galerkin meshfree methods \cite{Hillman2016}. The simplest choice is direct nodal integration. Though efficient, and no background mesh required, the direct nodal integration may have stability issues and can be less accurate or non-convergent. High order quadrature, on the other hand, can achieve stability and better convergence. But it is computationally expensive and prohibitive to be used in practice. A number of modified nodal integration methods are then proposed, in order to ensure the accuracy and stability, which in turn sacrifice the efficiency \cite{Hillman2016}. Some other newer-type quadrature, e.g., a three-point integration scheme using background triangle elements \cite{Duan2012}, are also adopted. Nevertheless, these techniques are all applied to remedy the complicated integration in the weak form for the MLS trial functions.

Alternatively, the difficulty of numerical integration can be simply solved by employing simple shape functions like polynomials. The Fragile Points Method (FPM), in contrast to all the previous meshless methods, is based on simple, polynomial, piecewise-continuous trial and test functions \cite{Dong2019, Guan2020}. Consequently, the classic Gaussian quadrature can be employed, and numerical integrations are trivially simple. Only one integration point in each subdomain is sufficient most of the time in the FPM, whereas in methods like EFG and MLPG, even a large number of Gaussian integration points cannot guarantee a reliable numerical integration. Furthermore, the shape functions in the FPM pass through the nodal data, thus the Dirichlet boundary conditions can be enforced directly. On the other hand, when compared with the element-based methods, the trial functions of the FPM are generated using scattered points with the support of each point and always involve the same Cartesian strains in each subdomain. Thus, it has benefits in avoiding the mesh distortion and locking problems associated with element-based methods. Besides, unlike the FEM and other element-based methods which require triangular or quadrilateral element meshing, the FPM can be implemented with random point distributions and arbitrary subdomain shapes, and thus can deal with problems with complicated overall geometries. The FPM has already shown great accuracy and efficiency in solving 2D heat conduction \cite{Guan2020, Guan2020a} and elasticity problems \cite{Dong2019, Yang2019}. In this study, we formulate and apply the FPM for analyzing flexoelectric problems using both full and reduced theories. Furthermore, we can also develop mixed formulations for the FPM in which the mechanical strain and strain gradients are used as independent variables. As in the mixed MLPG method, the additional higher order variables can all be eliminated at each Point level. Therefore, the final mixed FPM formulations simply have nodal displacement and electric potentials, i.e., only 3 DoFs are retained at each node for a 2D problem in the final analysis. Comparisons of all the representative primal and mixed numerical methods are shown in Tables~\ref{table:Comp_P} and \ref{table:Comp_M} respectively.

Otherwise, the flexoelectric effect can also be analyzed by the IGA approach \cite{Liu2019, Do2019}, in which Non-Uniform Rational B-Spline (NURBS) basis functions with $C^1$ or higher order continuity are chosen to approximate both the geometry and field variables (displacement and electric potential fields). Yet the NURBS basis function also lacks Kronecker-delta property and thus, as in EFG, MLPG, etc., the essential boundary conditions need to be imposed by Lagrange multipliers \cite{Nguyen2015}. The IGA approach also requires a  large number of DoFs in each patch and thus making it rather unfriendly for users.

All the above studies just concentrate on continuous domains. Yet a significant strain gradient can be expected in the near-tip fields for a crack in elastic materials. Thus the strain gradient and flexoelectric effects should be taken into consideration for fracture analysis at micro and nano scales. \citeauthor{Sladek2018a} has analyzed the flexoelectric effect for stationary cracks in dielectric solids with the FEM \cite{Sladek2018a} and MLPG method \cite{Sladek2018b}. \citeauthor{Abdollahi2012} has studied the crack propagation for a four-point bending test in piezoelectric material [46], as well as the influence of flexoelectricity on a stationary crack \cite{Abdollahi2015} based on a meshless formulation. However, few studies have been reported on the crack propagation simulations considering the flexoelectric effect. This may stem from the difficulties in simulating crack propagation while using continuous trial functions. Approaches to simulate crack and rupture initiation and propagations, e.g. Generalized FEM (GFEM) \cite{OHara2016}, Extended FEM (XFEM) \cite{Fries2010, Giovanardi2017}, or Zencrack \cite{Hou2001}, usually involve mesh refinement or trial function enrichment. Nevertheless, in the present FPM, as a result of the discontinuous trial and test functions, the algorithm for simulating crack propagation is relatively simple and intuitive, and does not involve either mesh refinement or trial function enrichment.

In this paper, we focus on analyzing the flexoelectric effect in dielectric solids using the FPM. Both full and reduced flexoelectric theories are considered and their corresponding equations are introduced in section~\ref{sec:theory}. The primal and mixed FPM formulations are developed in sections~\ref{sec:primal} and \ref{sec:mixed} respectively. Section~\ref{sec:Crack} presents the algorithm for simulating crack initiation and propagation while using the FPM. Numerical results and validation of the overall methodologies, as well as a parametric discussion, are given in Part II of this 2-part paper.

\begin{table}[htbp]
\caption{A comparison of primal numerical methods.}
\centering
{
\begin{tabular}{ | m{2.25cm} | m{2.25cm} | m{2.25cm} |  m{2.25cm} | m{2.25cm} | m{2.25cm} |}
\hline 
Method & Primal FPM & Primal FEM & Primal EFG & Primal MLPG & IGA \\
\hline
Trial functions & Discontinuous polynomial & $C^1$ 5th-order polynomial & MLS, etc. & MLS, etc. & NURBS \\
\hline
Test functions & As above & As above & As above & MLS, weight function, local fundamental solution, Dirac delta function, Heaviside function, $\cdots$ & As above \\
\hline
Delta function property & Yes & Yes & Yes/No & Yes/No & No \\
\hline
Weak form & Galerkin (present paper) or Petrov-Galerkin (future) & Galerkin & Galerkin & Petrov-Galerkin & Galerkin \\
\hline
Weak form integration & Simple (Gaussian quadrature) & Simple (Gaussian quadrature) & Very complicated & Very complicated & Very complicated\\
\hline
DoF at each node & 3 & 9 & 6 & 3 & 3\\
\hline
\end{tabular}}
\label{table:Comp_P}
\end{table}

\begin{table}[htbp]
\caption{A comparison of mixed numerical methods.}
\centering
{
\begin{tabular}{ | m{2.8cm} | m{3.73cm} | m{3.73cm} |  m{3.73cm} |  }
\hline 
Method & Mixed FPM & Mixed FEM & Mixed MLPG \\
\hline
Trial functions & Discontinuous linear polynomial & $C^0$ linear/quadratic polynomial &  MLS, etc. \\
\hline
Test functions &  As above & As above &  MLS, weight function, local fundamental solution, Dirac delta function, Heaviside function, etc. \\
\hline
Delta function property & Yes & Yes & Yes/No \\
\hline
Weak form & Galerkin (present paper) or Petrov-Galerkin (future) & Galerkin & Petrov-Galerkin \\
\hline
Weak form integration & Very simple (one-point quadrature) & Simple (Gaussian quadrature) & Very complicated\\
\hline
DoF at each node & 3 & 3 - 13  & 3\\
\hline
\end{tabular}}
\label{table:Comp_M}
\end{table}

\section{Flexoelectricity theories and boundary value problems} \label{sec:theory}

\subsection{Full theory}

In the current work, we consider a homogeneous elastic dielectric 2D domain $\Omega$ with boundary $\partial \Omega$. The spatial coordinate is given as $\mathbf{x} = \left[ x_1, x_2 \right]^\mathrm{T}$. Subjected to mechanical and electrical loadings, the responses of the material are described by the displacement vector field $\mathbf{u} \left( x_1, x_2 \right) = \left[ u_1, u_2\right]^\mathrm{T}$ and the electric potential field $\phi \left( x_1, x_2 \right)$.

Phenomenologically, the flexoelectric effect describes an electric polarization generated by the mechanical strain gradient:
\begin{align}
\begin{split}
P_i = \overline{\mu}_{ijkl} \kappa_{jkl}.
\end{split}
\end{align}
Einstein summation convention is used here. $\mathbf{P} = \left[ P_1, P_2\right]^\text{T}$ is the electrical polarization vector, $\kappa_{jkl} = u_{j,kl}$ is the second gradient of displacement (i.e., the strain gradient), and $\overline{\mu}_{ijkl}$ are elements of a fourth order flexoelectric tensor ($i, j, k, l = 1,2$).

The present boundary value problem is formulated from an electric Gibbs free energy $G$, in which the mechanical strain $\boldsymbol{\varepsilon}$, strain gradient $\boldsymbol{\kappa}$, electric field $\mathbf{E}$ and electric field gradient $\mathbf{V}$ are independent variables:
\begin{align} \label{eqn:define}
\begin{split}
\boldsymbol{\varepsilon} & = \left[ \varepsilon_{11}, \varepsilon_{22}, 2 \varepsilon_{12} \right]^\mathrm{T}, \\
\boldsymbol{\kappa} & = \left[ \kappa_{111}, \kappa_{222}, \kappa_{122}, \kappa_{211}, 2 \kappa_{112}, 2 \kappa_{212} \right]^\mathrm{T}, \\
\mathbf{E} & = \left[ E_1, E_2 \right]^\mathrm{T}, \\
\mathbf{V} & = \left[ V_{11}, V_{21}, V_{12}, V_{22} \right]^\mathrm{T}.
\end{split}
\begin{split}
& \varepsilon_{ij} = \frac{1}{2} \left( u_{i,j} + u_{j,i} \right), \\
& \kappa_{ijk} = u_{i,jk}, \\
& E_i = -\phi_{,i}, \\
& V_{ij} = E_{i,j}.
\end{split}
\end{align}

The governing equations are given as:
\begin{align}
\left( \sigma_{ij} - \mu_{ijk,k} + \sigma_{ij}^{ES} \right)_{,j} + b_i & = 0, \label{eq:gov_full_dis} \\
\left( D_i - Q_{ij,j} \right)_{,i} & = q, \label{eq:gov_full_ele}
\end{align}
where $b_i$ is the body force per volume, $q$ is the free charge per volume, $\sigma_{ij}$, $\mu_{ijk}$, $D_i$ and $Q_{ij}$ are the (local) stress, higher-order stress, electric displacement, and higher-order electric displacement, i.e., conjugate variables of $\varepsilon_{ij}$, $\kappa_{ijk}$, $E_i$ and $V_{ij}$ respectively. $\sigma_{ij}^{ES}$ is the generalized electrostatic stress, which can be written as:
\begin{align}
\begin{split}
\sigma_{ij}^{ES} & = \sigma_{ij}^M + \mu_{ijm,m}^M, \\
\sigma_{ij}^M & = E_i D_j - \frac{1}{2} E_k D_k \delta_{ij}, \\
\mu_{ijm,m}^M & = V_{ik} Q_{kj} - \frac{1}{2} V_{kl} Q_{kl} \delta_{ij}.
\end{split}
\end{align}
where $\sigma_{ij}^M$ and $\mu_{ijm,m}^M$ are the Maxwell stress and the higher-order electrostatic stress respectively.

The corresponding Dirichlet boundary conditions are:
\begin{align}
u_i & = \widetilde{u}_i, & \text{on} \; \partial \Omega_u, \label{eqn:BC_full_1} \\
\phi & = \widetilde{\phi}, & \text{on} \; \partial \Omega_{\phi}, \\
\nabla^n u_i = u_{i,l} n_l & = \widetilde{d}_i, & \text{on} \; \partial \Omega_{d}, \\
\nabla^n \phi = \phi_{,l} n_l & = \widetilde{P}, & \text{on} \; \partial \Omega_{P}.
\end{align}

The Neumann boundary conditions are:
\begin{align}
\left( \sigma_{ij} - \mu_{ijk,k} + \sigma_{ij}^{ES} \right) n_j + \nabla_k^t \left( n_l \right) n_m n_j \mu_{ijm} - \nabla_j^t \left( n_m \mu_{ijm} \right) & = \widetilde{Q}_i, & \text{on} \; \partial \Omega_Q,  \label{eqn:BC_full_5}\\
\left( D_i - Q_{ij,j} \right) n_i + \nabla_l^t \left( n_l \right) n_i n_j Q_{ij} - \nabla_i^t \left( Q_{ij} n_j\right) & = - \widetilde{\omega}, & \text{on} \; \partial \Omega_{\omega}, \\
\mu_{ijm} n_m n_j & = \widetilde{R}_i, & \text{on} \; \partial \Omega_{R}, \\
Q_{ij} n_i n_j & = \widetilde{Z}, & \text{on} \; \partial \Omega_{Z}.  \label{eqn:BC_full_8}
\end{align}

In the above equations, $\widetilde{\mathbf{u}}$, $\widetilde{\phi}$, $\widetilde{\mathbf{d}}$, $\widetilde{P}$, $\widetilde{\mathbf{Q}}$, $\widetilde{\omega}$, $\widetilde{\mathbf{R}}$ and $\widetilde{Z}$ are known functions. $\mathbf{n} = \left[ n_1, n_2 \right]^\mathrm{T}$ is the outward
unit normal vector to $\partial \Omega$. $\nabla^n = \mathbf{n} \cdot \boldsymbol{\nabla}$ and $\boldsymbol{\nabla}^t = \boldsymbol{\nabla} - \mathbf{n} \nabla^n$ are the normal and tangential surface differential operators respectively. $\partial \Omega_u \cup \partial \Omega_Q = \partial \Omega_{\phi} \cup \partial \Omega_{\omega} = \partial \Omega_d \cup \partial \Omega_R = \partial \Omega_P \cup \partial \Omega_Z = \partial \Omega$, and $\partial \Omega_u \cap \partial \Omega_Q = \partial \Omega_{\phi} \cap \partial \Omega_{\omega} = \partial \Omega_d \cap \partial \Omega_R = \partial \Omega_P \cap \partial \Omega_Z = \varnothing$.

The constitutive relations are given in the matrix forms:
\begin{align}
\begin{split}
\boldsymbol{\sigma} & = \mathbf{C}_{\sigma \varepsilon} \boldsymbol{\varepsilon} - \mathbf{e} \mathbf{E} - \mathbf{b} \mathbf{V}, \\
\boldsymbol{\mu} & = \mathbf{C}_{\mu \kappa} \boldsymbol{\kappa} - \mathbf{a} \mathbf{E}, \\
\mathbf{D} & = \boldsymbol{\Lambda} \mathbf{E} + \mathbf{e}^\mathrm{T} \boldsymbol{\varepsilon} + \mathbf{a}^\mathrm{T} \boldsymbol{\kappa}, \\
\mathbf{Q} & = \boldsymbol{\Phi} \mathbf{V} + \mathbf{b}^\mathrm{T} \boldsymbol{\varepsilon}.
\label{eqn:cons_full}
\end{split}
\end{align}
where
\begin{align}\nonumber
\begin{split}
\boldsymbol{\sigma} & = \left[ \sigma_{11}, \sigma_{22}, \sigma_{12} \right]^\mathrm{T}, \\
\boldsymbol{\mu} & = \left[ \mu_{111}, \mu_{222}, \mu_{122}, \mu_{211}, \mu_{112}, \mu_{212} \right]^\mathrm{T}, \\
\mathbf{D} & = \left[ D_1, D_2\right]^\mathrm{T}, \\
\mathbf{Q} & = \left[ Q_{11}, Q_{21}, Q_{12}, Q_{22} \right]^\mathrm{T}.
\end{split}
\end{align}
$\mathbf{C}_{\sigma \varepsilon}$, $\mathbf{C}_{\mu \kappa}$, $\boldsymbol{\Lambda}$, $\boldsymbol{\Phi}$, $\mathbf{a}$, $\mathbf{b}$ and $\mathbf{e}$ are matrices of material properties. The generalized electrostatic stress $\boldsymbol{\sigma}^{ES}$ can also be written in a matrix form:
\begin{align} \label{eqn:ele_stress}
\begin{split}
\boldsymbol{\sigma}^{ES} & = \left[ \sigma_{11}^{ES}, \sigma_{22}^{ES}, \sigma_{12}^{ES}, \sigma_{21}^{ES} \right]^\mathrm{T} \\
& = \widehat{\mathbf{D}} \mathbf{E} + \widehat{\mathbf{Q}} \mathbf{V},
\end{split}
\end{align}
where
\begin{align}\nonumber
\begin{split}
\widehat{\mathbf{D}} = \left[ \begin{matrix} \frac{1}{2} D_1 &  -\frac{1}{2} D_2 \\   -\frac{1}{2} D_1 & \frac{1}{2} D_2 \\ D_2 & 0 \\ 0 & D_1 \end{matrix} \right], \quad
 \widehat{\mathbf{Q}} = \left[ \begin{matrix} \frac{1}{2}Q_{11} & -\frac{1}{2} Q_{21} & Q_{21} - \frac{1}{2} Q_{12} & -\frac{1}{2} Q_{22} \\
 -\frac{1}{2} Q_{11} & Q_{12} - \frac{1}{2} Q_{21} & -\frac{1}{2} Q_{12} & \frac{1}{2} Q_{22} \\
 Q_{12} & 0 & Q_{22} & 0 \\
 0 & Q_{11} & 0 & Q_{21} \end{matrix} \right].
\end{split}
\end{align}

\subsection{Reduced theory}

Now, some different formulations can be reduced from the above-mentioned generalized formulation in the absence of electrostatic stress $\boldsymbol{\sigma}^{ES}$ and / or of the electric field gradient $\mathbf{V}$. Here we present a reduced formulation \cite{Mao2014} in which both $\boldsymbol{\sigma}^{ES}$ and $\mathbf{V}$ are omitted. The governing equations and boundary conditions in the reduced form are:
\begin{align}
\left( \sigma_{ij} - \mu_{ijk,k} \right)_{,j} + b_i & = 0,  \label{eqn:gov_red_dis}\\
 D_{i,i} & = q, \label{eqn:gov_red_ele}
\end{align}
Dirichlet boundary conditions:
\begin{align}
u_i & = \widetilde{u}_i, & \text{on} \; \partial \Omega_u, \\
\phi & = \widetilde{\phi}, & \text{on} \; \partial \Omega_{\phi}, \\
\nabla^n u_i = u_{i,l} n_l & = \widetilde{d}_i, & \text{on} \; \partial \Omega_{d},
\end{align}
Neumann boundary conditions:
\begin{align}
\left( \sigma_{ij} - \mu_{ijk,k} + \sigma_{ij}^{ES} \right) n_j + \nabla_k^t \left( n_l \right) n_m n_j \mu_{ijm} - \nabla_j^t \left( n_m \mu_{ijm} \right) & = \widetilde{Q}_i, & \text{on} \; \partial \Omega_Q, \\
D_i n_i & = - \widetilde{\omega}, & \text{on} \; \partial \Omega_{\omega}, \\
\mu_{ijm} n_m n_j & = \widetilde{R}_i, & \text{on} \; \partial \Omega_{R}.  \label{eqn:BC_red_6}
\end{align}
And the corresponding constitutive relations are:
\begin{align}
\begin{split}
\boldsymbol{\sigma} & = \mathbf{C}_{\sigma \varepsilon} \boldsymbol{\varepsilon} - \mathbf{e} \mathbf{E}, \\
\boldsymbol{\mu} & = \mathbf{C}_{\mu \kappa} \boldsymbol{\kappa} - \mathbf{a} \mathbf{E}, \\
\mathbf{D} & = \boldsymbol{\Lambda} \mathbf{E} + \mathbf{e}^\mathrm{T} \boldsymbol{\varepsilon} + \mathbf{a}^\mathrm{T} \boldsymbol{\kappa}.
\end{split}
\end{align}

As the electrostatic stress $\boldsymbol{\sigma}^{ES}$ is omitted, this reduced form is a linear formulation, and thus is convenient for numerical implementation. It remains accurate in micro or larger scale structures with negligible electric field gradient. However, studies have shown that the significance of the electric force increases dramatically with the minimization of the system size \cite{Hu2010}. Therefore, for nano-structures or structures with significant electric field gradient, the full theory including the electrostatic stress should be taken into consideration.

Note that an alternative formulation (``P-formulation'') based on an internal energy $\widetilde{U}$ which used mechanical strain $\boldsymbol{\varepsilon}$, strain gradient $\boldsymbol{\kappa}$, electrical polarization $\mathbf{P}$ and polarization field $\mathbf{P}_{,i}$ as independent variables are also commonly used in previous literatures \cite{Maranganti2006, Yvonnet2017}. Which formulation being used is simply a matter of choice, since the two formulations are essentially equivalent. And the relation between the electric displacement, polarization and field are given as:
\begin{align}
\begin{split}
\mathbf{D} = \epsilon_0 \mathbf{E} + \mathbf{P},
\end{split}
\end{align}
where $\epsilon_0$ is the permittivity of free space. 

\section{The Primal Fragile Points Method (primal FPM)}\label{sec:primal}

\subsection{Trial and test functions}\label{sec:func}

\subsubsection{Points and domain partition}

\begin{figure}[htbp] 
  \centering 
    \subfigure[]{ 
    \label{fig:Schem_Voro} 
    \includegraphics[width=0.48\textwidth]{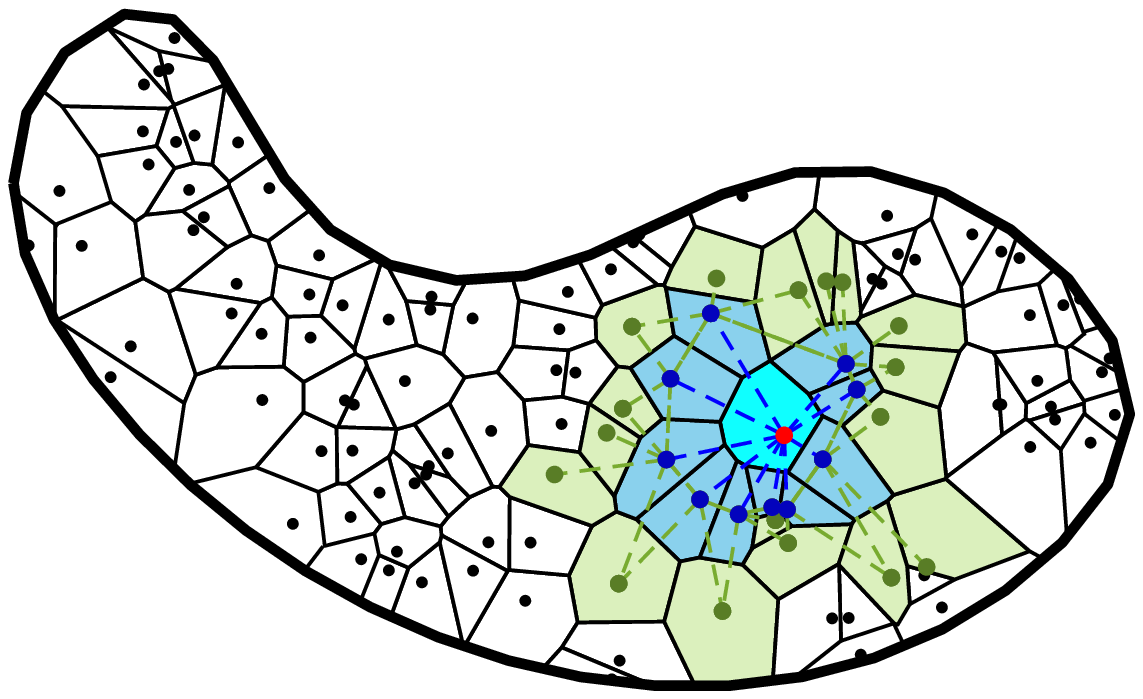}}  
    \subfigure[]{ 
    \label{fig:Schem_ABAQUS} 
    \includegraphics[width=0.48\textwidth]{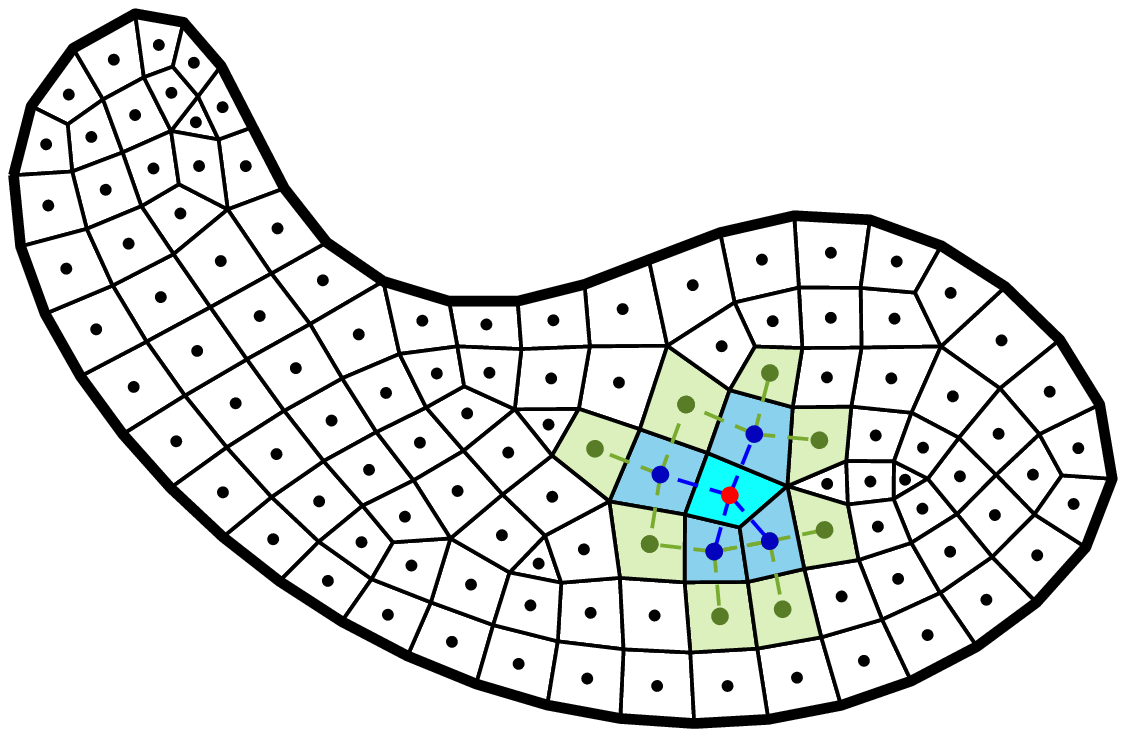}}  
  \caption{The domain $\Omega$ and its partitions. (a) Voronoi Diagram partition. (b) Quadrilateral and triangular partition (converted from FEA meshing).} 
  \label{fig:Schem} 
\end{figure}

In the Fragile Points Method (FPM), first, a set of random Points are scattered in the domain. The entire domain can then be partitioned into several nonoverlapping subdomains. Within each subdomain, only one Fragile Point exists. Multiple partitioning schemes may be employed, e.g., the Voronoi Diagram partition (see Fig.~\ref{fig:Schem_Voro}), as well as quadrilateral and triangular partition. The traditional FEM meshing can also be employed. The FEM elements can be converted into FPM subdomains, while the internal Fragile Point is defined as the centroid of each  FEM geometrical element. Thus, the preprocessing module of any commercial FEA software like ABAQUS can be helpful to generate the Points and subdomains in this proposed FPM. For example, Fig.~\ref{fig:Schem_ABAQUS} exhibits a typical point distribution and subdomain partition in FPM as converted form the ABAQUS meshing.

However, unlike in the FEM, the trial and test functions in the present FPM are {\it point-based}. In the primal FPM, in order to describe the high-order gradient-dependent  behavior, the trial functions for the displacement and electric potential (($\mathbf{u}^h$ and $\phi^h$) in each subdomain are written in the form of a third-order Taylor expansion at the corresponding internal point. For instance, in the FPM subdomain $E^0$ which contains an internal point $P_0$:
\begin{align} \label{eqn:Shape}
\begin{split}
\mathbf{u}^h \left( x, y \right) = \mathbf{u}^0 & + \left( x - x_0 \right) \mathbf{u}_{,1} \Big|_{P_0} + \left( y - y_0 \right) \mathbf{u}_{,2} \Big|_{P_0} \\
& + \frac{1}{2} \left( x - x_0 \right)^2 \mathbf{u}_{,11} \Big|_{P_0}  + \left( x  - x_0 \right) \left( y - y_0 \right) \mathbf{u}_{,12} \Big|_{P_0}  + \frac{1}{2} \left( y - y_0 \right)^2 \mathbf{u}_{,22} \Big|_{P_0} \\
& + \frac{1}{6} \left( x - x_0 \right)^3 \mathbf{u}_{,111} \Big|_{P_0} + \frac{1}{2} \left( x - x_0 \right)^2 \left( y - y_0 \right) \mathbf{u}_{,112} \Big|_{P_0} \\
& + \frac{1}{2} \left( x - x_0 \right) \left( y - y_0 \right)^2 \mathbf{u}_{,122} \Big|_{P_0}  + \frac{1}{6} \left( y - y_0 \right)^3 \mathbf{u}_{,222} \Big|_{P_0} 
\end{split}
\end{align}
where $\left[ x_0, y_0 \right]$ are the coordinates of point $P_0$. $\mathbf{u}^0$ is the value of $\mathbf{u}^h$ at $P_0$. The first to third derivatives $\left[ \overline{\mathbf{D}} \mathbf{u} \right]  \Big|_{P_0} =  \left[ \mathbf{u}_{,1}^\mathrm{T}, \mathbf{u}_{,2}^\mathrm{T}, \mathbf{u}_{,11}^\mathrm{T}, \mathbf{u}_{,12}^\mathrm{T}, \mathbf{u}_{,22}^\mathrm{T}, \mathbf{u}_{,111}^\mathrm{T}, \mathbf{u}_{,112}^\mathrm{T}, \mathbf{u}_{,122}^\mathrm{T}, \mathbf{u}_{,222}^\mathrm{T} \right]^\mathrm{T} \Big|_{P_0}$ are as yet unknown. Here we introduce the local Differential Quadrature method to determine these high-order derivatives in terms of the displacements at a few surrounding support Points.

\subsubsection{Local radial basis function-based Differential Quadrature method} \label{sec:DQ}

The conventional differential quadrature (DQ) method is a widely used numerical discretization technique. However, usually based on 1D basis functions (e.g., polynomials), the conventional DQ method can only be applied along a mesh-line and cannot be implemented directly with meshless schemes. In recent years, several multi-dimensional DQ methods have been developed based on the meshless Radial Basis Functions (RBFs). In the current work, we introduce and modify the local RBF-DQ method proposed by \citet{Shu2003} in 2003 which approximates the 2D derivatives directly using a limited number of supporting points.

In the current methodology, as shown in Fig.~\ref{fig:Schem}, for a given point $P_0  \in E_0$ (red), its support is defined to involve all its nearest (blue) and second (green) neighboring points. Here the nearest neighboring points are the points in subdomains sharing boundaries with $E_0$, while the second neighboring points are nearest neighbors of the former. For the points on or close to the domain boundary $\partial \Omega$, the third neighboring points are also taken into consideration to remedy the lack of effective supporting points. These supporting points are named as $P_1$, $P_2$, $\cdots$, $P_m$.

Many RBFs can be used as basis function in the local RBF-DQ method. The multiquadric (MQ), inverse-MQ and Gaussians RBFs are all available to incorporate within the FPM approach, and achieve similar accuracy with appropriate parameters. Here we use the MQ-RBF as an example:
\begin{align}
\begin{split}
\psi \left( r \right) = \sqrt{r^2 + c^2},
\end{split}
\end{align}
where $r$ is the radial length from the conference point. $c = c_0 d_0$ is a constant parameter, where $d_0$ is the diameter of the minimal-diameter circle enclosing all the supporting Points.

Thus, for supporting Points $P_i \left( \mathbf{x}_i \right)$, the approximation of the displacement field $u_j \left( x, y \right) \; \left( j = 1, 2\right)$ can be written as:
\begin{align}
\begin{split}
u_j \left( x, y \right) = \sum_{i = 0}^m \lambda_i \psi \left( \left\| \mathbf{x} - \mathbf{x}_i \right\|_2 \right) + \lambda_{m+1} + \lambda_{m+2} x + \lambda_{m+3} y,
\end{split}
\end{align}
in which the coefficients $\lambda_i$ are determined by using the $\left( m + 1 \right)$ collocating equations at $P_i$ and the following conditions:
\begin{align}
\begin{split}
\sum_{i = 0}^m \lambda_i = 0, \quad \sum_{i = 0}^m \lambda_i x_i = 0, \quad \sum_{i = 0}^m \lambda_i y_i = 0.
\end{split}
\end{align}

In the DQ method, any partial derivative of displacement at point $P_0$ can be approximated by a weighted linear sum of the value of $\mathbf{u}^h$  at all its supporting points:
\begin{align}
\begin{split}
\frac{\partial^{s+t}}{\partial x^s \partial y^t} \mathbf{u}^h \Big|_{P_0} = \sum_{i=0}^m W_i^{\left( s, t \right)} \mathbf{u}^h \big|_{P_i},
\end{split}
\end{align}
where $W_i^{\left( s,t \right)}$ is the weighting coefficient corresponding to point $P_i$. After rearrangement, the derivatives under study at point $P_0$ can be approximated as:
\begin{align}
\begin{split}
\left[ \overline{\mathbf{D}} \mathbf{u} \right]  \Big|_{P_0} =  \left( \overline{\mathbf{B}} \otimes \mathbf{I}_{2 \times 2} \right) \mathbf{u}_E ,
\end{split}
\end{align}
where
\begin{align} \nonumber
\begin{split}
\mathbf{u}_E=\left[ u_1^0, u_2^0, u_1^1, u_2^1, \cdots , u_1^m, u_2^m \right]^\mathrm{T}.
\end{split}
\end{align}
$u^i = \left[ u_1^i, u_2^i \right]^\mathrm{T}$ is the value of $\mathbf{u}^h$  at $P_i$. $\overline{\mathbf{B}}$ is the matrix of weighting coefficients $W_i^{\left( s, t \right)}$, $\mathbf{I}_{2 \times 2}$ is a unit matrix, and $\otimes$ denotes the Kronecker product. The weighting coefficient matrix $\overline{\mathbf{B}}$ can be obtained from:
\begin{align}
\begin{split}
\overline{\mathbf{B}}=\mathbf{G}^{-1} \left[ \overline{\mathbf{D}} \mathbf{G} \right],
\end{split}
\end{align}
where
\begin{align} \nonumber
\begin{split}
\mathbf{G} & = \left[ \begin{matrix} 1 & 1 & \cdots & 1 \\
x_0 & x_1 & \cdots & x_m \\
y_0 & y_1 & \cdots & y_m \\
g_3 \left( x_0, y_0 \right) & g_3 \left( x_1, y_1 \right) & \cdots & g_3 \left( x_m, y_m \right) \\
\vdots & \vdots & \ddots & \vdots \\
g_m \left( x_0, y_0 \right) & g_m \left( x_1, y_1 \right) & \cdots & g_m \left( x_m, y_m \right) \end{matrix} \right], 
\end{split} \\
\begin{split} \nonumber
\left[ \overline{\mathbf{D}} \mathbf{G} \right] & = \left[ \begin{matrix} 0 & 1 & 0 & g_{3,1} \left( x_0, y_0 \right) & \cdots & g_{m,1} \left( x_0, y_0 \right) \\
0 & 0 & 1 & g_{3,2} \left( x_0, y_0 \right) & \cdots & g_{m,2} \left( x_0, y_0 \right) \\
0 & 0 & 0 & g_{3,11} \left( x_0, y_0 \right) & \cdots & g_{m,11} \left( x_0, y_0 \right) \\
0 & 0 & 0 & g_{3,12} \left( x_0, y_0 \right) & \cdots & g_{m,12} \left( x_0, y_0 \right) \\
0 & 0 & 0 & g_{3,22} \left( x_0, y_0 \right) & \cdots & g_{m,22} \left( x_0, y_0 \right) \\
0 & 0 & 0 & g_{3,111} \left( x_0, y_0 \right) & \cdots & g_{m,111} \left( x_0, y_0 \right) \\
0 & 0 & 0 & g_{3,112} \left( x_0, y_0 \right) & \cdots & g_{m,112} \left( x_0, y_0 \right) \\
0 & 0 & 0 & g_{3,122} \left( x_0, y_0 \right) & \cdots & g_{m,122} \left( x_0, y_0 \right) \\
0 & 0 & 0 & g_{3,222} \left( x_0, y_0 \right) & \cdots & g_{m,222} \left( x_0, y_0 \right) \end{matrix} \right]^\mathrm{T}, \\
g_i \left( x, y \right) & = \psi_i \left( x, y \right) - \frac{ x_1 y_2 - x_2 y_1 + x_i \left( y_1 - y_2 \right) + y_i \left( x_2 - x_1 \right) }{\left( x_1 - x_0 \right) \left( y_2 - y_0 \right) - \left( x_2 - x_0 \right) \left( y_1 - y_0 \right)} \psi_0 \left( x, y \right)  \\
& \qquad \qquad \; \;  \, - \frac{ x_2 y_0 - x_0 y_2 + x_i \left( y_2 - y_0 \right) + y_i \left( x_0 - x_2 \right) }{\left( x_1 - x_0 \right) \left( y_2 - y_0 \right) - \left( x_2 - x_0 \right) \left( y_1 - y_0 \right)} \psi_1 \left( x, y \right) \\ 
& \qquad \qquad \; \;  \, - \frac{ x_0 y_1 - x_1 y_0 + x_i \left( y_0 - y_1 \right) + y_i \left( x_1 - x_0 \right) }{\left( x_1 - x_0 \right) \left( y_2 - y_0 \right) - \left( x_2 - x_0 \right) \left( y_1 - y_0 \right)} \psi_2 \left( x, y \right), \\
\psi_i \left( x, y \right) & = \sqrt{\left( x - x_i \right)^2 + \left( y - y_i \right)^2 + c^2}.
\end{split}
\end{align}

Note that the accuracy of local RBF-DQ approximation is excellent for the first derivative; whereas, for higher-order approximations, the accuracy decreases. This implies that a mixed formulation may help to improve the accuracy of the current primal FPM.

\subsubsection{Local, polynomial, discontinuous test and trial functions}

\begin{figure}[htbp] 
  \centering 
    \subfigure[]{ 
    \label{fig:Shape_3rd} 
    \includegraphics[width=0.48\textwidth]{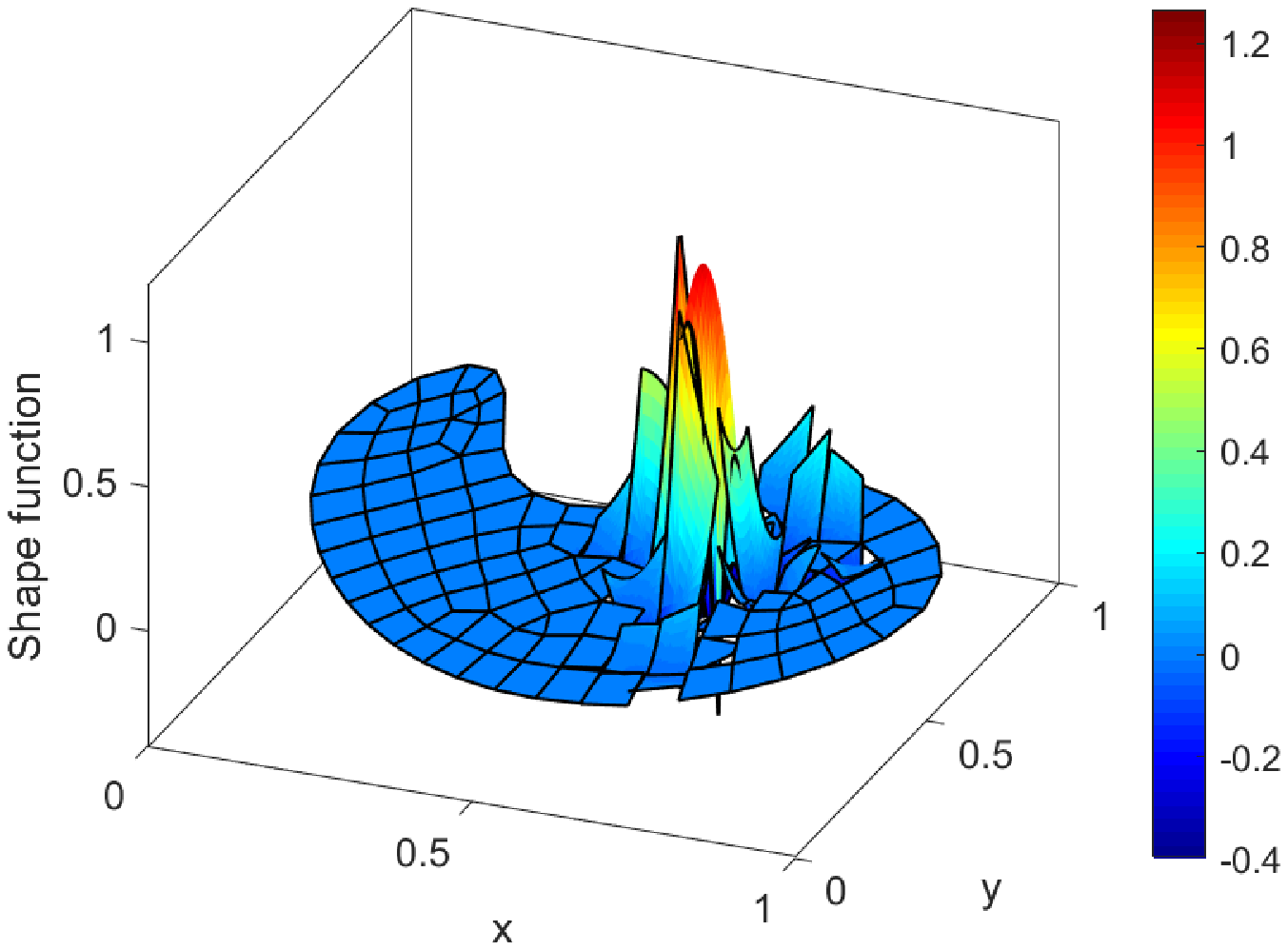}}  
    \subfigure[]{ 
    \label{fig:Trial_3rd} 
    \includegraphics[width=0.48\textwidth]{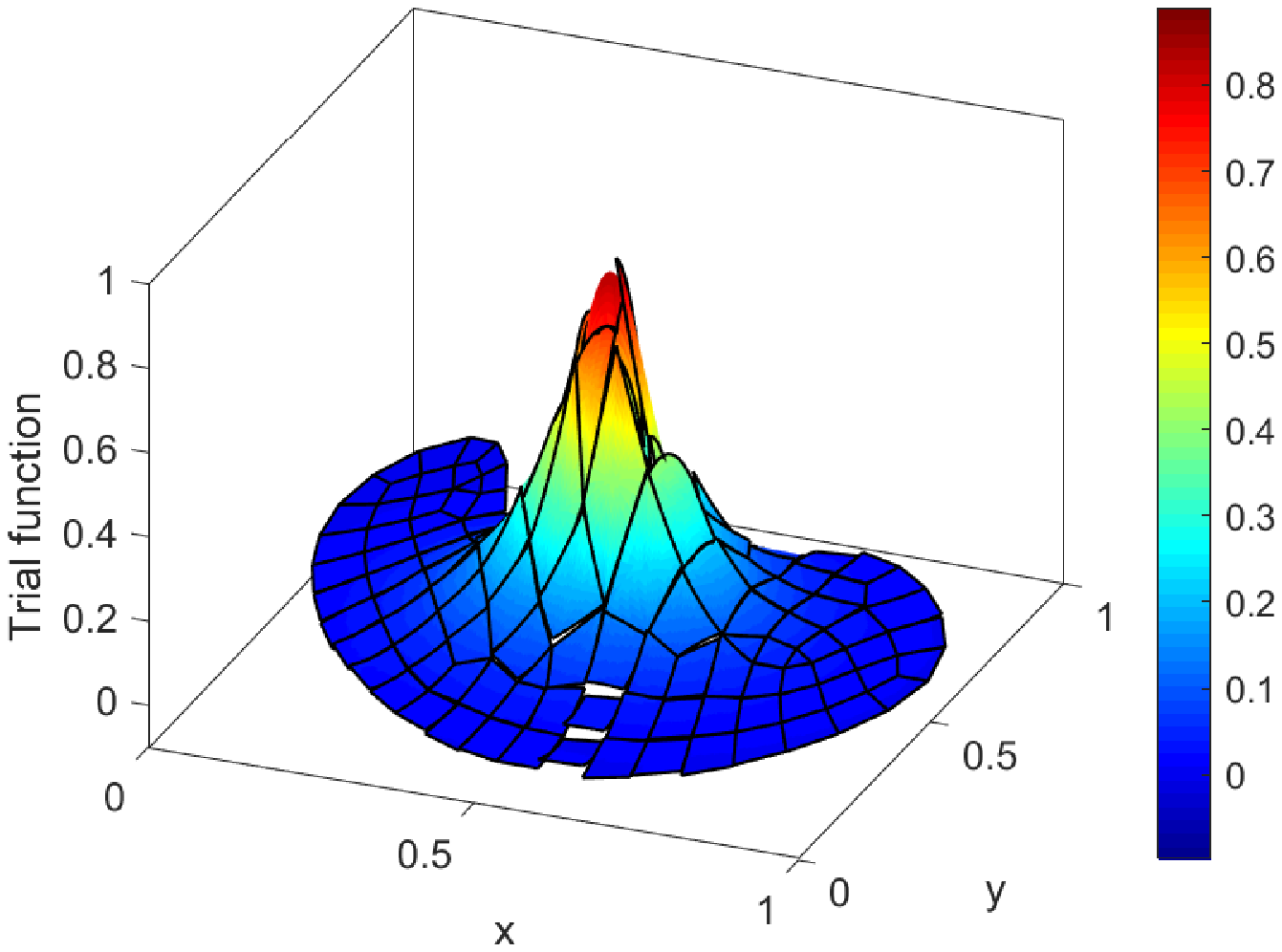}}  
  \caption{The shape and trial function in the primal FPM. (a) The shape function. (b) The trial function for $u_a = e^{-10 \left[ \left( x – 0.5 \right)^2 + \left( y – 0.5 \right)^2 \right]^{1/2}}$.} 
  \label{fig:Func_3rd} 
\end{figure}

Substituting the approximation of the derivatives into Eqn.~\ref{eqn:Shape}, the trial function for displacement $\mathbf{u}^h$ in the subdomain $E_0$ is achieved (see Eqn.~\ref{eq:Dis_shape}). Similarly, we can also generate the shape and trial functions for the electric potential $\phi^h$:
\begin{align}
\mathbf{u}^h \left( x, y \right) & = \mathbf{N}_u \left( x, y \right) \mathbf{u}_E, \label{eq:Dis_shape} \\
\phi^h \left( x, y \right) & = \mathbf{N}_\phi \left( x, y \right) \boldsymbol{\phi}_E,
\end{align}
where
\begin{align} \nonumber
\begin{split}
& \boldsymbol{\phi}_E = \left[ \phi^0, \phi^1, \cdots , \phi^m \right]^\mathrm{T}, \\
& \mathbf{N}_u \left( x, y \right) = \mathbf{N}_\phi \left( x, y \right) \otimes \mathbf{I}_{2 \times 2}, \quad \mathbf{N}_\phi \left( x, y \right) = \overline{\mathbf{N}} \left( x, y \right) \cdot \overline{\mathbf{B}} + \left[ 1, 0, \cdots, 0 \right]_{1 \times \left( m + 1 \right)}, \\
& \overline{\mathbf{N}} \left( x, y \right) = \Big[ x - x_0, y - y_0, \frac{1}{2} \left( x - x_0 \right)^2,  \left( x - x_0 \right) \left( y - y_0 \right), \frac{1}{2} \left( y - y_0 \right)^2, \cdots \\
& \qquad \qquad \qquad \frac{1}{6} \left( x - x_0 \right)^3,  \frac{1}{2} \left( x - x_0 \right)^2  \left( y - y_0 \right), \frac{1}{2} \left( x - x_0 \right)  \left( y - y_0 \right)^2, \frac{1}{6} \left( y - y_0 \right)^3 \Big].
\end{split}
\end{align}
The trial function and shape function are defined in each subdomain, following the same process. Since no continuity requirement is as yet enforced at the internal boundaries, the shape and trial functions can be discontinuous. Fig.~\ref{fig:Shape_3rd} shows the graph of a shape function for the domain and partition shown in Fig.~\ref{fig:Schem_ABAQUS} with 116 points. The trial function simulating an exponential function $u_a = e^{-10 \left[ \left( x – 0.5 \right)^2 + \left( y – 0.5 \right)^2 \right]^{1/2}}$ is presented in Fig.~\ref{fig:Trial_3rd}. As can be seen, the trial function is a cubic polynomial in each subdomain. It is piecewise-continuous. However, as the shape function in each subdomain also depends on the values of the neighboring Points, the entire shape function shows a weakly continuous tendency. The test function for displacement $\mathbf{v}^h$ and electrical potential $\tau_h$ in the Galerkin weak-form in the FPM are prescribed to possess the same shape as $\mathbf{u}^h$ and $\phi_h$ respectively. As a result of the discontinuous trial functions, the FPM also has a great natural potential in analyzing systems involving cracks, ruptures and fragmentations.

Unfortunately, if the above trial and test functions are applied in the traditional Galerkin weak form directly, the discontinuity will lead to inconsistent and inaccurate results. In order to resolve that problem, we introduce the Numerical Flux Corrections to the present FPM.

\subsection{Weak-form formulation and numerical flux corrections}

\subsubsection{Full theory} \label{sec:PFPM_full}

We multiply the governing equations (Eqn.~\ref{eq:gov_full_dis} and \ref{eq:gov_full_ele}) by the test functions $\mathbf{v}_i^h$ and $\tau^i$ respectively in each subdomain $E$, then apply the Gauss divergence theorem:
\begin{align}
\begin{split}
&\int_E v_{i,j} \sigma_{ij} \mathrm{d} \Omega +  \int_E v_{i,jk} \mu_{ijk} \mathrm{d} \Omega + \int_E v_{i,j} \sigma_{ij}^{ES} \mathrm{d} \Omega \\
& \qquad = \int_E v_{i} b_i \mathrm{d} \Omega +  \int_{\partial E} v_{i} \left( \sigma_{ij} - \mu_{ijk,k} + \sigma_{ij}^{ES} \right) n_j \mathrm{d} \Gamma + \int_{\partial E} v_{i,j} \mu_{ijk} n_k \mathrm{d} \Gamma,
\end{split} \\
& - \int_E \tau_{,i} D_i \mathrm{d} \Omega - \int_E \tau_{,ij} Q_{ij} \mathrm{d} \Omega  \\
& \qquad = \int_E \tau q \mathrm{d} \Omega - \int_{\partial E} \tau \left(  D_i - Q_{ij,j} \right) n_i \mathrm{d} \Gamma - \int_{\partial E} \tau_{,i} Q_{ij} n_j \mathrm{d} \Gamma.
\end{align}
where $\partial E$ is the boundary of the subdomain $E$.

Let $\Gamma_h$ denotes the set of all internal boundaries, and $\Gamma = \Gamma_h + \partial \Omega = \Gamma_h + \partial \Omega_u +\partial \Omega_Q = \Gamma_h + \partial \Omega_\phi +\partial \Omega_\omega = \Gamma_h + \partial \Omega_d +\partial \Omega_R = \Gamma_h + \partial \Omega_P +\partial \Omega_Z$ is the set of all internal and external boundaries. Each internal boundary $e \in \Gamma_h$ is shared by two subdomains $E^1$ and $E^2$, i.e., $e = \partial E^1 \cap \partial E^2$ (the order does not affect the formulation of FPM). With $\mathbf{n}_i$ being the unit vector normal to $e$ and pointing outward from $E^i$, we define $\mathbf{n}^e = \mathbf{n}^1 = -\mathbf{n}^2$ for $e \in \Gamma_h$. And for $e \in \partial \Omega$, $\mathbf{n}^e = \mathbf{n}^1$, where $\mathbf{n}^1$ is the outward unit normal vector of the only one neighboring subdomain $E^1$. By summing the above equations over the entire domain, we rewrite them in the matrix’s forms:
\begin{align}
\begin{split}
& \int_\Omega \boldsymbol{\varepsilon}^\mathrm{T} \left( \mathbf{v} \right) \boldsymbol{\sigma} \left( \mathbf{u}, \phi \right) \mathrm{d} \Omega +  \int_\Omega \boldsymbol{\kappa}^\mathrm{T} \left( \mathbf{v} \right) \boldsymbol{\mu} \left( \mathbf{u}, \phi \right) \mathrm{d} \Omega + \int_\Omega \hat{\boldsymbol{\varepsilon}}^\mathrm{T} \left( \mathbf{v} \right) \boldsymbol{\sigma}^{ES} \left( \mathbf{u}, \phi \right) \mathrm{d} \Omega \\
& \qquad -\int_{\Gamma_h} \left[ \! \left[ \hat{\boldsymbol{\varepsilon}}^\mathrm{T} \left( \mathbf{v} \right) \right] \! \right] \overline{\mathbf{n}}_3^e \left\{ \boldsymbol{\mu} \left( \mathbf{u}, \phi \right) \right\}  \mathrm{d} \Gamma -\int_{\Gamma_h} \left\{ \hat{\boldsymbol{\varepsilon}}^\mathrm{T} \left( \mathbf{v} \right) \right\} \overline{\mathbf{n}}_3^e \left[ \! \left[ \boldsymbol{\mu} \left( \mathbf{u}, \phi \right) \right] \! \right]  \mathrm{d} \Gamma\\
& \qquad   -\int_{\Gamma_h} \left[ \! \left[ \mathbf{v}^\mathrm{T} \right] \! \right] \left\{ \overline{\mathbf{n}}_1^e \boldsymbol{\sigma} \left( \mathbf{u}, \phi \right) - \overline{\mathbf{n}}_{21}^e \boldsymbol{\mu}_{,1} \left( \mathbf{u}, \phi \right) - \overline{\mathbf{n}}_{22}^e \boldsymbol{\mu}_{,2} \left( \mathbf{u}, \phi \right) + \overline{\mathbf{n}}_{4}^e \boldsymbol{\sigma}^{ES} \left( \mathbf{u}, \phi \right)  \right\}  \mathrm{d} \Gamma \\
& \qquad   -\int_{\Gamma_h} \left\{ \mathbf{v}^\mathrm{T} \right\} \left[ \! \left[ \overline{\mathbf{n}}_1^e \boldsymbol{\sigma} \left( \mathbf{u}, \phi \right) - \overline{\mathbf{n}}_{21}^e \boldsymbol{\mu}_{,1} \left( \mathbf{u}, \phi \right) - \overline{\mathbf{n}}_{22}^e \boldsymbol{\mu}_{,2} \left( \mathbf{u}, \phi \right) + \overline{\mathbf{n}}_{4}^e \boldsymbol{\sigma}^{ES} \left( \mathbf{u}, \phi \right)  \right] \! \right]  \mathrm{d} \Gamma \\
& \qquad  = \int_\Omega \mathbf{v}^\mathrm{T} \mathbf{b} \mathrm{d} \Omega + \int_{\partial \Omega}  \left[ \! \left[ \hat{\boldsymbol{\varepsilon}}^\mathrm{T} \left( \mathbf{v} \right) \right] \! \right] \overline{\mathbf{n}}_3^e \left\{ \boldsymbol{\mu} \left( \mathbf{u}, \phi \right) \right\}  \mathrm{d} \Gamma \\
& \qquad  + \int_{\partial \Omega} \left[ \! \left[ \mathbf{v}^\mathrm{T} \right] \! \right] \left\{ \overline{\mathbf{n}}_1^e \boldsymbol{\sigma} \left( \mathbf{u}, \phi \right) - \overline{\mathbf{n}}_{21}^e \boldsymbol{\mu}_{,1} \left( \mathbf{u}, \phi \right) - \overline{\mathbf{n}}_{22}^e \boldsymbol{\mu}_{,2} \left( \mathbf{u}, \phi \right) + \overline{\mathbf{n}}_{4}^e \boldsymbol{\sigma}^{ES} \left( \mathbf{u}, \phi \right)  \right\}  \mathrm{d} \Gamma,
\label{eqn:dis_1}
\end{split} \\
\begin{split}
& \int_\Omega \mathbf{E}^\mathrm{T} \left( \tau \right) \mathbf{D} \left( \mathbf{u}, \phi \right) \mathrm{d} \Omega + \int_\Omega \mathbf{V}^\mathrm{T} \left( \tau \right) \mathbf{Q} \left( \mathbf{u}, \phi \right) \mathrm{d} \Omega \\
& \qquad  + \int_{\Gamma_h} \left[ \! \left[ \tau \right] \! \right] \left\{ \mathbf{n}^{e \mathrm{T}} \mathbf{D} \left( \mathbf{u}, \phi \right) - \overline{\mathbf{n}}_{51}^e  \mathbf{Q}_{,1} \left( \mathbf{u}, \phi \right) - \overline{\mathbf{n}}_{52}^e  \mathbf{Q}_{,2} \left( \mathbf{u}, \phi \right) \right\}  \mathrm{d} \Gamma \\
& \qquad  + \int_{\Gamma_h} \left\{ \tau \right\} \left[ \! \left[ \mathbf{n}^{e \mathrm{T}} \mathbf{D} \left( \mathbf{u}, \phi \right) - \overline{\mathbf{n}}_{51}^e  \mathbf{Q}_{,1} \left( \mathbf{u}, \phi \right) - \overline{\mathbf{n}}_{52}^e  \mathbf{Q}_{,2} \left( \mathbf{u}, \phi \right) \right] \! \right]  \mathrm{d} \Gamma \\
& \qquad  - \int_{\Gamma_h} \left[ \! \left[ \mathbf{E}^\mathrm{T} \left( \tau \right) \right] \! \right] \overline{\mathbf{n}}_{6}^e \left\{ \mathbf{Q} \left( \mathbf{u}, \phi \right) \right\}  \mathrm{d} \Gamma - \int_{\Gamma_h} \left\{ \mathbf{E}^\mathrm{T} \left( \tau \right) \right\} \overline{\mathbf{n}}_{6}^e \left[ \! \left[ \mathbf{Q} \left( \mathbf{u}, \phi \right) \right] \! \right]  \mathrm{d} \Gamma \\
& \qquad = \int_\Omega \tau q \mathrm{d} \Omega - \int_{\partial \Omega} \left[ \! \left[ \tau \right] \! \right] \left\{ \mathbf{n}^{e \mathrm{T}} \mathbf{D} \left( \mathbf{u}, \phi \right) - \overline{\mathbf{n}}_{51}^e  \mathbf{Q}_{,1} \left( \mathbf{u}, \phi \right) - \overline{\mathbf{n}}_{52}^e  \mathbf{Q}_{,2} \left( \mathbf{u}, \phi \right) \right\}  \mathrm{d} \Gamma \\
& \qquad + \int_{\partial \Omega} \left[ \! \left[ \mathbf{E}^\mathrm{T} \left( \tau \right) \right] \! \right] \overline{\mathbf{n}}_{6}^e \left\{ \mathbf{Q} \left( \mathbf{u}, \phi \right) \right\}  \mathrm{d} \Gamma.
\label{eqn:ele_1}
\end{split}
\end{align}
where
\begin{align} \label{eqn:define_0}
\hat{\boldsymbol{\varepsilon}}^\mathrm{T} \left( \mathbf{u} \right) = \left[ u_{1,1}, u_{2,2}, u_{1,2}, u_{2,1} \right]^\mathrm{T},
\end{align}
The jump operator $\left[ \! \left[ \cdot \right] \! \right]$ and average operator $\left\{ \cdot \right\}$ are defined as (for $\forall w$):
\begin{align}\label{eqn:operator}
\begin{split}
\left[ \! \left[ w \right] \! \right] = \begin{cases}
w \Big|_e^{E^1} -  w \Big|_e^{E^2} & e \in \Gamma_h\\
w \Big|_e & e \in \partial \Omega
\end{cases}, \quad
\left\{ w \right\} = \begin{cases}
\frac{1}{2} \left( w \Big|_e^{E^1} +  w \Big|_e^{E^2} \right) & e \in \Gamma_h\\
w \Big|_e & e \in \partial \Omega
\end{cases}.
\end{split}
\end{align}

The following boundary terms in Eqn.~\ref{eqn:dis_1} and \ref{eqn:ele_1} can be divided into normal and tangent components and then integrated by parts:
\begin{align}
\begin{split}
& \int_{\partial \Omega}  \left[ \! \left[ \hat{\boldsymbol{\varepsilon}}^\mathrm{T} \left( \mathbf{v} \right) \right] \! \right] \overline{\mathbf{n}}_3^e \left\{ \boldsymbol{\mu} \left( \mathbf{u}, \phi \right) \right\}  \mathrm{d} \Gamma = \int_{\partial \Omega_d \cup \partial \Omega_R} \hat{\boldsymbol{\varepsilon}}^\mathrm{T} \left( \mathbf{v} \right) \overline{\mathbf{n}}_4^{e \mathrm{T}} \overline{\mathbf{n}}_4^e \overline{\mathbf{n}}_3^e \boldsymbol{\mu} \left( \mathbf{u}, \phi \right) \mathrm{d} \Gamma \\
& \qquad - \int_{\partial \Omega_u \cup \partial \Omega_Q} \mathbf{v}^\mathrm{T} \overline{\mathbf{c}}_1 \overline{\mathbf{s}}_4^{e \mathrm{T}} \overline{\mathbf{s}}_4^e \overline{\mathbf{n}}_3^e \boldsymbol{\mu}_{,1} \left( \mathbf{u}, \phi \right) \mathrm{d} \Gamma - \int_{\partial \Omega_u \cup \partial \Omega_Q} \mathbf{v}^\mathrm{T} \overline{\mathbf{c}}_2 \overline{\mathbf{s}}_4^{e \mathrm{T}} \overline{\mathbf{s}}_4^e \overline{\mathbf{n}}_3^e \boldsymbol{\mu}_{,2} \left( \mathbf{u}, \phi \right) \mathrm{d} \Gamma \\
& \qquad + \sum_{C_h} \left[ \! \left[ \mathbf{v}^\mathrm{T} \right] \! \right] \overline{\mathbf{s}}_4^e \overline{\mathbf{n}}_3^e  \left\{ \boldsymbol{\mu} \left( \mathbf{u}, \phi \right)  \right\} + \sum_{C_h} \left\{ \mathbf{v}^\mathrm{T} \right\} \overline{\mathbf{s}}_4^e \overline{\mathbf{n}}_3^e  \left[ \! \left[ \boldsymbol{\mu} \left( \mathbf{u}, \phi \right)  \right] \! \right] ,
\end{split} \label{eqn:dis_2} \\
&  \int_{\partial \Omega} \left[ \! \left[ \mathbf{E}^\mathrm{T} \left( \tau \right) \right] \! \right] \overline{\mathbf{n}}_{6}^e \left\{ \mathbf{Q} \left( \mathbf{u}, \phi \right) \right\}  \mathrm{d} \Gamma = \int_{\partial \Omega_P \cup \partial \Omega_Z}  \mathbf{E}^\mathrm{T} \left( \tau \right) \mathbf{n}^e \mathbf{n}^{e \mathrm{T}} \overline{\mathbf{n}}_6^e \mathbf{Q} \left( \mathbf{u}, \phi \right) \mathrm{d} \Gamma \\
& \qquad + \int_{\partial \Omega_\phi \cup \partial \Omega_\omega} \tau \overline{\mathbf{c}}_3 \mathbf{s}^e \mathbf{s}^{e \mathrm{T}} \overline{\mathbf{n}}_6^e \mathbf{Q}_{,1} \left( \mathbf{u}, \phi \right) \mathrm{d} \Gamma + \int_{\partial \Omega_\phi \cup \partial \Omega_\omega} \tau \overline{\mathbf{c}}_4 \mathbf{s}^e \mathbf{s}^{e \mathrm{T}} \overline{\mathbf{n}}_6^e \mathbf{Q}_{,2} \left( \mathbf{u}, \phi \right) \mathrm{d} \Gamma \\
& \qquad - \sum_{C_h} \left[ \! \left[ \tau \right] \! \right] \mathbf{s}^{e \mathrm{T}} \overline{\mathbf{n}}_6^e  \left\{ \mathbf{Q} \left( \mathbf{u}, \phi \right)  \right\} - \sum_{C_h} \left\{ \tau \right\}  \mathbf{s}^{e \mathrm{T}} \overline{\mathbf{n}}_6^e \left[ \! \left[ \mathbf{Q} \left( \mathbf{u}, \phi \right)  \right] \! \right]. \label{eqn:ele_2}
\end{align}
where $C_h$ denotes the set of all edges of the external boundaries, i.e., $C_h = \underset{e \in \partial \Omega} {\bigcup} \partial e$, $\overline{\mathbf{s}}^e$ is the unit tangent vector pointing from $E^1$ to $E^2$ (the order is the same as the jump and average operators defined in Eqn.~\ref{eqn:operator}), and $\overline{\mathbf{s}}^{e} \cdot \overline{\mathbf{n}}^{e} = 0$. The corresponding matrices in Eqn.~\ref{eqn:dis_1} -- \ref{eqn:ele_2} are shown in \ref{app: Mats}.

When $\left( \mathbf{u}, \phi \right)$ are exact solutions, due to continuity conditions, we have:
\begin{align} \nonumber
\begin{split}
& \left[ \! \left[ \boldsymbol{\mu} \left( \mathbf{u}, \phi \right) \right] \! \right] = \mathbf{0}, \\ 
& \left[ \! \left[ \overline{\mathbf{n}}_1^e \boldsymbol{\sigma} \left( \mathbf{u}, \phi \right) - \overline{\mathbf{n}}_{21}^e \boldsymbol{\mu}_{,1} \left( \mathbf{u}, \phi \right) - \overline{\mathbf{n}}_{22}^e \boldsymbol{\mu}_{,2} \left( \mathbf{u}, \phi \right) + \overline{\mathbf{n}}_{4}^e \boldsymbol{\sigma}^{ES} \left( \mathbf{u}, \phi \right) \right] \! \right] = \mathbf{0}, \\ 
& \left[ \! \left[ \mathbf{n}^{e \mathrm{T}} \mathbf{D} \left( \mathbf{u}, \phi \right) - \overline{\mathbf{n}}_{51}^e \mathbf{Q}_{,1} \left( \mathbf{u}, \phi \right) - \overline{\mathbf{n}}_{52}^e \mathbf{Q}_{,2} \left( \mathbf{u}, \phi \right) \right] \! \right]  = \mathbf{0} \\
& \left[ \! \left[ \mathbf{Q} \left( \mathbf{u}, \phi \right) \right] \! \right] = \mathbf{0},
\end{split},\quad \text{for} \; \forall e \in \Gamma_h.
\end{align} 
Hence, we can eliminate the corresponding terms in Eqn.~\ref{eqn:dis_1} and \ref{eqn:ele_1}. Similarly, $\left[ \! \left[ \mathbf{u} \right] \! \right] = \mathbf{0}$, $\left[ \! \left[ \hat{\boldsymbol{\varepsilon}}^\mathrm{T} \left( \mathbf{u} \right) \right] \! \right] = \mathbf{0}$, $\left[ \! \left[ \phi \right] \! \right] = 0$ and $\left[ \! \left[ \mathbf{E}^\mathrm{T} \left( \phi \right) \right] \! \right] = \mathbf{0}$. We substitute Eqn.~\ref{eqn:ele_stress} into Eqn.~\ref{eqn:dis_1}. Then, by adding antithetic terms and equivalent terms to both sides, the previous equations can be written in a symmetric form:
\begin{align}
\begin{split}
& \int_\Omega \boldsymbol{\varepsilon}^\mathrm{T} \left( \mathbf{v} \right) \boldsymbol{\sigma} \left( \mathbf{u}, \phi \right) \mathrm{d} \Omega +  \int_\Omega \boldsymbol{\kappa}^\mathrm{T} \left( \mathbf{v} \right) \boldsymbol{\mu} \left( \mathbf{u}, \phi \right) \mathrm{d} \Omega \\
& \qquad   -\int_{\Gamma_h \cup \partial \Omega_u} \left[ \! \left[ \mathbf{v}^\mathrm{T} \right] \! \right] \left\{ \begin{matrix} \overline{\mathbf{n}}_1^e \boldsymbol{\sigma} \left( \mathbf{u}, \phi \right) - \overline{\mathbf{n}}_{21}^e \boldsymbol{\mu}_{,1} \left( \mathbf{u}, \phi \right) - \overline{\mathbf{n}}_{22}^e \boldsymbol{\mu}_{,2} \left( \mathbf{u}, \phi \right) \\ + \overline{\mathbf{n}}_{4}^e \widehat{\mathbf{D}} \left( \mathbf{u}, \phi \right) \mathbf{E} \left( \phi \right) + \overline{\mathbf{n}}_{4}^e  \widehat{\mathbf{Q}} \left( \mathbf{u}, \phi \right) \mathbf{V} \left( \phi \right)  \end{matrix} \right\}  \mathrm{d} \Gamma \\
& \qquad -\int_{\Gamma_h \cup \partial \Omega_u}  \left\{ \boldsymbol{\sigma}^\mathrm{T} \left( \mathbf{v}, 0 \right) \overline{\mathbf{n}}_1^{e \mathrm{T}} - \boldsymbol{\mu}_{,1}^\mathrm{T} \left( \mathbf{v}, 0 \right) \overline{\mathbf{n}}_{21}^{e \mathrm{T}} - \boldsymbol{\mu}_{,2}^\mathrm{T} \left( \mathbf{v}, 0 \right) \overline{\mathbf{n}}_{22}^{e \mathrm{T}} \right\} \left[ \! \left[ \mathbf{u} \right] \! \right] \mathrm{d} \Gamma \\
& \qquad  - \int_{\Gamma_h \cup \partial \Omega_\phi} \left\{ \mathbf{D}^\mathrm{T} \left( \mathbf{v}, 0 \right) \mathbf{n}^{e} - \mathbf{Q}_{,1}^\mathrm{T} \left( \mathbf{v}, 0 \right) \overline{\mathbf{n}}_{51}^{e \mathrm{T}}  - \mathbf{Q}_{,2}^\mathrm{T} \left( \mathbf{v}, 0 \right) \overline{\mathbf{n}}_{52}^{e \mathrm{T}} \right\} \left[ \! \left[ \phi \right] \! \right] \mathrm{d} \Gamma \\
& \qquad -\int_{\Gamma_h} \left[ \! \left[ \hat{\boldsymbol{\varepsilon}}^\mathrm{T} \left( \mathbf{v} \right) \right] \! \right] \overline{\mathbf{n}}_3^e \left\{ \boldsymbol{\mu} \left( \mathbf{u}, \phi \right) \right\}  \mathrm{d} \Gamma -\int_{\Gamma_h} \left\{ \boldsymbol{\mu}^\mathrm{T} \left( \mathbf{v}, 0 \right) \right\}  \overline{\mathbf{n}}_3^{e \mathrm{T}} \left[ \! \left[ \hat{\boldsymbol{\varepsilon}} \left( \mathbf{u} \right) \right] \! \right]   \mathrm{d} \Gamma \\
& \qquad + \int_{\Gamma_h} \left\{ \mathbf{Q}^\mathrm{T} \left( \mathbf{v}, 0 \right) \right\} \overline{\mathbf{n}}_{6}^{e \mathrm{T}} \left[ \! \left[ \mathbf{E} \left( \phi \right) \right] \! \right] \mathrm{d} \Gamma \\
& \qquad + \int_{\partial \Omega_u} \mathbf{v}^\mathrm{T} \left[ \overline{\mathbf{c}}_1 \overline{\mathbf{s}}_4^{e \mathrm{T}} \overline{\mathbf{s}}_4^e \overline{\mathbf{n}}_3^e \boldsymbol{\mu}_{,1} \left( \mathbf{u}, \phi \right) + \overline{\mathbf{c}}_2 \overline{\mathbf{s}}_4^{e \mathrm{T}} \overline{\mathbf{s}}_4^e \overline{\mathbf{n}}_3^e \boldsymbol{\mu}_{,2} \left( \mathbf{u}, \phi \right) \right] \mathrm{d} \Gamma \\
& \qquad + \int_{\partial \Omega_u} \left[  \boldsymbol{\mu}_{,1}^\mathrm{T} \left( \mathbf{v}, 0 \right) \overline{\mathbf{n}}_3^{e \mathrm{T}} \overline{\mathbf{s}}_4^{e \mathrm{T}} \overline{\mathbf{s}}_4^e  \overline{\mathbf{c}}_1^\mathrm{T} + \boldsymbol{\mu}_{,2}^\mathrm{T} \left( \mathbf{v}, 0 \right) \overline{\mathbf{n}}_3^{e \mathrm{T}} \overline{\mathbf{s}}_4^{e \mathrm{T}} \overline{\mathbf{s}}_4^e  \overline{\mathbf{c}}_2^\mathrm{T} \right] \mathbf{u} \mathrm{d} \Gamma \\
& \qquad + \int_{\partial \Omega_\phi} \left[ \mathbf{Q}_{,1}^\mathrm{T} \left( \mathbf{v}, 0 \right) \overline{\mathbf{n}}_6^{e \mathrm{T}} \mathbf{s}^e \mathbf{s}^{e \mathrm{T}} \overline{\mathbf{c}}_3^\mathrm{T} + \mathbf{Q}_{,2}^\mathrm{T} \left( \mathbf{v}, 0 \right) \overline{\mathbf{n}}_6^{e \mathrm{T}} \mathbf{s}^e \mathbf{s}^{e \mathrm{T}} \overline{\mathbf{c}}_4^\mathrm{T} \right] \phi \mathrm{d} \Gamma \\
& \qquad - \int_{\partial \Omega_d} \hat{\boldsymbol{\varepsilon}}^\mathrm{T} \left( \mathbf{v} \right) \overline{\mathbf{n}}_4^{e \mathrm{T}} \overline{\mathbf{n}}_4^e \overline{\mathbf{n}}_3^e \boldsymbol{\mu} \left( \mathbf{u}, \phi \right) \mathrm{d} \Gamma  - \int_{\partial \Omega_d} \boldsymbol{\mu}^\mathrm{T} \left( \mathbf{v}, 0 \right) \overline{\mathbf{n}}_3^{e \mathrm{T}} \overline{\mathbf{n}}_4^{e \mathrm{T}} \overline{\mathbf{n}}_4^e \hat{\boldsymbol{\varepsilon}} \left( \mathbf{u} \right)  \mathrm{d} \Gamma \\
& \qquad +  \int_{\partial \Omega_P} \mathbf{Q}^\mathrm{T} \left( \mathbf{v}, 0 \right) \overline{\mathbf{n}}_6^{e \mathrm{T}} \mathbf{n}^e \mathbf{n}^{e \mathrm{T}}  \mathbf{E} \left( \phi \right) \mathrm{d} \Gamma + \sum_{C_h} \left\{ \mathbf{Q}^\mathrm{T} \left( \mathbf{v}, 0 \right) \right\} \overline{\mathbf{n}}_6^{e \mathrm{T}} \mathbf{s}^e  \left[ \! \left[ \phi  \right] \! \right] \\
& \qquad - \sum_{C_h} \left[ \! \left[ \mathbf{v}^\mathrm{T} \right] \! \right] \overline{\mathbf{s}}_4^e \overline{\mathbf{n}}_3^e  \left\{ \boldsymbol{\mu} \left( \mathbf{u}, \phi \right)  \right\} - \sum_{C_h} \left\{ \boldsymbol{\mu}^\mathrm{T} \left( \mathbf{v}, 0 \right) \right\} \overline{\mathbf{n}}_3^{e \mathrm{T}} \overline{\mathbf{s}}_4^{e \mathrm{T}} \left[ \! \left[ \mathbf{u}  \right] \! \right] \\
& \qquad  = \int_\Omega \mathbf{v}^\mathrm{T} \mathbf{b} \mathrm{d} \Omega - \int_\Omega \hat{\boldsymbol{\varepsilon}}^\mathrm{T} \left( \mathbf{v} \right) \left[ \widehat{\mathbf{D}} \left( \mathbf{u}, \phi \right) \mathbf{E} \left( \phi \right) + \widehat{\mathbf{Q}} \left( \mathbf{u}, \phi \right) \mathbf{V} \left( \phi \right) \right] \mathrm{d} \Omega \\
& \qquad  + \int_{\partial \Omega_Q} \mathbf{v}^\mathrm{T} \left[ \begin{matrix} \overline{\mathbf{n}}_1^e \boldsymbol{\sigma} \left( \mathbf{u}, \phi \right) - \left( \overline{\mathbf{n}}_{21}^e +  \overline{\mathbf{c}}_1 \overline{\mathbf{s}}_4^{e \mathrm{T}} \overline{\mathbf{s}}_4^e \overline{\mathbf{n}}_3^e \right) \boldsymbol{\mu}_{,1} \left( \mathbf{u}, \phi \right) \\ - \left( \overline{\mathbf{n}}_{22}^e +\overline{ \mathbf{c}}_2 \overline{\mathbf{s}}_4^{e \mathrm{T}} \overline{\mathbf{s}}_4^e \overline{\mathbf{n}}_3^e \right) \boldsymbol{\mu}_{,2} \left( \mathbf{u}, \phi \right) + \overline{\mathbf{n}}_{4}^e \boldsymbol{\sigma}^{ES} \left( \mathbf{u}, \phi \right) \end{matrix} \right] \mathrm{d} \Gamma \\
& \qquad + \int_{\partial \Omega_R} \hat{\boldsymbol{\varepsilon}}^\mathrm{T} \left( \mathbf{v} \right) \overline{\mathbf{n}}_4^{e \mathrm{T}} \overline{\mathbf{n}}_4^e \overline{\mathbf{n}}_3^e \boldsymbol{\mu} \left( \mathbf{u}, \phi \right) \mathrm{d} \Gamma \\
& \qquad -\int_{\partial \Omega_u}  \left[ \begin{matrix} \boldsymbol{\sigma}^\mathrm{T} \left( \mathbf{v}, 0 \right) \overline{\mathbf{n}}_1^{e \mathrm{T}} - \boldsymbol{\mu}_{,1}^\mathrm{T} \left( \mathbf{v}, 0 \right) \left( \overline{\mathbf{n}}_{21}^{e \mathrm{T}} + \overline{\mathbf{n}}_3^{e \mathrm{T}} \overline{\mathbf{s}}_4^{e \mathrm{T}} \overline{\mathbf{s}}_4^e  \overline{\mathbf{c}}_1^\mathrm{T} \right)  \\ - \boldsymbol{\mu}_{,2}^\mathrm{T} \left( \mathbf{v}, 0 \right) \left( \overline{\mathbf{n}}_{22}^{e \mathrm{T}} + \overline{\mathbf{n}}_3^{e \mathrm{T}} \overline{\mathbf{s}}_4^{e \mathrm{T}} \overline{\mathbf{s}}_4^e  \overline{\mathbf{c}}_2^\mathrm{T} \right) \end{matrix} \right]  \mathbf{u} \mathrm{d} \Gamma \\
& \qquad  - \int_{\partial \Omega_\phi} \left[ \begin{matrix} \mathbf{D}^\mathrm{T} \left( \mathbf{v}, 0 \right) \mathbf{n}^{e} - \mathbf{Q}_{,1}^\mathrm{T} \left( \mathbf{v}, 0 \right) \left( \overline{\mathbf{n}}_{51}^{e \mathrm{T}} + \overline{\mathbf{n}}_6^{e \mathrm{T}} \mathbf{s}^e \mathbf{s}^{e \mathrm{T}} \overline{\mathbf{c}}_3^\mathrm{T} \right) \\ - \mathbf{Q}_{,2}^\mathrm{T} \left( \mathbf{v}, 0 \right) \left( \overline{\mathbf{n}}_{52}^{e \mathrm{T}} +  \overline{\mathbf{n}}_6^{e \mathrm{T}} \mathbf{s}^e \mathbf{s}^{e \mathrm{T}} \overline{\mathbf{c}}_4^\mathrm{T} \right) \end{matrix} \right]  \phi \mathrm{d} \Gamma \\
& \qquad - \int_{\partial \Omega_d} \boldsymbol{\mu}^\mathrm{T} \left( \mathbf{v}, 0 \right) \overline{\mathbf{n}}_3^{e \mathrm{T}} \overline{\mathbf{n}}_4^{e \mathrm{T}} \overline{\mathbf{n}}_4^e \hat{\boldsymbol{\varepsilon}} \left( \mathbf{u} \right)  \mathrm{d} \Gamma +  \int_{\partial \Omega_P} \mathbf{Q}^\mathrm{T} \left( \mathbf{v}, 0 \right) \overline{\mathbf{n}}_6^{e \mathrm{T}} \mathbf{n}^e \mathbf{n}^{e \mathrm{T}}  \mathbf{E} \left( \phi \right) \mathrm{d} \Gamma,
\label{eqn:dis_3}
\end{split}
\end{align}
\begin{align}
\begin{split}
& - \int_\Omega \mathbf{E}^\mathrm{T} \left( \tau \right) \mathbf{D} \left( \mathbf{u}, \phi \right) \mathrm{d} \Omega - \int_\Omega \mathbf{V}^\mathrm{T} \left( \tau \right) \mathbf{Q} \left( \mathbf{u}, \phi \right) \mathrm{d} \Omega \\
& \qquad -\int_{\Gamma_h \cup \partial \Omega_u}  \left\{ \begin{matrix} \boldsymbol{\sigma}^\mathrm{T} \left( \mathbf{0}, \tau \right) \overline{\mathbf{n}}_1^{e \mathrm{T}} - \boldsymbol{\mu}_{,1}^\mathrm{T} \left( \mathbf{0}, \tau \right) \overline{\mathbf{n}}_{21}^{e \mathrm{T}} - \boldsymbol{\mu}_{,2}^\mathrm{T} \left( \mathbf{0}, \tau \right) \overline{\mathbf{n}}_{22}^{e \mathrm{T}} \\ + \mathbf{E}^\mathrm{T} \left( \tau \right) \widehat{\mathbf{D}}^\mathrm{T} \left( \mathbf{u}, \phi \right) \overline{\mathbf{n}}_4^{e \mathrm{T}} + \mathbf{V}^\mathrm{T} \left( \tau \right) \widehat{\mathbf{Q}}^\mathrm{T} \left( \mathbf{u}, \phi \right) \overline{\mathbf{n}}_4^{e \mathrm{T}} \end{matrix} \right\} \left[ \! \left[ \mathbf{u} \right] \! \right] \mathrm{d} \Gamma \\
& \qquad  - \int_{\Gamma_h \cup \partial \Omega_\phi} \left[ \! \left[ \tau \right] \! \right] \left\{ \mathbf{n}^{e \mathrm{T}} \mathbf{D} \left( \mathbf{u}, \phi \right) - \overline{\mathbf{n}}_{51}^e  \mathbf{Q}_{,1} \left( \mathbf{u}, \phi \right) - \overline{\mathbf{n}}_{52}^e  \mathbf{Q}_{,2} \left( \mathbf{u}, \phi \right) \right\}  \mathrm{d} \Gamma \\
& \qquad  - \int_{\Gamma_h \cup \partial \Omega_\phi} \left\{ \mathbf{D}^\mathrm{T} \left( \mathbf{0}, \tau \right) \mathbf{n}^{e} - \mathbf{Q}_{,1}^\mathrm{T} \left( \mathbf{0}, \tau \right) \overline{\mathbf{n}}_{51}^{e \mathrm{T}}  - \mathbf{Q}_{,2}^\mathrm{T} \left( \mathbf{0}, \tau \right) \overline{\mathbf{n}}_{52}^{e \mathrm{T}} \right\} \left[ \! \left[ \phi \right] \! \right] \mathrm{d} \Gamma \\
& \qquad  + \int_{\Gamma_h} \left[ \! \left[ \mathbf{E}^\mathrm{T} \left( \tau \right) \right] \! \right] \overline{\mathbf{n}}_{6}^e \left\{ \mathbf{Q} \left( \mathbf{u}, \phi \right) \right\}  \mathrm{d} \Gamma + \int_{\Gamma_h} \left\{ \mathbf{Q}^\mathrm{T} \left( \mathbf{0}, \tau \right) \right\} \overline{\mathbf{n}}_{6}^{e \mathrm{T}} \left[ \! \left[ \mathbf{E} \left( \phi \right) \right] \! \right] \mathrm{d} \Gamma \\
& \qquad  - \int_{\Gamma_h} \left\{ \boldsymbol{\mu}^\mathrm{T} \left( \mathbf{0}, \tau \right) \right\}  \overline{\mathbf{n}}_3^{e \mathrm{T}} \left[ \! \left[ \hat{\boldsymbol{\varepsilon}} \left( \mathbf{u} \right) \right] \! \right]   \mathrm{d} \Gamma \\
& \qquad + \int_{\partial \Omega_u} \left[  \boldsymbol{\mu}_{,1}^\mathrm{T} \left( \mathbf{0}, \tau \right) \overline{\mathbf{n}}_3^{e \mathrm{T}} \overline{\mathbf{s}}_4^{e \mathrm{T}} \overline{\mathbf{s}}_4^e  \overline{\mathbf{c}}_1^\mathrm{T} + \boldsymbol{\mu}_{,2}^\mathrm{T} \left( \mathbf{0}, \tau \right) \overline{\mathbf{n}}_3^{e \mathrm{T}} \overline{\mathbf{s}}_4^{e \mathrm{T}} \overline{\mathbf{s}}_4^e  \overline{\mathbf{c}}_2^\mathrm{T} \right] \mathbf{u} \mathrm{d} \Gamma \\
& \qquad + \int_{\partial \Omega_\phi} \tau \left[ \overline{\mathbf{c}}_3 \mathbf{s}^e \mathbf{s}^{e \mathrm{T}} \overline{\mathbf{n}}_6^e \mathbf{Q}_{,1} \left( \mathbf{u}, \phi \right) + \overline{\mathbf{c}}_4 \mathbf{s}^e \mathbf{s}^{e \mathrm{T}} \overline{\mathbf{n}}_6^e \mathbf{Q}_{,2} \left( \mathbf{u}, \phi \right) \right] \mathrm{d} \Gamma \\
& \qquad + \int_{\partial \Omega_\phi} \left[ \mathbf{Q}_{,1}^\mathrm{T} \left( \mathbf{0}, \tau \right) \overline{\mathbf{n}}_6^{e \mathrm{T}} \mathbf{s}^e \mathbf{s}^{e \mathrm{T}} \overline{\mathbf{c}}_3^\mathrm{T} + \mathbf{Q}_{,2}^\mathrm{T} \left( \mathbf{0}, \tau \right) \overline{\mathbf{n}}_6^{e \mathrm{T}} \mathbf{s}^e \mathbf{s}^{e \mathrm{T}} \overline{\mathbf{c}}_4^\mathrm{T} \right] \phi \mathrm{d} \Gamma \\
& \qquad + \int_{\partial \Omega_P}  \mathbf{E}^\mathrm{T} \left( \tau \right) \mathbf{n}^e \mathbf{n}^{e \mathrm{T}} \overline{\mathbf{n}}_6^e \mathbf{Q} \left( \mathbf{u}, \phi \right) \mathrm{d} \Gamma +  \int_{\partial \Omega_P} \mathbf{Q}^\mathrm{T} \left( \mathbf{0}, \tau \right) \overline{\mathbf{n}}_6^{e \mathrm{T}} \mathbf{n}^e \mathbf{n}^{e \mathrm{T}}  \mathbf{E} \left( \phi \right) \mathrm{d} \Gamma \\
& \qquad - \int_{\partial \Omega_d} \boldsymbol{\mu}^\mathrm{T} \left( \mathbf{0}, \tau \right) \overline{\mathbf{n}}_3^{e \mathrm{T}} \overline{\mathbf{n}}_4^{e \mathrm{T}} \overline{\mathbf{n}}_4^e \hat{\boldsymbol{\varepsilon}} \left( \mathbf{u} \right)  \mathrm{d} \Gamma - \sum_{C_h} \left\{ \boldsymbol{\mu}^\mathrm{T} \left( \mathbf{0}, \tau \right) \right\} \overline{\mathbf{n}}_3^{e \mathrm{T}} \overline{\mathbf{s}}_4^{e \mathrm{T}} \left[ \! \left[ \mathbf{u}  \right] \! \right] \\
& \qquad + \sum_{C_h} \left[ \! \left[ \tau \right] \! \right] \mathbf{s}^{e \mathrm{T}} \overline{\mathbf{n}}_6^e  \left\{ \mathbf{Q} \left( \mathbf{u}, \phi \right)  \right\} + \sum_{C_h} \left\{ \mathbf{Q}^\mathrm{T} \left( \mathbf{0}, \tau \right) \right\} \overline{\mathbf{n}}_6^{e \mathrm{T}} \mathbf{s}^e  \left[ \! \left[ \phi  \right] \! \right] \\
& \qquad = - \int_\Omega \tau q \mathrm{d} \Omega + \int_{\partial \Omega_\omega} \tau \left[ \begin{matrix} \mathbf{n}^{e \mathrm{T}} \mathbf{D} \left( \mathbf{u}, \phi \right) - \left( \overline{\mathbf{n}}_{51}^e + \overline{\mathbf{c}}_3 \mathbf{s}^e \mathbf{s}^{e \mathrm{T}} \overline{\mathbf{n}}_6^e \right) \mathbf{Q}_{,1} \left( \mathbf{u}, \phi \right) \\ - \left( \overline{\mathbf{n}}_{52}^e +  \overline{\mathbf{c}}_4 \mathbf{s}^e \mathbf{s}^{e \mathrm{T}} \overline{\mathbf{n}}_6^e \right)  \mathbf{Q}_{,2} \left( \mathbf{u}, \phi \right) \end{matrix} \right]  \mathrm{d} \Gamma \\
& \qquad - \int_{\partial \Omega_Z}  \mathbf{E}^\mathrm{T} \left( \tau \right) \mathbf{n}^e \mathbf{n}^{e \mathrm{T}} \overline{\mathbf{n}}_6^e \mathbf{Q} \left( \mathbf{u}, \phi \right) \mathrm{d} \Gamma \\
& \qquad -\int_{\partial \Omega_u}  \left[ \begin{matrix} \boldsymbol{\sigma}^\mathrm{T} \left( \mathbf{0}, \tau \right) \overline{\mathbf{n}}_1^{e \mathrm{T}} - \boldsymbol{\mu}_{,1}^\mathrm{T} \left( \mathbf{0}, \tau \right) \left( \overline{\mathbf{n}}_{21}^{e \mathrm{T}} + \overline{\mathbf{n}}_3^{e \mathrm{T}} \overline{\mathbf{s}}_4^{e \mathrm{T}} \overline{\mathbf{s}}_4^e  \overline{\mathbf{c}}_1^\mathrm{T} \right) \\ - \boldsymbol{\mu}_{,2}^\mathrm{T} \left( \mathbf{0}, \tau \right) \left( \overline{\mathbf{n}}_{22}^{e \mathrm{T}} + \overline{\mathbf{n}}_3^{e \mathrm{T}} \overline{\mathbf{s}}_4^{e \mathrm{T}} \overline{\mathbf{s}}_4^e  \overline{\mathbf{c}}_2^\mathrm{T} \right) \\ + \mathbf{E}^\mathrm{T} \left( \tau \right) \widehat{\mathbf{D}}^\mathrm{T} \left( \mathbf{u}, \phi \right) \overline{\mathbf{n}}_4^{e \mathrm{T}} + \mathbf{V}^\mathrm{T} \left( \tau \right) \widehat{\mathbf{Q}}^\mathrm{T} \left( \mathbf{u}, \phi \right) \overline{\mathbf{n}}_4^{e \mathrm{T}} \end{matrix} \right] \mathbf{u} \mathrm{d} \Gamma \\
& \qquad  - \int_{\partial \Omega_\phi} \left[ \begin{matrix} \mathbf{D}^\mathrm{T} \left( \mathbf{0}, \tau \right) \mathbf{n}^{e} - \mathbf{Q}_{,1}^\mathrm{T} \left( \mathbf{0}, \tau \right) \left( \overline{\mathbf{n}}_{51}^{e \mathrm{T}} +  \overline{\mathbf{n}}_6^{e \mathrm{T}} \mathbf{s}^e \mathbf{s}^{e \mathrm{T}} \overline{\mathbf{c}}_3^\mathrm{T} \right)  \\ - \mathbf{Q}_{,2}^\mathrm{T} \left( \mathbf{0}, \tau \right) \left( \overline{\mathbf{n}}_{52}^{e \mathrm{T}} + \overline{\mathbf{n}}_6^{e \mathrm{T}} \mathbf{s}^e \mathbf{s}^{e \mathrm{T}} \overline{\mathbf{c}}_4^\mathrm{T} \right) \end{matrix} \right] \phi \mathrm{d} \Gamma \\
& \qquad - \int_{\partial \Omega_d} \boldsymbol{\mu}^\mathrm{T} \left( \mathbf{0}, \tau \right) \overline{\mathbf{n}}_3^{e \mathrm{T}} \overline{\mathbf{n}}_4^{e \mathrm{T}} \overline{\mathbf{n}}_4^e \hat{\boldsymbol{\varepsilon}} \left( \mathbf{u} \right)  \mathrm{d} \Gamma \\
& \qquad +  \int_{\partial \Omega_P} \mathbf{Q}^\mathrm{T} \left( \mathbf{0}, \tau \right) \overline{\mathbf{n}}_6^{e \mathrm{T}} \mathbf{n}^e \mathbf{n}^{e \mathrm{T}}  \mathbf{E} \left( \phi \right) \mathrm{d} \Gamma.
\label{eqn:ele_3}
\end{split}
\end{align}

Furthermore, in order to impose the essential boundary conditions ($\partial \Omega_u$, $\partial \Omega_d$, $\partial \Omega_\phi$ and $\partial \Omega_P$) and the continuity condition across subdomains, we employ the Numerical Flux Correction, which is widely used in Discontinuous Galerkin (DG) Methods \cite{Bassi1997, Mozolevski2007} to ensure the consistency and stability of the method. Here we use the Interior Penalty (IP) Numerical Flux Corrections. A set of penalty parameters are introduced, in which $\boldsymbol{\eta}_1 = \left[ \eta_{11}, \eta_{12}, \eta_{13}, \eta_{14} \right]$ and $\boldsymbol{\eta}_2 = \left[ \eta_{21}, \eta_{22}, \eta_{23}, \eta_{24} \right]$ are related to the essential boundary conditions and the interior continuity conditions respectively. The FPM is only stable when the penalty parameters are large enough.

Finally, after substituting the constitutive equations and the following boundary conditions (see Eqn.~\ref{eqn:BC_full_1} -- \ref{eqn:BC_full_8}):
\begin{align}
\mathbf{u} & = \widetilde{\mathbf{u}}, & \text{on} \; \partial \Omega_u, \\
\phi & = \widetilde{\phi}, & \text{on} \; \partial \Omega_{\phi}, \\
\overline{\mathbf{n}}_4^e \hat{\boldsymbol{\varepsilon}} \left( \mathbf{u} \right) & = \widetilde{\mathbf{d}}, & \text{on} \; \partial \Omega_{d}, \\
\mathbf{n}^{e \mathrm{T}}  \mathbf{E} \left( \phi \right) & = - \widetilde{P}, & \text{on} \; \partial \Omega_{P}, \\
\begin{matrix} \overline{\mathbf{n}}_1^e \boldsymbol{\sigma} \left( \mathbf{u}, \phi \right) - \left( \overline{\mathbf{n}}_{21}^e +  \overline{\mathbf{c}}_1 \overline{\mathbf{s}}_4^{e \mathrm{T}} \overline{\mathbf{s}}_4^e \overline{\mathbf{n}}_3^e \right) \boldsymbol{\mu}_{,1} \left( \mathbf{u}, \phi \right) \\ \qquad \qquad - \left( \overline{\mathbf{n}}_{22}^e +\overline{ \mathbf{c}}_2 \overline{\mathbf{s}}_4^{e \mathrm{T}} \overline{\mathbf{s}}_4^e \overline{\mathbf{n}}_3^e \right) \boldsymbol{\mu}_{,2} \left( \mathbf{u}, \phi \right) + \overline{\mathbf{n}}_{4}^e \boldsymbol{\sigma}^{ES} \left( \mathbf{u}, \phi \right) \end{matrix} & = \widetilde{\mathbf{Q}}, & \text{on} \; \partial \Omega_Q, \\
\begin{matrix} \mathbf{n}^{e \mathrm{T}} \mathbf{D} \left( \mathbf{u}, \phi \right) - \left( \overline{\mathbf{n}}_{51}^e + \overline{\mathbf{c}}_3 \mathbf{s}^e \mathbf{s}^{e \mathrm{T}} \overline{\mathbf{n}}_6^e \right) \mathbf{Q}_{,1} \left( \mathbf{u}, \phi \right) \qquad \qquad \\ - \left( \overline{\mathbf{n}}_{52}^e +  \overline{\mathbf{c}}_4 \mathbf{s}^e \mathbf{s}^{e \mathrm{T}} \overline{\mathbf{n}}_6^e \right)  \mathbf{Q}_{,2} \left( \mathbf{u}, \phi \right) \end{matrix} & = - \widetilde{\omega}, & \text{on} \; \partial \Omega_{\omega}, \\
\overline{\mathbf{n}}_4^e \overline{\mathbf{n}}_3^e \boldsymbol{\mu} \left( \mathbf{u}, \phi \right) & = \widetilde{\mathbf{R}}, & \text{on} \; \partial \Omega_{R}, \\
\mathbf{n}^{e \mathrm{T}} \overline{\mathbf{n}}_6^e \mathbf{Q} \left( \mathbf{u}, \phi \right) & = \widetilde{Z}, & \text{on} \; \partial \Omega_{Z},
\end{align}
we can obtain the final symmetric formulation for the primal FPM for full flexoelectric theory:
\begin{align}
\begin{split}
\label{eqn:dis_full}
& \int_\Omega \boldsymbol{\varepsilon}^\mathrm{T} \left( \mathbf{v} \right) \mathbf{C}_{\sigma \varepsilon} \boldsymbol{\varepsilon} \left( \mathbf{u} \right) \mathrm{d} \Omega +  \int_\Omega \boldsymbol{\kappa}^\mathrm{T} \left( \mathbf{v} \right) \mathbf{C}_{\mu \kappa} \boldsymbol{\kappa} \left( \mathbf{u} \right) \mathrm{d} \Omega \\
& \qquad - \int_\Omega \boldsymbol{\varepsilon}^\mathrm{T} \left( \mathbf{v} \right) \left( \mathbf{e} \mathbf{E} \left( \phi \right) + \mathbf{b} \mathbf{V} \left( \phi \right) \right) \mathrm{d} \Omega -  \int_\Omega \boldsymbol{\kappa}^\mathrm{T} \left( \mathbf{v} \right) \mathbf{a} \mathbf{E} \left( \phi \right) \mathrm{d} \Omega \\
& \qquad   -\int_{\Gamma_h \cup \partial \Omega_u} \left( \left[ \! \left[ \mathbf{v}^\mathrm{T} \right] \! \right] \overline{\mathbf{n}}_1^e \mathbf{C}_{\sigma \varepsilon} \left\{ \boldsymbol{\varepsilon} \left( \mathbf{u} \right) \right\} + \left\{ \boldsymbol{\varepsilon}^\mathrm{T} \left( \mathbf{v} \right) \right\} \mathbf{C}_{\sigma \varepsilon} \overline{\mathbf{n}}_1^{e \mathrm{T}} \left[ \! \left[ \mathbf{u} \right] \! \right]  \right) \mathrm{d} \Gamma \\
& \qquad   + \int_{\Gamma_h \cup \partial \Omega_u} \left( \left[ \! \left[ \mathbf{v}^\mathrm{T} \right] \! \right] \overline{\mathbf{n}}_{21}^e \mathbf{C}_{\mu \kappa} \left\{ \boldsymbol{\kappa}_{,1} \left( \mathbf{u} \right) \right\} + \left\{ \boldsymbol{\kappa}_{,1}^\mathrm{T} \left( \mathbf{v} \right) \right\} \mathbf{C}_{\mu \kappa} \overline{\mathbf{n}}_{21}^{e \mathrm{T}} \left[ \! \left[ \mathbf{u} \right] \! \right]  \right) \mathrm{d} \Gamma \\
& \qquad   + \int_{\Gamma_h \cup \partial \Omega_u} \left( \left[ \! \left[ \mathbf{v}^\mathrm{T} \right] \! \right] \overline{\mathbf{n}}_{22}^e \mathbf{C}_{\mu \kappa} \left\{ \boldsymbol{\kappa}_{,2} \left( \mathbf{u} \right) \right\} + \left\{ \boldsymbol{\kappa}_{,2}^\mathrm{T} \left( \mathbf{v} \right) \right\} \mathbf{C}_{\mu \kappa} \overline{\mathbf{n}}_{22}^{e \mathrm{T}} \left[ \! \left[ \mathbf{u} \right] \! \right]  \right) \mathrm{d} \Gamma \\
& \qquad   + \int_{\Gamma_h \cup \partial \Omega_u} \left[ \! \left[ \mathbf{v}^\mathrm{T} \right] \! \right] \left\{ \left( \overline{\mathbf{n}}_1^e \mathbf{b} - \overline{\mathbf{n}}_{21}^e \mathbf{a} \overline{\mathbf{c}}_5 - \overline{\mathbf{n}}_{22}^e \mathbf{a} \overline{\mathbf{c}}_6 - \overline{\mathbf{n}}_4^e  \widehat{\mathbf{Q}} \left( \mathbf{u}, \phi \right) \right) \mathbf{V} \left( \phi \right) \right\} \mathrm{d} \Gamma \\
& \qquad   + \int_{\Gamma_h \cup \partial \Omega_u} \left[ \! \left[ \mathbf{v}^\mathrm{T} \right] \! \right] \left\{ \left( \overline{\mathbf{n}}_1^e \mathbf{e} - \overline{\mathbf{n}}_4^e  \widehat{\mathbf{D}} \left( \mathbf{u}, \phi \right) \right) \mathbf{E} \left( \phi \right) \right\} \mathrm{d} \Gamma  - \int_{\Gamma_h \cup \partial \Omega_\phi} \left\{ \boldsymbol{\varepsilon}^\mathrm{T} \left( \mathbf{v} \right) \right\}  \mathbf{e} \mathbf{n}^{e}   \left[ \! \left[ \phi \right] \! \right] \mathrm{d} \Gamma \\
& \qquad  - \int_{\Gamma_h \cup \partial \Omega_\phi} \left\{ \boldsymbol{\kappa}^\mathrm{T} \left( \mathbf{v} \right) \right\}  \left( \mathbf{a} \mathbf{n}^{e} -  \overline{\mathbf{c}}_7^\mathrm{T} \mathbf{b} \overline{\mathbf{n}}_{51}^{e \mathrm{T}} -  \overline{\mathbf{c}}_8^\mathrm{T} \mathbf{b} \overline{\mathbf{n}}_{52}^{e \mathrm{T}}  \right)  \left[ \! \left[ \phi \right] \! \right] \mathrm{d} \Gamma \\
& \qquad -\int_{\Gamma_h} \left( \left[ \! \left[ \hat{\boldsymbol{\varepsilon}}^\mathrm{T} \left( \mathbf{v} \right) \right] \! \right] \overline{\mathbf{n}}_3^e \mathbf{C}_{\mu \kappa} \left\{ \boldsymbol{\kappa} \left( \mathbf{u} \right) \right\} + \left\{ \boldsymbol{\kappa}^\mathrm{T} \left( \mathbf{v} \right) \right\} \mathbf{C}_{\mu \kappa} \overline{\mathbf{n}}_3^{e \mathrm{T}} \left[ \! \left[ \hat{\boldsymbol{\varepsilon}} \left( \mathbf{u} \right) \right] \! \right] \right)  \mathrm{d} \Gamma \\
& \qquad +\int_{\Gamma_h} \left[ \! \left[ \hat{\boldsymbol{\varepsilon}}^\mathrm{T} \left( \mathbf{v} \right) \right] \! \right] \overline{\mathbf{n}}_3^e \mathbf{a} \left\{ \mathbf{E} \left( \phi \right) \right\}  \mathrm{d} \Gamma + \int_{\Gamma_h} \left\{ \boldsymbol{\varepsilon}^\mathrm{T} \left( \mathbf{v} \right) \right\} \mathbf{b} \overline{\mathbf{n}}_{6}^{e \mathrm{T}} \left[ \! \left[ \mathbf{E} \left( \phi \right) \right] \! \right] \mathrm{d} \Gamma \\
& \qquad + \int_{\partial \Omega_u} \left( \mathbf{v}^\mathrm{T} \overline{\mathbf{c}}_1 \overline{\mathbf{s}}_4^{e \mathrm{T}} \overline{\mathbf{s}}_4^e \overline{\mathbf{n}}_3^e \mathbf{C}_{\mu \kappa} \boldsymbol{\kappa}_{,1} \left( \mathbf{u} \right) + \boldsymbol{\kappa}_{,1}^\mathrm{T} \left( \mathbf{v} \right) \mathbf{C}_{\mu \kappa} \overline{\mathbf{n}}_3^{e \mathrm{T}} \overline{\mathbf{s}}_4^{e \mathrm{T}} \overline{\mathbf{s}}_4^e  \overline{\mathbf{c}}_1^\mathrm{T} \mathbf{u}  \right)  \mathrm{d} \Gamma \\
& \qquad + \int_{\partial \Omega_u} \left( \mathbf{v}^\mathrm{T} \overline{\mathbf{c}}_2 \overline{\mathbf{s}}_4^{e \mathrm{T}} \overline{\mathbf{s}}_4^e \overline{\mathbf{n}}_3^e \mathbf{C}_{\mu \kappa} \boldsymbol{\kappa}_{,2} \left( \mathbf{u} \right) + \boldsymbol{\kappa}_{,2}^\mathrm{T} \left( \mathbf{v} \right) \mathbf{C}_{\mu \kappa} \overline{\mathbf{n}}_3^{e \mathrm{T}} \overline{\mathbf{s}}_4^{e \mathrm{T}} \overline{\mathbf{s}}_4^e  \overline{\mathbf{c}}_2^\mathrm{T} \mathbf{u}  \right)  \mathrm{d} \Gamma \\
& \qquad - \int_{\partial \Omega_u} \mathbf{v}^\mathrm{T} \left( \overline{\mathbf{c}}_1 \overline{\mathbf{s}}_4^{e \mathrm{T}} \overline{\mathbf{s}}_4^e \overline{\mathbf{n}}_3^e \mathbf{a} \overline{\mathbf{c}}_5 + \overline{\mathbf{c}}_2 \overline{\mathbf{s}}_4^{e \mathrm{T}} \overline{\mathbf{s}}_4^e \overline{\mathbf{n}}_3^e \mathbf{a} \overline{\mathbf{c}}_6 \right) \mathbf{V} \left( \phi \right)  \mathrm{d} \Gamma \\
& \qquad + \int_{\partial \Omega_\phi} \boldsymbol{\kappa}^\mathrm{T} \left( \mathbf{v} \right)  \left( \overline{\mathbf{c}}_7^\mathrm{T} \mathbf{b} \overline{\mathbf{n}}_6^{e \mathrm{T}} \mathbf{s}^e \mathbf{s}^{e \mathrm{T}} \overline{\mathbf{c}}_3^\mathrm{T} + \overline{\mathbf{c}}_8^\mathrm{T} \mathbf{b} \overline{\mathbf{n}}_6^{e \mathrm{T}} \mathbf{s}^e \mathbf{s}^{e \mathrm{T}} \overline{\mathbf{c}}_4^\mathrm{T} \right) \phi \mathrm{d} \Gamma \\
& \qquad - \int_{\partial \Omega_d} \left( \hat{\boldsymbol{\varepsilon}}^\mathrm{T} \left( \mathbf{v} \right) \overline{\mathbf{n}}_4^{e \mathrm{T}} \overline{\mathbf{n}}_4^e \overline{\mathbf{n}}_3^e \mathbf{C}_{\mu \kappa} \boldsymbol{\kappa} \left( \mathbf{u} \right) + \boldsymbol{\kappa}^\mathrm{T} \left( \mathbf{v} \right) \mathbf{C}_{\mu \kappa} \overline{\mathbf{n}}_3^{e \mathrm{T}} \overline{\mathbf{n}}_4^{e \mathrm{T}} \overline{\mathbf{n}}_4^e \hat{\boldsymbol{\varepsilon}} \left( \mathbf{u} \right) \right)  \mathrm{d} \Gamma \\
& \qquad + \int_{\partial \Omega_d} \hat{\boldsymbol{\varepsilon}}^\mathrm{T} \left( \mathbf{v} \right) \overline{\mathbf{n}}_4^{e \mathrm{T}} \overline{\mathbf{n}}_4^e \overline{\mathbf{n}}_3^e \mathbf{a} \mathbf{E} \left( \phi \right) \mathrm{d} \Gamma  +  \int_{\partial \Omega_P}  \boldsymbol{\varepsilon}^\mathrm{T} \left( \mathbf{v} \right) \mathbf{b} \overline{\mathbf{n}}_6^{e \mathrm{T}} \mathbf{n}^e \mathbf{n}^{e \mathrm{T}}  \mathbf{E} \left( \phi \right) \mathrm{d} \Gamma \\
& \qquad + \int_{\partial \Omega_u} \frac{\eta_{11}}{h_e} \mathbf{v}^\mathrm{T} \mathbf{u} \mathrm{d} \Gamma + \int_{\Gamma_h} \frac{\eta_{21}}{h_e} \left[ \! \left[ \mathbf{v}^\mathrm{T} \right] \! \right] \left[ \!\left[ \mathbf{u} \right] \! \right] \mathrm{d} \Gamma \\
& \qquad + \int_{\partial \Omega_d} \eta_{12} h_e \hat{\boldsymbol{\varepsilon}}^\mathrm{T} \left( \mathbf{v} \right) \overline{\mathbf{n}}_4^{e \mathrm{T}} \overline{\mathbf{n}}_4^e \hat{\boldsymbol{\varepsilon}} \left( \mathbf{u} \right) \mathrm{d} \Gamma + \int_{\Gamma_h} \eta_{22} h_e \left[ \! \left[ \hat{\boldsymbol{\varepsilon}}^\mathrm{T} \left( \mathbf{v} \right) \right] \! \right] \overline{\mathbf{n}}_4^{e \mathrm{T}} \overline{\mathbf{n}}_4^e \left[ \! \left[ \hat{\boldsymbol{\varepsilon}} \left( \mathbf{u} \right) \right] \! \right] \mathrm{d} \Gamma 
\end{split} \\
\begin{split} \nonumber
& \qquad - \sum_{C_h} \left[ \! \left[ \mathbf{v}^\mathrm{T} \right] \! \right] \overline{\mathbf{s}}_4^e \overline{\mathbf{n}}_3^e  \mathbf{C}_{\mu \kappa} \left\{ \boldsymbol{\kappa} \left( \mathbf{u} \right) \right\}  - \sum_{C_h} \left\{ \boldsymbol{\kappa}^\mathrm{T} \left( \mathbf{v} \right) \right\} \mathbf{C}_{\mu \kappa} \overline{\mathbf{n}}_3^{e \mathrm{T}} \overline{\mathbf{s}}_4^{e \mathrm{T}} \left[ \! \left[ \mathbf{u}  \right] \! \right] \\
& \qquad + \sum_{C_h} \left\{ \boldsymbol{\varepsilon}^\mathrm{T} \left( \mathbf{v} \right)  \right\} \mathbf{b} \overline{\mathbf{n}}_6^{e \mathrm{T}} \mathbf{s}^e  \left[ \! \left[ \phi  \right] \! \right] +  \sum_{C_h} \left[ \! \left[ \mathbf{v}^\mathrm{T} \right] \! \right] \overline{\mathbf{s}}_4^e \overline{\mathbf{n}}_3^e \mathbf{a} \left\{ \mathbf{E} \left( \phi \right)   \right\} \\
& \qquad  = \int_\Omega \mathbf{v}^\mathrm{T} \mathbf{b} \mathrm{d} \Omega - \int_\Omega \hat{\boldsymbol{\varepsilon}}^\mathrm{T} \left( \mathbf{v} \right) \left( \widehat{\mathbf{D}} \left( \mathbf{u}, \phi \right) \mathbf{E} \left( \phi \right) + \widehat{\mathbf{Q}} \left( \mathbf{u}, \phi \right) \mathbf{V} \left( \phi \right) \right) \mathrm{d} \Omega \\
& \qquad  + \int_{\partial \Omega_Q} \mathbf{v}^\mathrm{T} \widetilde{\mathbf{Q}} \mathrm{d} \Gamma + \int_{\partial \Omega_R} \hat{\boldsymbol{\varepsilon}}^\mathrm{T} \left( \mathbf{v} \right) \overline{\mathbf{n}}_4^{e \mathrm{T}} \widetilde{\mathbf{R}} \mathrm{d} \Gamma \\
& \qquad -\int_{\partial \Omega_u}  \left[ \begin{matrix} \boldsymbol{\varepsilon}^\mathrm{T} \left( \mathbf{v} \right) \mathbf{C}_{\sigma \varepsilon} \overline{\mathbf{n}}_1^{e \mathrm{T}} - \boldsymbol{\kappa}_{,1}^\mathrm{T} \left( \mathbf{v} \right) \mathbf{C}_{\mu \kappa} \left( \overline{\mathbf{n}}_{21}^{e \mathrm{T}} + \overline{\mathbf{n}}_3^{e \mathrm{T}} \overline{\mathbf{s}}_4^{e \mathrm{T}} \overline{\mathbf{s}}_4^e  \overline{\mathbf{c}}_1^\mathrm{T} \right)  \\ - \boldsymbol{\kappa}_{,2}^\mathrm{T} \left( \mathbf{v} \right) \mathbf{C}_{\mu \kappa} \left( \overline{\mathbf{n}}_{22}^{e \mathrm{T}} + \overline{\mathbf{n}}_3^{e \mathrm{T}} \overline{\mathbf{s}}_4^{e \mathrm{T}} \overline{\mathbf{s}}_4^e  \overline{\mathbf{c}}_2^\mathrm{T} \right) \end{matrix} \right]  \widetilde{\mathbf{u}} \mathrm{d} \Gamma \\
& \qquad  - \int_{\partial \Omega_\phi} \left[ \begin{matrix} \boldsymbol{\varepsilon}^\mathrm{T} \left( \mathbf{v} \right) \mathbf{e} \mathbf{n}^{e} \\ + \boldsymbol{\kappa}^\mathrm{T} \left( \mathbf{v} \right) \left[ \mathbf{a} \mathbf{n}^{e} -  \overline{\mathbf{c}}_7^\mathrm{T} \mathbf{b} \left( \overline{\mathbf{n}}_{51}^{e \mathrm{T}} + \overline{\mathbf{n}}_6^{e \mathrm{T}} \mathbf{s}^e \mathbf{s}^{e \mathrm{T}} \overline{\mathbf{c}}_3^\mathrm{T} \right) -  \overline{\mathbf{c}}_8^\mathrm{T} \mathbf{b} \left( \overline{\mathbf{n}}_{52}^{e \mathrm{T}} +  \overline{\mathbf{n}}_6^{e \mathrm{T}} \mathbf{s}^e \mathbf{s}^{e \mathrm{T}} \overline{\mathbf{c}}_4^\mathrm{T} \right)  \right] \end{matrix} \right]  \widetilde{\phi} \mathrm{d} \Gamma \\
& \qquad - \int_{\partial \Omega_d} \boldsymbol{\kappa}^\mathrm{T} \left( \mathbf{v} \right) \mathbf{C}_{\mu \kappa} \overline{\mathbf{n}}_3^{e \mathrm{T}} \overline{\mathbf{n}}_4^{e \mathrm{T}} \widetilde{\mathbf{d}}  \mathrm{d} \Gamma -  \int_{\partial \Omega_P} \boldsymbol{\varepsilon}^\mathrm{T} \left( \mathbf{v} \right) \mathbf{b} \overline{\mathbf{n}}_6^{e \mathrm{T}} \mathbf{n}^e \widetilde{P} \mathrm{d} \Gamma \\
& \qquad + \int_{\partial \Omega_u} \frac{\eta_{11}}{h_e} \mathbf{v}^\mathrm{T} \widetilde{\mathbf{u}} \mathrm{d} \Gamma + \int_{\partial \Omega_d} \eta_{12} h_e \hat{\boldsymbol{\varepsilon}}^\mathrm{T} \left( \mathbf{v} \right) \overline{\mathbf{n}}_4^{e \mathrm{T}} \widetilde{\mathbf{d}} \mathrm{d} \Gamma,
\end{split}  \\~\nonumber \\
\begin{split}
\label{eqn:ele_full}
& - \int_\Omega \mathbf{E}^\mathrm{T} \left( \tau \right) \boldsymbol{\Lambda} \mathbf{E} \left( \phi \right) \mathrm{d} \Omega - \int_\Omega \mathbf{V}^\mathrm{T} \left( \tau \right)  \boldsymbol{\Phi} \mathbf{V} \left( \phi \right) \mathrm{d} \Omega \\
& \qquad - \int_\Omega \left( \mathbf{E}^\mathrm{T} \left( \tau \right) \mathbf{e}^\mathrm{T} + \mathbf{V}^\mathrm{T} \left( \tau \right) \mathbf{b}^\mathrm{T} \right) \boldsymbol{\varepsilon} \left( \mathbf{u} \right) \mathrm{d} \Omega - \int_\Omega \mathbf{E}^\mathrm{T} \left( \tau \right)  \mathbf{a}^\mathrm{T}  \boldsymbol{\kappa} \left( \mathbf{u} \right) \mathrm{d} \Omega \\
& \qquad + \int_{\Gamma_h \cup \partial \Omega_u}  \left\{ \mathbf{E}^\mathrm{T} \left( \tau \right) \left( \mathbf{e}^\mathrm{T} \overline{\mathbf{n}}_1^{e \mathrm{T}} - \widehat{\mathbf{D}}^\mathrm{T} \left( \mathbf{u}, \phi \right) \overline{\mathbf{n}}_4^{e \mathrm{T}} \right) \right\} \left[ \! \left[ \mathbf{u} \right] \! \right] \mathrm{d} \Gamma \\
& \qquad + \int_{\Gamma_h \cup \partial \Omega_u}  \left\{ \mathbf{V}^\mathrm{T} \left( \tau \right) \left( \mathbf{b}^\mathrm{T} \overline{\mathbf{n}}_1^{e \mathrm{T}} - \overline{\mathbf{c}}_5^\mathrm{T} \mathbf{a}^\mathrm{T} \overline{\mathbf{n}}_{21}^{e \mathrm{T}} - \overline{\mathbf{c}}_6^\mathrm{T} \mathbf{a}^\mathrm{T} \overline{\mathbf{n}}_{22}^{e \mathrm{T}} - \widehat{\mathbf{Q}}^\mathrm{T} \left( \mathbf{u}, \phi \right) \overline{\mathbf{n}}_4^{e \mathrm{T}} \right) \right\} \left[ \! \left[ \mathbf{u} \right] \! \right] \mathrm{d} \Gamma \\
& \qquad  - \int_{\Gamma_h \cup \partial \Omega_\phi} \left(  \left[ \! \left[ \tau \right] \! \right] \mathbf{n}^{e \mathrm{T}} \boldsymbol{\Lambda} \left\{ \mathbf{E} \left( \phi \right) \right\} + \left\{ \mathbf{E}^\mathrm{T} \left( \tau \right) \right\} \boldsymbol{\Lambda}  \mathbf{n}^{e} \left[ \! \left[ \phi \right] \! \right] \right) \mathrm{d} \Gamma \\
& \qquad  + \int_{\Gamma_h \cup \partial \Omega_\phi} \left(  \left[ \! \left[ \tau \right] \! \right] \overline{\mathbf{n}}_{51}^e \boldsymbol{\Phi} \left\{ \mathbf{V}_{,1} \left( \phi \right) \right\} + \left\{ \mathbf{V}_{,1}^\mathrm{T} \left( \tau \right) \right\} \boldsymbol{\Phi} \overline{\mathbf{n}}_{51}^{e \mathrm{T}} \left[ \! \left[ \phi \right] \! \right] \right) \mathrm{d} \Gamma \\
& \qquad  + \int_{\Gamma_h \cup \partial \Omega_\phi} \left(  \left[ \! \left[ \tau \right] \! \right] \overline{\mathbf{n}}_{52}^e \boldsymbol{\Phi} \left\{ \mathbf{V}_{,2} \left( \phi \right) \right\} + \left\{ \mathbf{V}_{,2}^\mathrm{T} \left( \tau \right) \right\} \boldsymbol{\Phi} \overline{\mathbf{n}}_{52}^{e \mathrm{T}} \left[ \! \left[ \phi \right] \! \right] \right) \mathrm{d} \Gamma \\
& \qquad  - \int_{\Gamma_h \cup \partial \Omega_\phi} \left[ \! \left[ \tau \right] \! \right] \mathbf{n}^{e \mathrm{T}} \mathbf{e}^\mathrm{T} \left\{ \boldsymbol{\varepsilon} \left( \mathbf{u} \right) \right\} \mathrm{d} \Gamma \\
& \qquad  - \int_{\Gamma_h \cup \partial \Omega_\phi} \left[ \! \left[ \tau \right] \! \right] \left( \mathbf{n}^{e \mathrm{T}} \mathbf{a}^\mathrm{T} - \overline{\mathbf{n}}_{51}^e \mathbf{b}^\mathrm{T} \overline{\mathbf{c}}_7 - \overline{\mathbf{n}}_{52}^e \mathbf{b}^\mathrm{T} \overline{\mathbf{c}}_8 \right) \left\{ \boldsymbol{\kappa} \left( \mathbf{u} \right) \right\} \mathrm{d} \Gamma \\
& \qquad  + \int_{\Gamma_h} \left( \left[ \! \left[ \mathbf{E}^\mathrm{T} \left( \tau \right) \right] \! \right] \overline{\mathbf{n}}_{6}^e \boldsymbol{\Phi} \left\{ \mathbf{V} \left( \phi \right) \right\} + \left\{ \mathbf{V}^\mathrm{T} \left( \tau \right) \right\} \boldsymbol{\Phi} \overline{\mathbf{n}}_{6}^{e \mathrm{T}} \left[ \! \left[ \mathbf{E} \left( \phi \right) \right] \! \right] \right) \mathrm{d} \Gamma \\
& \qquad  + \int_{\Gamma_h} \left[ \! \left[ \mathbf{E}^\mathrm{T} \left( \tau \right) \right] \! \right] \overline{\mathbf{n}}_{6}^e \mathbf{b}^\mathrm{T} \left\{\boldsymbol{\varepsilon} \left( \mathbf{u} \right) \right\}  \mathrm{d} \Gamma  + \int_{\Gamma_h} \left\{ \mathbf{E}^\mathrm{T} \left( \tau \right) \right\} \mathbf{a}^\mathrm{T} \overline{\mathbf{n}}_3^{e \mathrm{T}} \left[ \! \left[ \hat{\boldsymbol{\varepsilon}} \left( \mathbf{u} \right) \right] \! \right]   \mathrm{d} \Gamma \\
& \qquad - \int_{\partial \Omega_u} \mathbf{V}^\mathrm{T} \left( \tau \right) \left( \overline{\mathbf{c}}_5^\mathrm{T} \mathbf{a}^\mathrm{T} \overline{\mathbf{n}}_3^{e \mathrm{T}} \overline{\mathbf{s}}_4^{e \mathrm{T}} \overline{\mathbf{s}}_4^e  \overline{\mathbf{c}}_1^\mathrm{T} +  \overline{\mathbf{c}}_6^\mathrm{T} \mathbf{a}^\mathrm{T} \overline{\mathbf{n}}_3^{e \mathrm{T}} \overline{\mathbf{s}}_4^{e \mathrm{T}} \overline{\mathbf{s}}_4^e  \overline{\mathbf{c}}_2^\mathrm{T} \right) \mathbf{u} \mathrm{d} \Gamma \\
& \qquad + \int_{\partial \Omega_\phi} \left( \tau \overline{\mathbf{c}}_3 \mathbf{s}^e \mathbf{s}^{e \mathrm{T}} \overline{\mathbf{n}}_6^e \boldsymbol{\Phi} \mathbf{V}_{,1} \left( \phi \right)  + \mathbf{V}_{,1}^\mathrm{T} \left( \tau \right) \boldsymbol{\Phi} \overline{\mathbf{n}}_6^{e \mathrm{T}} \mathbf{s}^e \mathbf{s}^{e \mathrm{T}} \overline{\mathbf{c}}_3^\mathrm{T} \phi \right) \mathrm{d} \Gamma \\
& \qquad + \int_{\partial \Omega_\phi} \left( \tau \overline{\mathbf{c}}_4 \mathbf{s}^e \mathbf{s}^{e \mathrm{T}} \overline{\mathbf{n}}_6^e \boldsymbol{\Phi} \mathbf{V}_{,2} \left( \phi \right)  + \mathbf{V}_{,2}^\mathrm{T} \left( \tau \right) \boldsymbol{\Phi} \overline{\mathbf{n}}_6^{e \mathrm{T}} \mathbf{s}^e \mathbf{s}^{e \mathrm{T}} \overline{\mathbf{c}}_4^\mathrm{T} \phi \right) \mathrm{d} \Gamma \\
& \qquad + \int_{\partial \Omega_\phi} \tau \left( \overline{\mathbf{c}}_3 \mathbf{s}^e \mathbf{s}^{e \mathrm{T}} \overline{\mathbf{n}}_6^e \mathbf{b}^\mathrm{T} \overline{\mathbf{c}}_7 + \overline{\mathbf{c}}_4 \mathbf{s}^e \mathbf{s}^{e \mathrm{T}} \overline{\mathbf{n}}_6^e \mathbf{b}^\mathrm{T} \overline{\mathbf{c}}_8 \right)  \boldsymbol{\kappa} \left( \mathbf{u} \right) \mathrm{d} \Gamma \\
& \qquad + \int_{\partial \Omega_d} \mathbf{E}^\mathrm{T} \left( \tau \right) \mathbf{a}^\mathrm{T} \overline{\mathbf{n}}_3^{e \mathrm{T}} \overline{\mathbf{n}}_4^{e \mathrm{T}} \overline{\mathbf{n}}_4^e \hat{\boldsymbol{\varepsilon}} \left( \mathbf{u} \right)  \mathrm{d} \Gamma + \int_{\partial \Omega_P}  \mathbf{E}^\mathrm{T} \left( \tau \right) \mathbf{n}^e \mathbf{n}^{e \mathrm{T}} \overline{\mathbf{n}}_6^e \mathbf{b}^\mathrm{T} \boldsymbol{\varepsilon} \left( \mathbf{u} \right) \mathrm{d} \Gamma \\
& \qquad + \int_{\partial \Omega_P}  \left( \mathbf{E}^\mathrm{T} \left( \tau \right) \mathbf{n}^e \mathbf{n}^{e \mathrm{T}} \overline{\mathbf{n}}_6^e \boldsymbol{\Phi} \mathbf{V} \left( \phi \right) + \mathbf{V}^\mathrm{T} \left( \tau \right) \boldsymbol{\Phi} \overline{\mathbf{n}}_6^{e \mathrm{T}} \mathbf{n}^e \mathbf{n}^{e \mathrm{T}}  \mathbf{E} \left( \phi \right) \right) \mathrm{d} \Gamma \\
& \qquad + \int_{\partial \Omega_\phi} \frac{\eta_{13}}{h_e} \tau \phi \mathrm{d} \Gamma + \int_{\Gamma_h} \frac{\eta_{23}}{h_e} \left[ \! \left[ \tau \right] \! \right] \left[ \!\left[ \phi \right] \! \right] \mathrm{d} \Gamma \\
& \qquad + \int_{\partial \Omega_P} \eta_{14} h_e \mathbf{E}^\mathrm{T} \left( \tau \right) \mathbf{n}^{e} \mathbf{n}^{e \mathrm{T}}  \mathbf{E} \left( \phi \right) \mathrm{d} \Gamma + \int_{\Gamma_h} \eta_{24} h_e \left[ \! \left[ \mathbf{E}^\mathrm{T} \left( \tau \right) \mathbf{n}^{e} \right] \! \right] \overline{\mathbf{n}}_4^{e \mathrm{T}} \overline{\mathbf{n}}_4^e \left[ \! \left[ \mathbf{n}^{e \mathrm{T}}  \mathbf{E} \left( \phi \right)  \right] \! \right] \mathrm{d} \Gamma 
\end{split} \\
\begin{split} \nonumber
& \qquad + \sum_{C_h} \left[ \! \left[ \tau \right] \! \right] \mathbf{s}^{e \mathrm{T}} \overline{\mathbf{n}}_6^e \boldsymbol{\Phi} \left\{ \mathbf{V} \left( \phi \right)  \right\} + \sum_{C_h} \left\{ \mathbf{V}^\mathrm{T} \left( \tau \right) \right\} \boldsymbol{\Phi} \overline{\mathbf{n}}_6^{e \mathrm{T}} \mathbf{s}^e  \left[ \! \left[ \phi  \right] \! \right] \\
& \qquad + \sum_{C_h} \left[ \! \left[ \tau \right] \! \right] \mathbf{s}^{e \mathrm{T}} \overline{\mathbf{n}}_6^e \mathbf{b}^\mathrm{T} \left\{ \boldsymbol{\varepsilon} \left(  \mathbf{u} \right)  \right\} + \sum_{C_h} \left\{ \mathbf{E}^\mathrm{T} \left( \tau \right) \right\} \mathbf{a}^\mathrm{T} \overline{\mathbf{n}}_3^{e \mathrm{T}} \overline{\mathbf{s}}_4^{e \mathrm{T}} \left[ \! \left[ \mathbf{u}  \right] \! \right] \\
& \qquad = - \int_\Omega \tau q \mathrm{d} \Omega - \int_{\partial \Omega_\omega} \tau \widetilde{\omega} \mathrm{d} \Gamma - \int_{\partial \Omega_Z}  \mathbf{E}^\mathrm{T} \left( \tau \right) \mathbf{n}^e \widetilde{Z} \mathrm{d} \Gamma \\
& \qquad + \int_{\partial \Omega_u}  \mathbf{E}^\mathrm{T} \left( \tau \right) \left( \mathbf{e}^\mathrm{T} \overline{\mathbf{n}}_1^{e \mathrm{T}} - \widehat{\mathbf{D}}^\mathrm{T} \left( \mathbf{u}, \phi \right) \overline{\mathbf{n}}_4^{e \mathrm{T}} \right) \widetilde{\mathbf{u}} \mathrm{d} \Gamma \\
& \qquad + \int_{\partial \Omega_u}  \mathbf{V}^\mathrm{T} \left(  \tau \right) \left[ \begin{matrix} \mathbf{b}^\mathrm{T} \overline{\mathbf{n}}_1^{e \mathrm{T}} - \overline{\mathbf{c}}_5^\mathrm{T} \mathbf{a}^\mathrm{T} \left( \overline{\mathbf{n}}_{21}^{e \mathrm{T}} + \overline{\mathbf{n}}_3^{e \mathrm{T}} \overline{\mathbf{s}}_4^{e \mathrm{T}} \overline{\mathbf{s}}_4^e  \overline{\mathbf{c}}_1^\mathrm{T} \right) \\ - \overline{\mathbf{c}}_6^\mathrm{T} \mathbf{a}^\mathrm{T} \left( \overline{\mathbf{n}}_{22}^{e \mathrm{T}} + \overline{\mathbf{n}}_3^{e \mathrm{T}} \overline{\mathbf{s}}_4^{e \mathrm{T}} \overline{\mathbf{s}}_4^e  \overline{\mathbf{c}}_2^\mathrm{T} \right) - \widehat{\mathbf{Q}}^\mathrm{T} \left( \mathbf{u}, \phi \right) \overline{\mathbf{n}}_4^{e \mathrm{T}} \end{matrix} \right] \widetilde{\mathbf{u}} \mathrm{d} \Gamma \\
& \qquad  - \int_{\partial \Omega_\phi} \left[ \begin{matrix} \mathbf{E}^\mathrm{T} \left( \tau \right) \boldsymbol{\Lambda} \mathbf{n}^{e} - \mathbf{V}_{,1}^\mathrm{T} \left( \tau \right) \boldsymbol{\Phi} \left( \overline{\mathbf{n}}_{51}^{e \mathrm{T}} +  \overline{\mathbf{n}}_6^{e \mathrm{T}} \mathbf{s}^e \mathbf{s}^{e \mathrm{T}} \overline{\mathbf{c}}_3^\mathrm{T} \right)  \\ - \mathbf{V}_{,2}^\mathrm{T} \left( \tau \right)  \boldsymbol{\Phi} \left( \overline{\mathbf{n}}_{52}^{e \mathrm{T}} + \overline{\mathbf{n}}_6^{e \mathrm{T}} \mathbf{s}^e \mathbf{s}^{e \mathrm{T}} \overline{\mathbf{c}}_4^\mathrm{T} \right) \end{matrix} \right] \widetilde{\phi} \mathrm{d} \Gamma \\
& \qquad + \int_{\partial \Omega_d} \mathbf{E}^\mathrm{T} \left( \tau \right) \mathbf{a}^\mathrm{T} \overline{\mathbf{n}}_3^{e \mathrm{T}} \overline{\mathbf{n}}_4^{e \mathrm{T}} \widetilde{\mathbf{d}} \mathrm{d} \Gamma - \int_{\partial \Omega_P} \mathbf{V}^\mathrm{T} \left( \tau \right) \boldsymbol{\Phi} \overline{\mathbf{n}}_6^{e \mathrm{T}} \mathbf{n}^e \widetilde{P} \mathrm{d} \Gamma \\
& \qquad + \int_{\partial \Omega_\phi} \frac{\eta_{13}}{h_e} \tau \widetilde{\phi} \mathrm{d} \Gamma - \int_{\partial \Omega_P} \eta_{14} h_e \mathbf{E}^\mathrm{T} \left( \tau \right) \mathbf{n}^{e} \widetilde{P} \mathrm{d} \Gamma.
\end{split}
\end{align}
where the normal and constant matrices $\overline{\mathbf{n}}_i$ and $\overline{\mathbf{c}}_i$ can be seen in \ref{app: Mats}. $h_e$ is a boundary-dependent parameter with the unit of length. Here it is defined as the distance between the internal points in subdomains sharing the boundary for $e \in \Gamma_h$, and the smallest distance between the centroid of the subdomain and the external boundary for $e \in \partial \Omega$. The penalty parameters are independent of the boundary size, in which $\left( \eta_{11}, \eta_{12}, \eta_{21}, \eta_{22} \right)$ have the same unit as of the Young's modulus, and $\left( \eta_{13}, \eta_{14}, \eta_{23}, \eta_{24} \right)$ have the same unit as of the permittivity of the material.

As has been stated, the method is only stable when the penalty parameter $\boldsymbol{\eta}_2$ is large enough. However, unlike the DG methods in which the element shape functions are completely independent, a small penalty parameter $\boldsymbol{\eta}_2$ is enough to stabilize the present FPM. On the other hand, the accuracy may decrease if $\boldsymbol{\eta}_2$ is too large. $\boldsymbol{\eta}_1$ on the other hand, should be larger, which corresponds to the essential boundary conditions.

It should be noted that the essential boundary conditions can also be imposed in an alternative way in the FPM. When boundary points are employed, i.e., for $\partial E_i \cap \partial \Omega \not= \varnothing$, $P_i \in \left( \partial E_i \cap \partial \Omega \right)$, we can enforce $\mathbf{u} = \widetilde{\mathbf{u}}$ and $\phi = \widetilde{\phi}$ strongly at the boundary points and thus the Internal Penalty (IP) terms relating to $\eta_{11}$ and $\eta_{13}$ in Eqn.~\ref{eqn:dis_full} and \ref{eqn:ele_full} vanish. However, when high-order essential boundary conditions are under study ($\partial \Omega_d \cup \partial \Omega_P \not= \varnothing$), the IP terms are still required. Note that the IP Numerical Flux Corrections can enforce the essential boundary conditions with any point distributions, regardless of whether there are points distributed on the boundary. 

At last, the formula of the FPM can be converted into the following form:
\begin{align} 
\begin{split}
\left[ \begin{matrix} \mathbf{K}_{uu}^P \left( \overline{\mathbf{u}}, \overline{\boldsymbol{\phi}} \right)  & \mathbf{K}_{u \phi}^P \left( \overline{\mathbf{u}}, \overline{\boldsymbol{\phi}} \right) \\ \left( \mathbf{K}_{u \phi}^P \left( \overline{\mathbf{u}}, \overline{\boldsymbol{\phi}} \right) \right)^\mathrm{T} & \mathbf{K}_{\phi \phi}^P \left( \overline{\mathbf{u}}, \overline{\boldsymbol{\phi}} \right) \end{matrix} \right] \left[ \begin{matrix} \overline{\mathbf{u}} \\  \overline{\boldsymbol{\phi}} \end{matrix} \right] = \left[ \begin{matrix} \mathbf{f}_u^P \left( \overline{\mathbf{u}}, \overline{\boldsymbol{\phi}} \right) \\  \mathbf{f}_\phi^P \left( \overline{\mathbf{u}}, \overline{\boldsymbol{\phi}} \right) \end{matrix} \right]
\end{split},\quad \text{or} \; \mathbf{K}^P \left( \overline{\mathbf{x}} \right) \overline{\mathbf{x}} = \mathbf{f}^P \left( \overline{\mathbf{x}} \right), \label{eqn:full}
\end{align} 
where $ \left( \overline{\mathbf{u}}, \overline{\boldsymbol{\phi}} \right)$ are the global nodal displacement and electric potential vectors. $\mathbf{K}^P$ is the global stiffness matrix for the primal FPM, assembled by a series of point stiffness matrices ($\mathbf{K}_E^P$) and boundary stiffness matrices ($\mathbf{K}_h^P$, $\mathbf{K}_u^P$, $\mathbf{K}_Q^P$, $\mathbf{K}_d^P$, $\mathbf{K}_R^P$, $\mathbf{K}_\phi^P$, $\mathbf{K}_\omega^P$, $\mathbf{K}_P^P$ and $\mathbf{K}_Z^P$ corresponding to $\Gamma_h$, $\partial \Omega_u$, $\partial \Omega_Q$, $\partial \Omega_d$, $\partial \Omega_R$, $\partial \Omega_\phi$, $\partial \Omega_\omega$, $\partial \Omega_P$ and $\partial \Omega_Z$ respectively). $\mathbf{f}^P$ is the global load vector. When evaluating these stiffness matrices and load vectors, we use simple Gaussian quadrature rules. Numerical examples have shown that only one integration point in each subdomain and only one integration point on each internal and external boundary can result in good solutions with acceptable accuracies. However, $2 \times 2$ or more integration points may help to further increase the accuracy, but only one point integration schemes are employed in the present paper. Further numerical implementation in details are omitted here.

It should be pointed out that Eqn.~\ref{eqn:full} based on the full flexoelectric theory is nonlinear, i.e., the global stiffness matrix $\mathbf{K}^P$ and load vector $\mathbf{f}^P$ are functions of the unknown displacement and electric potential fields. As can be seen from Eqn.~\ref{eqn:dis_full} and \ref{eqn:ele_full}, the Point and boundary stiffness matrices and load vectors are dependent on $\widetilde{\mathbf{D}} \left( \mathbf{u}, \phi \right)$ and $\widetilde{\mathbf{Q}} \left( \mathbf{u}, \phi \right)$ in each subdomain. As a result, the global stiffness matrix and load vector has a linear relationship with $\overline{\mathbf{x}}$. The Newton-Raphson method is then applied to solve the nonlinear equation. Notice that when the initial guess $\overline{\mathbf{x}}^{(0)} = \mathbf{0}$, the first approximation $\overline{\mathbf{x}}^{(1)} = \left[ \mathbf{K}^P \left( \mathbf{0} \right) \right]^{-1} \mathbf{f}^P \left( \mathbf{0} \right) $ is equivalent to the linear solution when the electrostatic stress $\boldsymbol{\sigma}^{ES}$ is neglected.

\subsubsection{Reduced theory}

In the reduced formulation (Eqn.~\ref{eqn:gov_red_dis} -- \ref{eqn:BC_red_6}), both the generalized electrostatic stress $\boldsymbol{\sigma}^{ES}$ and electric field gradient $\mathbf{V}$ are neglected. Thus the trial and test functions of the electric potential can be written in a quadratic form:
\begin{align}
\begin{split}
\phi^h \left( x, y \right) & = \mathbf{N}_\phi \left( x, y \right) \boldsymbol{\phi}_E, \qquad \mathbf{N}_\phi \left( x, y \right) = \overline{\mathbf{N}}_R \left( x, y \right) \cdot \overline{\mathbf{B}} + \left[ 1, 0, \cdots, 0 \right]_{1 \times \left( m + 1 \right)}, \\
\overline{\mathbf{N}}_R \left( x, y \right) & = \left[ x - x_0, y - y_0, \frac{1}{2} \left( x - x_0 \right)^2,  \left( x - x_0 \right) \left( y - y_0 \right), \frac{1}{2} \left( y - y_0 \right)^2, 0, 0, 0, 0 \right].
\end{split}
\end{align}

Following the same steps in section~\ref{sec:PFPM_full}, we can obtain the linear primal FPM formula for the reduced flexoelectric theory:

\begin{align}
\begin{split}
\label{eqn:dis_red}
& \int_\Omega \boldsymbol{\varepsilon}^\mathrm{T} \left( \mathbf{v} \right) \mathbf{C}_{\sigma \varepsilon} \boldsymbol{\varepsilon} \left( \mathbf{u} \right) \mathrm{d} \Omega +  \int_\Omega \boldsymbol{\kappa}^\mathrm{T} \left( \mathbf{v} \right) \mathbf{C}_{\mu \kappa} \boldsymbol{\kappa} \left( \mathbf{u} \right) \mathrm{d} \Omega - \int_\Omega \left( \boldsymbol{\varepsilon}^\mathrm{T} \left( \mathbf{v} \right) \mathbf{e} + \boldsymbol{\kappa}^\mathrm{T} \left( \mathbf{v} \right) \mathbf{a} \right) \mathbf{E} \left( \phi \right) \mathrm{d} \Omega \\
& \qquad   -\int_{\Gamma_h \cup \partial \Omega_u} \left( \left[ \! \left[ \mathbf{v}^\mathrm{T} \right] \! \right] \overline{\mathbf{n}}_1^e \mathbf{C}_{\sigma \varepsilon} \left\{ \boldsymbol{\varepsilon} \left( \mathbf{u} \right) \right\} + \left\{ \boldsymbol{\varepsilon}^\mathrm{T} \left( \mathbf{v} \right) \right\} \mathbf{C}_{\sigma \varepsilon} \overline{\mathbf{n}}_1^{e \mathrm{T}} \left[ \! \left[ \mathbf{u} \right] \! \right]  \right) \mathrm{d} \Gamma \\
& \qquad   + \int_{\Gamma_h \cup \partial \Omega_u} \left( \left[ \! \left[ \mathbf{v}^\mathrm{T} \right] \! \right] \overline{\mathbf{n}}_{21}^e \mathbf{C}_{\mu \kappa} \left\{ \boldsymbol{\kappa}_{,1} \left( \mathbf{u} \right) \right\} + \left\{ \boldsymbol{\kappa}_{,1}^\mathrm{T} \left( \mathbf{v} \right) \right\} \mathbf{C}_{\mu \kappa} \overline{\mathbf{n}}_{21}^{e \mathrm{T}} \left[ \! \left[ \mathbf{u} \right] \! \right]  \right) \mathrm{d} \Gamma \\
& \qquad   + \int_{\Gamma_h \cup \partial \Omega_u} \left( \left[ \! \left[ \mathbf{v}^\mathrm{T} \right] \! \right] \overline{\mathbf{n}}_{22}^e \mathbf{C}_{\mu \kappa} \left\{ \boldsymbol{\kappa}_{,2} \left( \mathbf{u} \right) \right\} + \left\{ \boldsymbol{\kappa}_{,2}^\mathrm{T} \left( \mathbf{v} \right) \right\} \mathbf{C}_{\mu \kappa} \overline{\mathbf{n}}_{22}^{e \mathrm{T}} \left[ \! \left[ \mathbf{u} \right] \! \right]  \right) \mathrm{d} \Gamma \\
& \qquad   - \int_{\Gamma_h \cup \partial \Omega_u} \left[ \! \left[ \mathbf{v}^\mathrm{T} \right] \! \right]  \left\{ \overline{\mathbf{n}}_{21}^e \mathbf{a} \mathbf{E}_{,1} \left( \phi \right) + \overline{\mathbf{n}}_{22}^e \mathbf{a} \mathbf{E}_{,2} \left( \phi \right)  \right\} \mathrm{d} \Gamma \\
& \qquad + \int_{\Gamma_h \cup \partial \Omega_u} \left[ \! \left[ \mathbf{v}^\mathrm{T} \right] \! \right] \overline{\mathbf{n}}_1^e \mathbf{e} \left\{ \mathbf{E} \left( \phi \right) \right\} \mathrm{d} \Gamma +\int_{\Gamma_h} \left[ \! \left[ \hat{\boldsymbol{\varepsilon}}^\mathrm{T} \left( \mathbf{v} \right) \right] \! \right] \overline{\mathbf{n}}_3^e \mathbf{a} \left\{ \mathbf{E} \left( \phi \right) \right\}  \mathrm{d} \Gamma \\
& \qquad -\int_{\Gamma_h} \left( \left[ \! \left[ \hat{\boldsymbol{\varepsilon}}^\mathrm{T} \left( \mathbf{v} \right) \right] \! \right] \overline{\mathbf{n}}_3^e \mathbf{C}_{\mu \kappa} \left\{ \boldsymbol{\kappa} \left( \mathbf{u} \right) \right\} + \left\{ \boldsymbol{\kappa}^\mathrm{T} \left( \mathbf{v} \right) \right\} \mathbf{C}_{\mu \kappa} \overline{\mathbf{n}}_3^{e \mathrm{T}} \left[ \! \left[ \hat{\boldsymbol{\varepsilon}} \left( \mathbf{u} \right) \right] \! \right] \right)  \mathrm{d} \Gamma \\
& \qquad - \int_{\Gamma_h \cup \partial \Omega_\phi} \left\{ \boldsymbol{\varepsilon}^\mathrm{T} \left( \mathbf{v} \right) \right\}  \mathbf{e} \mathbf{n}^{e}   \left[ \! \left[ \phi \right] \! \right] \mathrm{d} \Gamma - \int_{\Gamma_h \cup \partial \Omega_\phi} \left\{ \boldsymbol{\kappa}^\mathrm{T} \left( \mathbf{v} \right) \right\}  \mathbf{a} \mathbf{n}^{e}  \left[ \! \left[ \phi \right] \! \right] \mathrm{d} \Gamma \\
& \qquad + \int_{\partial \Omega_u} \left( \mathbf{v}^\mathrm{T} \overline{\mathbf{c}}_1 \overline{\mathbf{s}}_4^{e \mathrm{T}} \overline{\mathbf{s}}_4^e \overline{\mathbf{n}}_3^e \mathbf{C}_{\mu \kappa} \boldsymbol{\kappa}_{,1} \left( \mathbf{u} \right) + \boldsymbol{\kappa}_{,1}^\mathrm{T} \left( \mathbf{v} \right) \mathbf{C}_{\mu \kappa} \overline{\mathbf{n}}_3^{e \mathrm{T}} \overline{\mathbf{s}}_4^{e \mathrm{T}} \overline{\mathbf{s}}_4^e  \overline{\mathbf{c}}_1^\mathrm{T} \mathbf{u}  \right)  \mathrm{d} \Gamma \\
& \qquad + \int_{\partial \Omega_u} \left( \mathbf{v}^\mathrm{T} \overline{\mathbf{c}}_2 \overline{\mathbf{s}}_4^{e \mathrm{T}} \overline{\mathbf{s}}_4^e \overline{\mathbf{n}}_3^e \mathbf{C}_{\mu \kappa} \boldsymbol{\kappa}_{,2} \left( \mathbf{u} \right) + \boldsymbol{\kappa}_{,2}^\mathrm{T} \left( \mathbf{v} \right) \mathbf{C}_{\mu \kappa} \overline{\mathbf{n}}_3^{e \mathrm{T}} \overline{\mathbf{s}}_4^{e \mathrm{T}} \overline{\mathbf{s}}_4^e  \overline{\mathbf{c}}_2^\mathrm{T} \mathbf{u}  \right)  \mathrm{d} \Gamma \\
& \qquad - \int_{\partial \Omega_u} \mathbf{v}^\mathrm{T} \left( \overline{\mathbf{c}}_1 \overline{\mathbf{s}}_4^{e \mathrm{T}} \overline{\mathbf{s}}_4^e \overline{\mathbf{n}}_3^e \mathbf{a} \mathbf{E}_{,1} \left( \phi \right) + \overline{\mathbf{c}}_2 \overline{\mathbf{s}}_4^{e \mathrm{T}} \overline{\mathbf{s}}_4^e \overline{\mathbf{n}}_3^e \mathbf{a}  \mathbf{E}_{,2} \left( \phi \right) \right) \mathrm{d} \Gamma \\
& \qquad - \int_{\partial \Omega_d} \left( \hat{\boldsymbol{\varepsilon}}^\mathrm{T} \left( \mathbf{v} \right) \overline{\mathbf{n}}_4^{e \mathrm{T}} \overline{\mathbf{n}}_4^e \overline{\mathbf{n}}_3^e \mathbf{C}_{\mu \kappa} \boldsymbol{\kappa} \left( \mathbf{u} \right) + \boldsymbol{\kappa}^\mathrm{T} \left( \mathbf{v} \right) \mathbf{C}_{\mu \kappa} \overline{\mathbf{n}}_3^{e \mathrm{T}} \overline{\mathbf{n}}_4^{e \mathrm{T}} \overline{\mathbf{n}}_4^e \hat{\boldsymbol{\varepsilon}} \left( \mathbf{u} \right) \right)  \mathrm{d} \Gamma \\
& \qquad + \int_{\partial \Omega_d} \hat{\boldsymbol{\varepsilon}}^\mathrm{T} \left( \mathbf{v} \right) \overline{\mathbf{n}}_4^{e \mathrm{T}} \overline{\mathbf{n}}_4^e \overline{\mathbf{n}}_3^e \mathbf{a} \mathbf{E} \left( \phi \right) \mathrm{d} \Gamma + \int_{\partial \Omega_u} \frac{\eta_{11}}{h_e} \mathbf{v}^\mathrm{T} \mathbf{u} \mathrm{d} \Gamma + \int_{\Gamma_h} \frac{\eta_{21}}{h_e} \left[ \! \left[ \mathbf{v}^\mathrm{T} \right] \! \right] \left[ \!\left[ \mathbf{u} \right] \! \right] \mathrm{d} \Gamma \\
& \qquad + \int_{\partial \Omega_d} \eta_{12} h_e \hat{\boldsymbol{\varepsilon}}^\mathrm{T} \left( \mathbf{v} \right) \overline{\mathbf{n}}_4^{e \mathrm{T}} \overline{\mathbf{n}}_4^e \hat{\boldsymbol{\varepsilon}} \left( \mathbf{u} \right) \mathrm{d} \Gamma + \int_{\Gamma_h} \eta_{22} h_e \left[ \! \left[ \hat{\boldsymbol{\varepsilon}}^\mathrm{T} \left( \mathbf{v} \right) \right] \! \right] \overline{\mathbf{n}}_4^{e \mathrm{T}} \overline{\mathbf{n}}_4^e \left[ \! \left[ \hat{\boldsymbol{\varepsilon}} \left( \mathbf{u} \right) \right] \! \right] \mathrm{d} \Gamma \\
& \qquad - \sum_{C_h} \left[ \! \left[ \mathbf{v}^\mathrm{T} \right] \! \right] \overline{\mathbf{s}}_4^e \overline{\mathbf{n}}_3^e  \mathbf{C}_{\mu \kappa} \left\{ \boldsymbol{\kappa} \left( \mathbf{u} \right) \right\}  - \sum_{C_h} \left\{ \boldsymbol{\kappa}^\mathrm{T} \left( \mathbf{v} \right) \right\} \mathbf{C}_{\mu \kappa} \overline{\mathbf{n}}_3^{e \mathrm{T}} \overline{\mathbf{s}}_4^{e \mathrm{T}} \left[ \! \left[ \mathbf{u}  \right] \! \right] \\
& \qquad +  \sum_{C_h} \left[ \! \left[ \mathbf{v}^\mathrm{T} \right] \! \right] \overline{\mathbf{s}}_4^e \overline{\mathbf{n}}_3^e \mathbf{a} \left\{ \mathbf{E} \left( \phi \right)   \right\}
\end{split} \\
\begin{split} \nonumber
& \qquad  = \int_\Omega \mathbf{v}^\mathrm{T} \mathbf{b} \mathrm{d} \Omega + \int_{\partial \Omega_Q} \mathbf{v}^\mathrm{T} \widetilde{\mathbf{Q}} \mathrm{d} \Gamma + \int_{\partial \Omega_R} \hat{\boldsymbol{\varepsilon}}^\mathrm{T} \left( \mathbf{v} \right) \overline{\mathbf{n}}_4^{e \mathrm{T}} \widetilde{\mathbf{R}} \mathrm{d} \Gamma \\
& \qquad -\int_{\partial \Omega_u}  \left[ \begin{matrix} \boldsymbol{\varepsilon}^\mathrm{T} \left( \mathbf{v} \right) \mathbf{C}_{\sigma \varepsilon} \overline{\mathbf{n}}_1^{e \mathrm{T}} - \boldsymbol{\kappa}_{,1}^\mathrm{T} \left( \mathbf{v} \right) \mathbf{C}_{\mu \kappa} \left( \overline{\mathbf{n}}_{21}^{e \mathrm{T}} + \overline{\mathbf{n}}_3^{e \mathrm{T}} \overline{\mathbf{s}}_4^{e \mathrm{T}} \overline{\mathbf{s}}_4^e  \overline{\mathbf{c}}_1^\mathrm{T} \right)  \\ - \boldsymbol{\kappa}_{,2}^\mathrm{T} \left( \mathbf{v} \right) \mathbf{C}_{\mu \kappa} \left( \overline{\mathbf{n}}_{22}^{e \mathrm{T}} + \overline{\mathbf{n}}_3^{e \mathrm{T}} \overline{\mathbf{s}}_4^{e \mathrm{T}} \overline{\mathbf{s}}_4^e  \overline{\mathbf{c}}_2^\mathrm{T} \right) \end{matrix} \right]  \widetilde{\mathbf{u}} \mathrm{d} \Gamma \\
& \qquad  - \int_{\partial \Omega_\phi} \left( \boldsymbol{\varepsilon}^\mathrm{T} \left( \mathbf{v} \right) \mathbf{e} \mathbf{n}^{e} + \boldsymbol{\kappa}^\mathrm{T} \left( \mathbf{v} \right) \mathbf{a} \mathbf{n}^{e} \right)  \widetilde{\phi} \mathrm{d} \Gamma - \int_{\partial \Omega_d} \boldsymbol{\kappa}^\mathrm{T} \left( \mathbf{v} \right) \mathbf{C}_{\mu \kappa} \overline{\mathbf{n}}_3^{e \mathrm{T}} \overline{\mathbf{n}}_4^{e \mathrm{T}} \widetilde{\mathbf{d}}  \mathrm{d} \Gamma \\
& \qquad + \int_{\partial \Omega_u} \frac{\eta_{11}}{h_e} \mathbf{v}^\mathrm{T} \widetilde{\mathbf{u}} \mathrm{d} \Gamma + \int_{\partial \Omega_d} \eta_{12} h_e \hat{\boldsymbol{\varepsilon}}^\mathrm{T} \left( \mathbf{v} \right) \overline{\mathbf{n}}_4^{e \mathrm{T}} \widetilde{\mathbf{d}} \mathrm{d} \Gamma,
\end{split}  \\~\nonumber \\
\begin{split}
\label{eqn:ele_red}
& - \int_\Omega \mathbf{E}^\mathrm{T} \left( \tau \right) \boldsymbol{\Lambda} \mathbf{E} \left( \phi \right) \mathrm{d} \Omega - \int_\Omega  \mathbf{E}^\mathrm{T} \left( \tau \right) \left( \mathbf{e}^\mathrm{T}  \boldsymbol{\varepsilon} \left( \mathbf{u} \right) + \mathbf{a}^\mathrm{T}  \boldsymbol{\kappa} \left( \mathbf{u} \right) \right) \mathrm{d} \Omega \\
& \qquad - \int_{\Gamma_h \cup \partial \Omega_u}  \left\{ \mathbf{E}_{,1}^\mathrm{T} \left( \tau \right) \mathbf{a}^\mathrm{T} \overline{\mathbf{n}}_{21}^{e \mathrm{T}} + \mathbf{E}_{,2}^\mathrm{T}  \left( \tau \right) \mathbf{a}^\mathrm{T} \overline{\mathbf{n}}_{22}^{e \mathrm{T}} \right\} \left[ \! \left[ \mathbf{u} \right] \! \right] \mathrm{d} \Gamma \\
& \qquad + \int_{\Gamma_h \cup \partial \Omega_u}  \left\{ \mathbf{E}^\mathrm{T} \left( \tau \right) \right\}  \mathbf{e}^\mathrm{T} \overline{\mathbf{n}}_1^{e \mathrm{T}} \left[ \! \left[ \mathbf{u} \right] \! \right] \mathrm{d} \Gamma + \int_{\Gamma_h} \left\{ \mathbf{E}^\mathrm{T} \left( \tau \right) \right\} \mathbf{a}^\mathrm{T} \overline{\mathbf{n}}_3^{e \mathrm{T}} \left[ \! \left[ \hat{\boldsymbol{\varepsilon}} \left( \mathbf{u} \right) \right] \! \right]   \mathrm{d} \Gamma  \\
& \qquad  - \int_{\Gamma_h \cup \partial \Omega_\phi} \left[ \! \left[ \tau \right] \! \right] \mathbf{n}^{e \mathrm{T}} \mathbf{e}^\mathrm{T} \left\{ \boldsymbol{\varepsilon} \left( \mathbf{u} \right) \right\} \mathrm{d} \Gamma  - \int_{\Gamma_h \cup \partial \Omega_\phi} \left[ \! \left[ \tau \right] \! \right]  \mathbf{n}^{e \mathrm{T}} \mathbf{a}^\mathrm{T} \left\{ \boldsymbol{\kappa} \left( \mathbf{u} \right) \right\} \mathrm{d} \Gamma \\
& \qquad  - \int_{\Gamma_h \cup \partial \Omega_\phi} \left(  \left[ \! \left[ \tau \right] \! \right] \mathbf{n}^{e \mathrm{T}} \boldsymbol{\Lambda} \left\{ \mathbf{E} \left( \phi \right) \right\} + \left\{ \mathbf{E}^\mathrm{T} \left( \tau \right) \right\} \boldsymbol{\Lambda}  \mathbf{n}^{e} \left[ \! \left[ \phi \right] \! \right] \right) \mathrm{d} \Gamma \\
& \qquad - \int_{\partial \Omega_u} \left( \mathbf{E}_{,1}^\mathrm{T} \left( \tau \right) \mathbf{a}^\mathrm{T} \overline{\mathbf{n}}_3^{e \mathrm{T}} \overline{\mathbf{s}}_4^{e \mathrm{T}} \overline{\mathbf{s}}_4^e  \overline{\mathbf{c}}_1^\mathrm{T} +   \mathbf{E}_{,2}^\mathrm{T} \left( \tau \right)  \mathbf{a}^\mathrm{T} \overline{\mathbf{n}}_3^{e \mathrm{T}} \overline{\mathbf{s}}_4^{e \mathrm{T}} \overline{\mathbf{s}}_4^e  \overline{\mathbf{c}}_2^\mathrm{T} \right) \mathbf{u} \mathrm{d} \Gamma \\
& \qquad + \int_{\partial \Omega_d} \mathbf{E}^\mathrm{T} \left( \tau \right) \mathbf{a}^\mathrm{T} \overline{\mathbf{n}}_3^{e \mathrm{T}} \overline{\mathbf{n}}_4^{e \mathrm{T}} \overline{\mathbf{n}}_4^e \hat{\boldsymbol{\varepsilon}} \left( \mathbf{u} \right)  \mathrm{d} \Gamma+ \int_{\partial \Omega_\phi} \frac{\eta_{13}}{h_e} \tau \phi \mathrm{d} \Gamma + \int_{\Gamma_h} \frac{\eta_{23}}{h_e} \left[ \! \left[ \tau \right] \! \right] \left[ \!\left[ \phi \right] \! \right] \mathrm{d} \Gamma\\
& \qquad + \sum_{C_h} \left\{ \mathbf{E}^\mathrm{T} \left( \tau \right) \right\} \mathbf{a}^\mathrm{T} \overline{\mathbf{n}}_3^{e \mathrm{T}} \overline{\mathbf{s}}_4^{e \mathrm{T}} \left[ \! \left[ \mathbf{u}  \right] \! \right] \\
& \qquad = - \int_\Omega \tau q \mathrm{d} \Omega - \int_{\partial \Omega_\omega} \tau \widetilde{\omega} \mathrm{d} \Gamma + \int_{\partial \Omega_u}  \mathbf{E}^\mathrm{T} \left( \tau \right) \mathbf{e}^\mathrm{T} \overline{\mathbf{n}}_1^{e \mathrm{T}}  \widetilde{\mathbf{u}} \mathrm{d} \Gamma \\
& \qquad - \int_{\partial \Omega_u} \left[ \mathbf{E}_{,1}^\mathrm{T} \left(  \tau \right)  \mathbf{a}^\mathrm{T} \left( \overline{\mathbf{n}}_{21}^{e \mathrm{T}} + \overline{\mathbf{n}}_3^{e \mathrm{T}} \overline{\mathbf{s}}_4^{e \mathrm{T}} \overline{\mathbf{s}}_4^e  \overline{\mathbf{c}}_1^\mathrm{T} \right) + \mathbf{E}_{,2}^\mathrm{T} \left(  \tau \right)  \mathbf{a}^\mathrm{T} \left( \overline{\mathbf{n}}_{22}^{e \mathrm{T}} + \overline{\mathbf{n}}_3^{e \mathrm{T}} \overline{\mathbf{s}}_4^{e \mathrm{T}} \overline{\mathbf{s}}_4^e  \overline{\mathbf{c}}_2^\mathrm{T} \right)  \right] \widetilde{\mathbf{u}} \mathrm{d} \Gamma \\
& \qquad  - \int_{\partial \Omega_\phi}  \mathbf{E}^\mathrm{T} \left( \tau \right) \boldsymbol{\Lambda} \mathbf{n}^{e} \widetilde{\phi} \mathrm{d} \Gamma + \int_{\partial \Omega_d} \mathbf{E}^\mathrm{T} \left( \tau \right) \mathbf{a}^\mathrm{T} \overline{\mathbf{n}}_3^{e \mathrm{T}} \overline{\mathbf{n}}_4^{e \mathrm{T}} \widetilde{\mathbf{d}} \mathrm{d} \Gamma + \int_{\partial \Omega_\phi} \frac{\eta_{13}}{h_e} \tau \widetilde{\phi} \mathrm{d} \Gamma.
\end{split}
\end{align}

In a matrix form:
\begin{align} 
\begin{split}
\left[ \begin{matrix} \mathbf{K}_{uu}^P & \mathbf{K}_{u \phi}^P \\ \left( \mathbf{K}_{u \phi}^P \right)^\mathrm{T} & \mathbf{K}_{\phi \phi}^P \end{matrix} \right] \left[ \begin{matrix} \overline{\mathbf{u}} \\  \overline{\boldsymbol{\phi}} \end{matrix} \right] = \left[ \begin{matrix} \mathbf{f}_u^P \\  \mathbf{f}_\phi^P \end{matrix} \right]
\end{split},\quad \text{or} \; \mathbf{K}^P \overline{\mathbf{x}} = \mathbf{f}^P. \label{eqn:red}
\end{align} 
In the reduced theory, the global stiffness matrix and load vector $\mathbf{K}^P$ and $\mathbf{f}^P$ are independent of the displacement and electric potential field, i.e., Eqn.~\ref{eqn:red} is linear. Besides, the final stiffness matrix is sparse and symmetric, which is beneficial for modeling complex problems with a large number of DoFs.

Whereas the trial and test functions for displacement given in this section are cubic polynomial, numerical study shows that simplified quadratic trial and test functions can also by employed and can achieve good accuracy (with relative error \textless 1\%) for most problems. Thus the third-order terms in Eqn.~\ref{eqn:Shape}, as well as the terms containing $\boldsymbol{\kappa}_{,i} \; \left( i = 1, 2 \right)$ on the internal boundaries in Eqn.~\ref{eqn:dis_full} -- \ref{eqn:ele_red} are neglected. Similarly, the trial and test functions for the electric potential can also be simplified to quadratic (in the full theory) or linear (in the reduced theory). And the corresponding terms containing $\mathbf{V}_{,i} \; \left( i = 1, 2 \right)$ or $\mathbf{E}_{,i} \; \left( i = 1, 2 \right)$ in the above formulations are then absent. The cubic trial and test functions can improve the accuracy of the primal FPM with more terms in the Taylor’s expansion.

\section{Mixed Fragile Points Method (mixed FPM)}\label{sec:mixed}

\subsection{Independent variables and linear trial and test functions}

In this section, we develop a mixed FPM in which the displacement $\mathbf{u}$, electric potential $\phi$, high-order stress $\boldsymbol{\mu}$ and high-order electric displacement $\mathbf{Q}$ are all interpolated independently.

Here we introduce the following independent variables and their trial functions in subdomain $E_0$:
\begin{align}
\mathbf{u}^h \left( x, y \right) & = \mathbf{u}^0 + \left( x - x_0 \right) \mathbf{u}_{,1} \big|_{P_0} + \left( y - y_0 \right) \mathbf{u}_{,2} \big|_{P_0} = \mathbf{N}_u^M \left( x, y \right) \mathbf{u}_E, \label{mixed_Var_1} \\
\phi^h \left( x, y \right) & =  \phi^0 + \left( x - x_0 \right) \phi_{,1} \big|_{P_0} + \left( y - y_0 \right) \phi_{,2} \big|_{P_0} = \mathbf{N}_\phi^M \left( x, y \right) \boldsymbol{\phi}_E, \\
{\boldsymbol{\mu}}^h \left( x, y \right) & = {\boldsymbol{\mu}}^0 + \left( x - x_0 \right) {\boldsymbol{\mu}}_{,1} \big|_{P_0} + \left( y - y_0 \right) {\boldsymbol{\mu}}_{,2} \big|_{P_0} = \mathbf{N}_{{{\mu}}}^M \left( x, y \right) {\boldsymbol{\mu}}_E, \\
\mathbf{Q}^h \left( x, y \right) &= \mathbf{Q}^0 + \left( x - x_0 \right) \mathbf{Q}_{,1} \big|_{P_0} + \left( y - y_0 \right) \mathbf{Q}_{,2} \big|_{P_0} = \mathbf{N}_Q^M \left( x, y \right) \mathbf{Q}_E. \label{mixed_Var_6} 
\end{align}
where
\begin{align} \label{eqn:mixed_B}
\begin{split}
& {\boldsymbol{\mu}}_E = \left[ {\boldsymbol{\mu}}^{0 \mathrm{T}}, {\boldsymbol{\mu}}^{1 \mathrm{T}}, \cdots, {\boldsymbol{\mu}}^{m \mathrm{T}} \right]^\mathrm{T}, \quad \mathbf{Q}_E= \left[  \mathbf{Q}^{0 \mathrm{T}}, \mathbf{Q}^{1 \mathrm{T}}, \cdots , \mathbf{Q}^{m \mathrm{T}} \right]^\mathrm{T},  \\
& \mathbf{N}_u^M \left( x, y \right) = \mathbf{N}_\phi^M \left( x, y \right) \otimes \mathbf{I}_{2 \times 2}, \quad \mathbf{N}_{{{\mu}}}^M \left( x, y \right) =  \mathbf{N}_\phi^M \left( x, y \right) \otimes \mathbf{I}_{6 \times 6}, \\
& \mathbf{N}_{{{Q}}}^M \left( x, y \right) = \mathbf{N}_\phi^M \left( x, y \right) \otimes \mathbf{I}_{4 \times 4}, \quad \mathbf{N}_\phi^M \left( x, y \right) = \overline{\mathbf{N}}^M \left( x, y \right) \cdot \overline{\mathbf{B}}^M + \left[ 1, 0, \cdots, 0 \right]_{1 \times \left( m + 1 \right)}, \\
& \overline{\mathbf{N}}^M \left( x, y \right) = \left[ x - x_0, y - y_0 \right], \quad \overline{\mathbf{B}}^M =\overline{\mathbf{B}} \cdot  \left[ \mathbf{I}_{2 \times 2}, \mathbf{0}_{2 \times 7} \right]^\mathrm{T} .
\end{split}
\end{align}

\begin{figure}[htbp] 
  \centering 
    \subfigure[]{ 
    \label{fig:Shape_1st} 
    \includegraphics[width=0.48\textwidth]{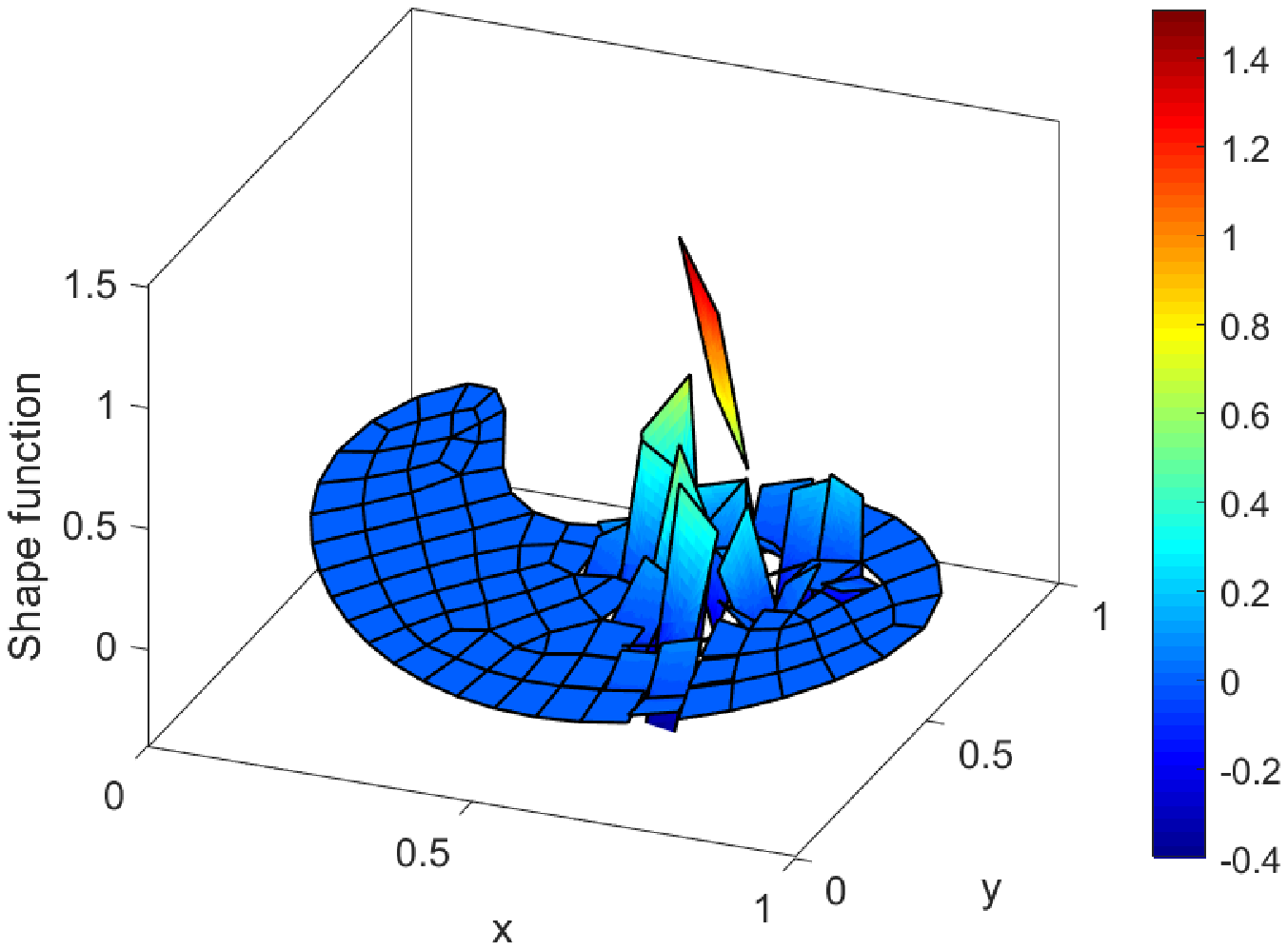}}  
    \subfigure[]{ 
    \label{fig:Trial_1st} 
    \includegraphics[width=0.48\textwidth]{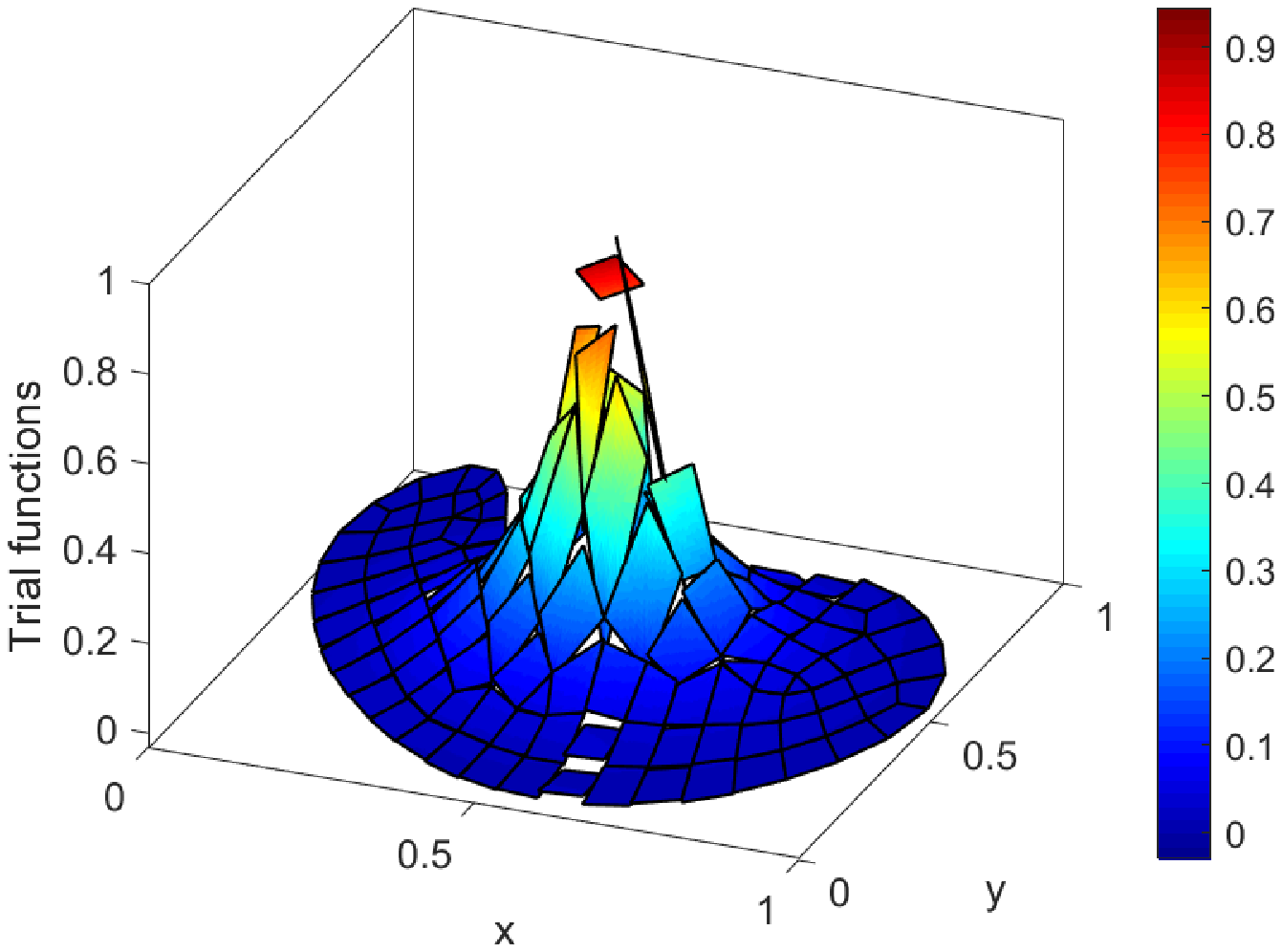}}  
  \caption{The shape and trial function in the mixed FPM. (a) The shape function. (b) The trial function for $u_a = e^{-10 \sqrt{ \left( x – 0.5 \right)^2 + \left( y – 0.5 \right)^2 }}$.} 
  \label{fig:Func_1st} 
\end{figure}

Similar to the primal FPM, the trial functions in the mixed FPM are also piecewise-continuous. Yet for simplicity, these trial functions are just written as a linear Taylor expansion at the corresponding Fragile Point and all the high-order terms are omitted. The graph of the linear piecewise-continuous shape function is shown in Fig.~\ref{fig:Shape_1st}. And the trial function simulating an exponential function $u_a = e^{-10 \sqrt{ \left( x – 0.5 \right)^2 + \left( y – 0.5 \right)^2 }}$ is presented in Fig.~\ref{fig:Trial_1st}. There are various approaches to approximate the first derivatives based on the values of the supporting points. In addition to the local RBF-DQ method shown in section~\ref{sec:DQ}, Generalized Finite Difference (GFD) method can also be applied to generate $\overline{\mathbf{B}}^M$ in Eqn.~\ref{eqn:mixed_B}. The corresponding matrices can be seen in \cite{Guan2020} and are omitted here. The test functions for displacement $\mathbf{v}^h$, electric potential $\tau^h$, high-order stress $\mathbf{w}^h$ and high-order electric displacement $\boldsymbol{\xi}^h$ have the same shape as their corresponding trial functions in the mixed FPM.

\subsection{Weak-form formulation}

\subsubsection{Full theory}

First, we multiply the governing equations (Eqn.~\ref{eq:gov_full_dis} and \ref{eq:gov_full_ele}) by the test functions $\mathbf{v}_i^h$ and $\tau^i$ respectively and integrate in subdomain $E$ by parts:
\begin{align}
\begin{split}
& \int_E \boldsymbol{\varepsilon}^\mathrm{T} \left( \mathbf{v} \right) \boldsymbol{\sigma} \left( \mathbf{u}, \phi, \mathbf{Q} \right) \mathrm{d} \Omega + \int_E \hat{\boldsymbol{\varepsilon}}^\mathrm{T} \left( \mathbf{v} \right) \boldsymbol{\sigma}^{ES} \left( \mathbf{u}, \phi, \mathbf{Q} \right) \mathrm{d} \Omega \\
& \qquad -  \int_E \hat{\boldsymbol{\varepsilon}}^\mathrm{T} \left( \mathbf{v} \right) \overline{\mathbf{c}}_{9}^\mathrm{T} \boldsymbol{\mu}_{,1} \mathrm{d} \Omega -  \int_E \hat{\boldsymbol{\varepsilon}}^\mathrm{T} \left( \mathbf{v} \right) \overline{\mathbf{c}}_{10}^\mathrm{T} \boldsymbol{\mu}_{,2} \mathrm{d} \Omega \\
& \qquad  = \int_E \mathbf{v}^\mathrm{T} \mathbf{b} \mathrm{d} \Omega + \int_{\partial E}  \mathbf{v}^\mathrm{T}  \left( \overline{\mathbf{n}}_1^e \boldsymbol{\sigma} \left( \mathbf{u}, \phi, \mathbf{Q} \right) - \overline{\mathbf{n}}_{21}^e \boldsymbol{\mu}_{,1} - \overline{\mathbf{n}}_{22}^e \boldsymbol{\mu}_{,2} + \overline{\mathbf{n}}_{4}^e \boldsymbol{\sigma}^{ES} \left( \mathbf{u}, \phi, \mathbf{Q} \right)  \right)  \mathrm{d} \Gamma,
\end{split} \\
\begin{split}
& \int_E \mathbf{E}^\mathrm{T} \left( \tau \right) \mathbf{D} \left( \mathbf{u}, \phi, \boldsymbol{\mu} \right) \mathrm{d} \Omega -  \int_E \mathbf{E}^\mathrm{T} \left( \tau \right) 
\overline{\mathbf{c}}_5 \mathbf{Q}_{,1} \mathrm{d} \Omega -  \int_E \mathbf{E}^\mathrm{T} \left( \tau \right) 
\overline{\mathbf{c}}_6 \mathbf{Q}_{,2} \mathrm{d} \Omega \\
& \qquad = \int_E \tau q \mathrm{d} \Omega - \int_{\partial E}  \tau  \left( \mathbf{n}^{e \mathrm{T}} \mathbf{D} \left( \mathbf{u}, \phi, \boldsymbol{\mu} \right) - \overline{\mathbf{n}}_{51}^e  \mathbf{Q}_{,1} - \overline{\mathbf{n}}_{52}^e  \mathbf{Q}_{,2}  \right)  \mathrm{d} \Gamma.
\end{split}
\end{align}
Similarly, the constitutive relations for $\boldsymbol{\mu}$ and $\mathbf{Q}$ (see Eqn.~\ref{eqn:cons_full}) are also multiplied by their corresponding test functions $\mathbf{w}$ and $\boldsymbol{\xi}$ and then integrated by parts in $E$:
\begin{align}
\begin{split}
& \int_E \mathbf{w}^\mathrm{T} \mathbf{C}_{\mu \kappa}^{-1} \boldsymbol{\mu} \mathrm{d} \Omega + \int_E \mathbf{w}_{,1}^\mathrm{T} \overline{\mathbf{c}}_9 \hat{\boldsymbol{\varepsilon}} \left( \mathbf{u} \right) \mathrm{d} \Omega + \int_E \mathbf{w}_{,2}^\mathrm{T} \overline{\mathbf{c}}_{10} \hat{\boldsymbol{\varepsilon}} \left( \mathbf{u} \right) \mathrm{d} \Omega \\
& \qquad - \int_E \mathbf{w}^\mathrm{T} \mathbf{C}_{\mu \kappa}^{-1} \mathbf{a} \mathbf{E} \left( \phi \right) \mathrm{d} \Omega = \int_{\partial E} \mathbf{w}^\mathrm{T} \overline{\mathbf{n}}_3^{e \mathrm{T}} \hat{\boldsymbol{\varepsilon}} \left( \mathbf{u} \right) \mathrm{d} \Gamma,
\end{split} \\
\begin{split}
& \int_E \boldsymbol{\xi}^\mathrm{T} \boldsymbol{\Phi}^{-1} \mathbf{Q} \mathrm{d} \Omega + \int_E \boldsymbol{\xi}_{,1}^\mathrm{T} \overline{\mathbf{c}}_5^\mathrm{T} \mathbf{E} \left( \phi \right) \mathrm{d} \Omega + \int_E \boldsymbol{\xi}_{,2}^\mathrm{T} \overline{\mathbf{c}}_6^\mathrm{T} \mathbf{E} \left( \phi \right) \mathrm{d} \Omega \\
& \qquad - \int_E \boldsymbol{\xi}^\mathrm{T} \boldsymbol{\Phi}^{-1} \mathbf{b}^\mathrm{T} \boldsymbol{\varepsilon} \left( \mathbf{u} \right) \mathrm{d} \Omega = \int_{\partial E} \boldsymbol{\xi}^\mathrm{T} \overline{\mathbf{n}}_6^{e \mathrm{T}} \mathbf{E} \left( \phi \right) \mathrm{d} \Gamma.
\end{split}
\end{align}

The high-order boundary conditions for $\boldsymbol{\mu}$ and $\mathbf{Q}$ on $\partial \Omega_R$ and $\partial \Omega_Z$ are enforced by Lagrange multipliers $\boldsymbol{\lambda}_R = \left[ \lambda_{R1}, \lambda_{R2} \right]$ and $\lambda_Z$ respectively. These Lagrange multipliers are defined on  $\partial \Omega_R$ and $\partial \Omega_Z$ with constant shape function in each subdomain boundary:
\begin{align}
\begin{split}
\boldsymbol{\lambda}_R^h \left( x , y \right) = \boldsymbol{\lambda}_R^i , \quad \lambda_Z^h \left( x , y \right) = \lambda_Z^i, \; \text{for} \left( x , y \right) \in e_i.
\end{split}
\end{align}
The nodal Lagrange multiplier vectors can be written as $\overline{\boldsymbol{\lambda}}_R = \left[ \lambda_{R1}^1, \lambda_{R2}^1, \lambda_{R1}^2, \lambda_{R2}^2, \cdots, \lambda_{R1}^{n_R}, \lambda_{R2}^{n_R} \right]^\mathrm{T}$ and $\overline{\boldsymbol{\lambda}}_Z = \left[ \lambda_Z^1, \lambda_Z^2, \cdots, \lambda_Z^{n_Z} \right] ^\mathrm{T}$, where $n_R$ and $n_Z$ are the number of subdomain boundaries on $\partial \Omega_R$ and $\partial \Omega_Z$ respectively. The test functions of the Lagrange multipliers $\hat{\boldsymbol{\lambda}}_R$ and $\hat{\lambda}_Z$ take the same shape as ${\boldsymbol{\lambda}}_R$ and ${\lambda}_Z$.

Following the same steps in the primal FPM, we can obtain the basic symmetric formula for the mixed FPM:
\begin{align}
\begin{split}
\label{eqn:dis_full_mixed}
& \int_\Omega \boldsymbol{\varepsilon}^\mathrm{T} \left( \mathbf{v} \right) \left( \mathbf{C}_{\sigma \varepsilon} + \mathbf{b} \boldsymbol{\Phi}^{-1} \mathbf{b}^{\mathrm{T}} \right) \boldsymbol{\varepsilon} \left( \mathbf{u} \right) \mathrm{d} \Omega - \int_\Omega \boldsymbol{\varepsilon}^\mathrm{T} \left( \mathbf{v} \right) \mathbf{e} \mathbf{E} \left( \phi \right) \mathrm{d} \Omega  \\
& \qquad - \int_\Omega \boldsymbol{\varepsilon}^\mathrm{T} \left( \mathbf{v} \right) \mathbf{b} \boldsymbol{\Phi}^{-1} \mathbf{Q} \mathrm{d} \Omega - \int_\Omega \hat{\boldsymbol{\varepsilon}}^\mathrm{T} \left( \mathbf{v} \right) \overline{\mathbf{c}}_9^\mathrm{T} \boldsymbol{\mu}_{,1} \mathrm{d} \Omega - \int_\Omega \hat{\boldsymbol{\varepsilon}}^\mathrm{T} \left( \mathbf{v} \right) \overline{\mathbf{c}}_{10}^\mathrm{T} \boldsymbol{\mu}_{,2} \mathrm{d} \Omega \\
& \qquad - \int_{\Gamma_h \cup \partial \Omega_u} \left[ \! \left[ \mathbf{v}^\mathrm{T} \right] \! \right] \left\{ \left( \overline{\mathbf{n}}_1^e \left( \mathbf{C}_{\sigma \varepsilon} + \mathbf{b} \boldsymbol{\Phi}^{-1} \mathbf{b}^\mathrm{T} \right) - \overline{\mathbf{n}}_4^e \widehat{\mathbf{Q}} \boldsymbol{\Phi}^{-1} \mathbf{b}^\mathrm{T} \right) \boldsymbol{\varepsilon} \left( \mathbf{u} \right) \right\} \mathrm{d} \Gamma \\
& \qquad -\int_{\Gamma_h \cup \partial \Omega_u} \left\{ \boldsymbol{\varepsilon}^\mathrm{T} \left( \mathbf{v} \right) \left(\left( \mathbf{C}_{\sigma \varepsilon} + \mathbf{b} \boldsymbol{\Phi}^{-1} \mathbf{b}^\mathrm{T} \right) \overline{\mathbf{n}}_1^{e \mathrm{T}} - \mathbf{b} \boldsymbol{\Phi}^{-1} \widehat{\mathbf{Q}}^\mathrm{T} \overline{\mathbf{n}}_4^{e \mathrm{T}} \right) \right\} \left[ \! \left[ \mathbf{u} \right] \! \right] \mathrm{d} \Gamma \\
& \qquad + \int_{\Gamma_h \cup \partial \Omega_u} \left[ \! \left[ \mathbf{v}^\mathrm{T} \right] \! \right] \left\{ \left( \overline{\mathbf{n}}_1^e \mathbf{e} - \overline{\mathbf{n}}_4^e \widehat{\mathbf{D}} \left( \mathbf{u}, \phi, \boldsymbol{\mu} \right) \right) \mathbf{E} \left( \phi \right) \right\}  \mathrm{d} \Gamma  \\
& \qquad + \int_{\Gamma_h \cup \partial \Omega_u} \left[ \! \left[ \mathbf{v}^\mathrm{T} \right] \! \right] \left\{ \overline{\mathbf{n}}_{21}^e \boldsymbol{\mu}_{,1} + \overline{\mathbf{n}}_{22}^e \boldsymbol{\mu}_{,2} \right\} \mathrm{d} \Gamma + \int_{\Gamma_h} \left\{ \hat{\boldsymbol{\varepsilon}}^\mathrm{T} \left( \mathbf{v} \right) \right\} \overline{\mathbf{n}}_3^e \left[ \! \left[ \boldsymbol{\mu} \right] \! \right] \mathrm{d} \Gamma \\
& \qquad + \int_{\Gamma_h \cup \partial \Omega_u} \left[ \! \left[ \mathbf{v}^\mathrm{T} \right] \! \right] \left\{ \left( \overline{\mathbf{n}}_1^e \mathbf{b} - \overline{\mathbf{n}}_4^e \widehat{\mathbf{Q}} \right)  \boldsymbol{\Phi}^{-1} \mathbf{Q} \right\} \mathrm{d} \Gamma - \int_{\Gamma_h \cup \partial \Omega_\phi} \left\{ \boldsymbol{\varepsilon}^\mathrm{T} \left( \mathbf{v} \right) \right\} \mathbf{e} \mathbf{n}^e \left[ \! \left[ \phi \right] \! \right] \mathrm{d} \Gamma \\
& \qquad - \int_{\partial \Omega_Q} \mathbf{v}^\mathrm{T} \left( \overline{\mathbf{c}}_1 \overline{\mathbf{s}}_4^{e \mathrm{T}} \overline{\mathbf{s}}_4^{e} \boldsymbol{\mu}_{,1} + \overline{\mathbf{c}}_2 \overline{\mathbf{s}}_4^{e \mathrm{T}} \overline{\mathbf{s}}_4^{e} \boldsymbol{\mu}_{,2}  \right) \mathrm{d} \Gamma + \int_{\partial \Omega_R} \hat{\boldsymbol{\varepsilon}}^\mathrm{T} \left( \mathbf{v} \right)  \overline{\mathbf{n}}_4^{e \mathrm{T}} \overline{\mathbf{n}}_4^{e} \overline{\mathbf{n}}_3^{e} \boldsymbol{\mu} \mathrm{d} \Gamma \\
& \qquad  + \int_{\partial \Omega_u} \frac{\eta_{11}}{h_e} \mathbf{v}^\mathrm{T} \mathbf{u} \mathrm{d} \Gamma + \int_{\Gamma_h} \frac{\eta_{21}}{h_e} \left[ \! \left[ \mathbf{v}^\mathrm{T} \right] \! \right] \left[ \!\left[ \mathbf{u} \right] \! \right] \mathrm{d} \Gamma  + \sum_{C_h} \left\{ \mathbf{v}^\mathrm{T} \right\} \overline{\mathbf{s}}_4^e \overline{\mathbf{n}}_3^e \left[ \! \left[ \boldsymbol{\mu} \right] \! \right]
\end{split} \\
\begin{split} \nonumber
& \qquad  = \int_\Omega \mathbf{v}^\mathrm{T} \mathbf{b} \mathrm{d} \Omega + \int_{\partial \Omega_Q} \mathbf{v}^\mathrm{T} \widetilde{\mathbf{Q}} \mathrm{d} \Gamma - \int_\Omega \hat{\boldsymbol{\varepsilon}}^\mathrm{T} \left( \mathbf{v} \right) \widehat{\mathbf{D}} \left( \mathbf{u}, \phi, \boldsymbol{\mu} \right) \mathbf{E} \left( \phi \right) \mathrm{d} \Omega - \int_\Omega \hat{\boldsymbol{\varepsilon}}^\mathrm{T} \left( \mathbf{v} \right) \widehat{\mathbf{Q}} \boldsymbol{\Phi}^{-1} \mathbf{Q} \mathrm{d} \Omega \\
& \qquad + \int_\Omega \hat{\boldsymbol{\varepsilon}}^\mathrm{T} \left( \mathbf{v} \right) \widehat{\mathbf{Q}} \boldsymbol{\Phi}^{-1} \mathbf{b}^\mathrm{T} \boldsymbol{\varepsilon} \left( \mathbf{u} \right) \mathrm{d} \Omega  - \int_{\partial \Omega_\phi} \boldsymbol{\varepsilon}^\mathrm{T} \left( \mathbf{v} \right) \mathbf{e} \mathbf{n}^e \widetilde{\phi} \mathrm{d} \Gamma + \int_{\partial \Omega_R} \hat{\boldsymbol{\varepsilon}}^\mathrm{T} \left( \mathbf{v} \right) \overline{\mathbf{n}}_4^{e \mathrm{T}} \widetilde{\mathbf{R}} \mathrm{d} \Gamma  \\
& \qquad - \int_{\partial \Omega_u} \boldsymbol{\varepsilon}^\mathrm{T} \left( \mathbf{v} \right) \left(\left( \mathbf{C}_{\sigma \varepsilon} + \mathbf{b} \boldsymbol{\Phi}^{-1} \mathbf{b}^\mathrm{T} \right) \overline{\mathbf{n}}_1^{e \mathrm{T}} - \mathbf{b} \boldsymbol{\Phi}^{-1} \widehat{\mathbf{Q}}^\mathrm{T} \overline{\mathbf{n}}_4^{e \mathrm{T}} \right) \widetilde{\mathbf{u}} \mathrm{d} \Gamma + \int_{\partial \Omega_u} \frac{\eta_{11}}{h_e} \mathbf{v}^\mathrm{T} \widetilde{\mathbf{u}} \mathrm{d} \Gamma,
\end{split} \\~\nonumber \\
\begin{split}
\label{eqn:ele_full_mixed}
& - \int_\Omega \mathbf{E}^\mathrm{T} \left( \tau \right) \left( \boldsymbol{\Lambda} + \mathbf{a}^\mathrm{T} \mathbf{C}_{\mu \kappa}^{-1} \mathbf{a} \right) \mathbf{E} \left( \phi \right) \mathrm{d} \Omega - \int_\Omega \mathbf{E}^\mathrm{T} \left( \tau \right) \mathbf{e}^\mathrm{T} \boldsymbol{\varepsilon} \left( \mathbf{u} \right) \mathrm{d} \Omega \\
& \qquad - \int_\Omega \mathbf{E}^\mathrm{T} \left( \tau \right) \mathbf{a}^\mathrm{T} \mathbf{C}_{\mu \kappa}^{-1} \boldsymbol{\mu} \mathrm{d} \Omega + \int_\Omega \mathbf{E}^\mathrm{T} \left( \tau \right)  \overline{\mathbf{n}}_5 \mathbf{Q}_{,1} \mathrm{d} \Omega + \int_\Omega \mathbf{E}^\mathrm{T} \left( \tau \right)  \overline{\mathbf{n}}_6 \mathbf{Q}_{,2} \mathrm{d} \Omega \\
& \qquad - \int_{\Gamma_h \cup \partial \Omega_\phi} \left[ \! \left[ \tau \right] \! \right] \mathbf{n}^{e \mathrm{T}} \left( \boldsymbol{\Lambda} + \mathbf{a}^\mathrm{T} \mathbf{C}_{\mu \kappa}^{-1} \mathbf{a} \right) \left\{ \mathbf{E} \left( \phi \right) \right\} \mathrm{d} \Gamma - \int_{\Gamma_h \cup \partial \Omega_\phi} \left[ \! \left[ \tau \right] \! \right] \mathbf{n}^{e \mathrm{T}} \mathbf{e}^\mathrm{T} \left\{ \boldsymbol{\varepsilon} \left( \mathbf{u} \right) \right\} \mathrm{d} \Gamma \\
& \qquad - \int_{\Gamma_h \cup \partial \Omega_\phi} \left\{ \mathbf{E}^\mathrm{T} \left( \tau \right) \right\} \left( \boldsymbol{\Lambda} + \mathbf{a}^\mathrm{T} \mathbf{C}_{\mu \kappa}^{-1} \mathbf{a} \right) \mathbf{n}^e \left[ \! \left[ \phi \right] \! \right] \mathrm{d} \Gamma - \int_{\Gamma_h \cup \partial \Omega_\phi} \left[ \! \left[ \tau \right] \! \right] \mathbf{n}^{e \mathrm{T}} \mathbf{a}^\mathrm{T} \mathbf{C}_{\mu \kappa}^{-1} \left\{ \boldsymbol{\mu} \right\} \mathrm{d} \Gamma \\
& \qquad + \int_{\Gamma_h \cup \partial \Omega_\phi} \left[ \! \left[ \tau \right] \! \right] \overline{\mathbf{n}}_{51}^e \left\{ \mathbf{Q}_{,1} \right\} \mathrm{d} \Gamma + \int_{\Gamma_h \cup \partial \Omega_\phi} \left[ \! \left[ \tau \right] \! \right] \overline{\mathbf{n}}_{52}^e \left\{ \mathbf{Q}_{,2} \right\} \mathrm{d} \Gamma \\
& \qquad + \int_{\Gamma_h \cup \partial \Omega_u} \left\{ \mathbf{E}^\mathrm{T} \left( \tau \right) \left( \mathbf{e}^\mathrm{T} \overline{\mathbf{n}}_1^{e \mathrm{T}} - \widehat{\mathbf{D}}^\mathrm{T} \left( \mathbf{u}, \phi, \boldsymbol{\mu} \right) \overline{\mathbf{n}}_4^{e \mathrm{T}} \right) \right\} \left[ \! \left[ \mathbf{u} \right] \! \right] \mathrm{d} \Gamma \\
& \qquad - \int_{\Gamma_h} \left\{ \mathbf{E}^\mathrm{T} \left( \tau \right) \right\} \overline{\mathbf{n}}_6^e \left[ \! \left[ \mathbf{Q} \right] \! \right] \mathrm{d} \Gamma  + \int_{\Gamma_h} \frac{\eta_{23}}{h_e} \left[ \! \left[ \tau \right] \! \right] \left[ \!\left[ \phi \right] \! \right] \mathrm{d} \Gamma - \int_{\partial \Omega_Z} \mathbf{E}^\mathrm{T} \left( \tau \right) \mathbf{n}^e \mathbf{n}^{e \mathrm{T}} \overline{\mathbf{n}}_6^e \mathbf{Q} \mathrm{d} \Gamma \\
& \qquad - \int_{\partial \Omega_\omega} \tau \left( \overline{\mathbf{c}}_3 \mathbf{s}^e \mathbf{s}^{e \mathrm{T}} \overline{\mathbf{n}}_6^e \mathbf{Q}_{,1} + \overline{\mathbf{c}}_4 \mathbf{s}^e \mathbf{s}^{e \mathrm{T}} \overline{\mathbf{n}}_6^e \mathbf{Q}_{,2} \right) \mathrm{d} \Gamma + \int_{\partial \Omega_\phi} \frac{\eta_{13}}{h_e} \tau \phi \mathrm{d} \Gamma  + \sum_{C_h} \left\{ \tau \right\} \mathbf{s}^{e \mathrm{T}} \overline{\mathbf{n}}_6^e \left[ \! \left[ \mathbf{Q} \right] \! \right]
\end{split} \\
\begin{split} \nonumber
& \qquad = - \int_\Omega \tau q \mathrm{d} \Omega - \int_{\partial \Omega_\omega} \tau \widetilde{\omega} \mathrm{d} \Gamma + \int_{\partial \Omega_u} \mathbf{E}^\mathrm{T} \left( \tau \right) \left( \mathbf{e}^\mathrm{T} \overline{\mathbf{n}}_1^{e \mathrm{T}} - \widehat{\mathbf{D}}^\mathrm{T} \left( \mathbf{u}, \phi, \boldsymbol{\mu} \right) \overline{\mathbf{n}}_4^{e \mathrm{T}} \right) \widetilde{\mathbf{u}} \mathrm{d} \Gamma   \\
& \qquad  - \int_{\partial \Omega_Z} \mathbf{E}^\mathrm{T} \left( \tau \right) \mathbf{n}^e \widetilde{Z} \mathrm{d} \Gamma - \int_{\partial \Omega_\phi} \mathbf{E}^\mathrm{T} \left( \tau \right) \left( \boldsymbol{\Lambda} + \mathbf{a}^\mathrm{T} \mathbf{C}_{\mu \kappa}^{-1} \mathbf{a} \right) \mathbf{n}^e \widetilde{\phi} \mathrm{d} \Gamma + \int_{\partial \Omega_\phi} \frac{\eta_{13}}{h_e} \tau \widetilde{\phi} \mathrm{d} \Gamma,
\end{split} \\~\nonumber \\
\begin{split}
\label{eqn:mu_full_mixed}
& -\int_\Omega \mathbf{w}^\mathrm{T} \mathbf{C}_{\mu \kappa}^{-1} \boldsymbol{\mu} \mathrm{d} \Omega - \int_\Omega \mathbf{w}^\mathrm{T} \mathbf{C}_{\mu \kappa}^{-1} \mathbf{a} \mathbf{E} \left( \phi \right) \mathrm{d} \Omega - \int_\Omega \mathbf{w}_{,1}^{\mathrm{T}} \overline{\mathbf{c}}_9 \hat{\boldsymbol{\varepsilon}} \left( \mathbf{u} \right) \mathrm{d} \Omega - \int_\Omega \mathbf{w}_{,2}^{\mathrm{T}} \overline{\mathbf{c}}_{10} \hat{\boldsymbol{\varepsilon}} \left( \mathbf{u} \right) \mathrm{d} \Omega \\
& \qquad + \int_{\Gamma_h \cup \partial \Omega_u} \left\{ \mathbf{w}_{,1}^\mathrm{T} \overline{\mathbf{n}}_{21}^{e \mathrm{T}} + \mathbf{w}_{,2}^\mathrm{T} \overline{\mathbf{n}}_{22}^{e \mathrm{T}} \right\} \left[ \! \left[ \mathbf{u} \right] \! \right] \mathrm{d} \Gamma  - \int_{\Gamma_h \cup \partial \Omega_\phi} \left\{ \mathbf{w}^\mathrm{T} \right\} \mathbf{C}_{\mu \kappa}^{-1} \mathbf{a} \mathbf{n}^e \left[ \! \left[ \phi \right] \! \right] \mathrm{d} \Gamma  \\
& \qquad + \int_{\Gamma_h} \left[ \! \left[ \mathbf{w}^\mathrm{T} \right] \! \right] \overline{\mathbf{n}}_3^{e \mathrm{T}} \left\{  \hat{\boldsymbol{\varepsilon}} \left( \mathbf{u} \right)  \right\} \mathrm{d} \Gamma + \int_{\partial \Omega_R} \mathbf{w}^\mathrm{T} \overline{\mathbf{n}}_3^{e \mathrm{T}} \overline{\mathbf{n}}_4^{e \mathrm{T}} \overline{\mathbf{n}}_4^{e} \hat{\boldsymbol{\varepsilon}} \left( \mathbf{u} \right) \mathrm{d} \Gamma + \int_{\partial \Omega_R} \mathbf{w}^\mathrm{T} \overline{\mathbf{n}}_3^{e \mathrm{T}}  \overline{\mathbf{n}}_4^{e \mathrm{T}} \boldsymbol{\lambda}_R \mathrm{d} \Gamma \\
& \qquad - \int_{\partial \Omega_Q} \mathbf{w}_{,1}^\mathrm{T}  \overline{\mathbf{n}}_3^{e \mathrm{T}} \overline{\mathbf{s}}_4^{e \mathrm{T}} \overline{\mathbf{s}}_4^{e} \overline{\mathbf{c}}_1^\mathrm{T} \mathbf{u} \mathrm{d} \Gamma - \int_{\partial \Omega_Q} \mathbf{w}_{,2}^\mathrm{T}  \overline{\mathbf{n}}_3^{e \mathrm{T}} \overline{\mathbf{s}}_4^{e \mathrm{T}} \overline{\mathbf{s}}_4^{e} \overline{\mathbf{c}}_2^\mathrm{T} \mathbf{u} \mathrm{d} \Gamma  + \sum_{C_h} \left[ \! \left[ \mathbf{w}^\mathrm{T} \right] \! \right] \overline{\mathbf{n}}_3^{e \mathrm{T}} \overline{\mathbf{s}}_4^{e \mathrm{T}} \left\{ \mathbf{u} \right\}
\end{split} \\
\begin{split} \nonumber
& \qquad = - \int_{\partial \Omega_d} \mathbf{w}^\mathrm{T} \overline{\mathbf{n}}_3^{e \mathrm{T}} \overline{\mathbf{n}}_4^{e \mathrm{T}} \widetilde{\mathbf{d}} \mathrm{d} \Gamma - \int_{\partial \Omega_\phi} \mathbf{w}^\mathrm{T} \mathbf{C}_{\mu \kappa}^{-1} \mathbf{a} \mathbf{n}^e  \widetilde{\phi} \mathrm{d} \Gamma  \\
& \qquad + \int_{\partial \Omega_u} \mathbf{w}_{,1}^\mathrm{T} \left( \overline{\mathbf{n}}_{21}^{e \mathrm{T}} + \overline{\mathbf{n}}_3^{e \mathrm{T}} \overline{\mathbf{s}}_4^{e \mathrm{T}} \overline{\mathbf{s}}_4^{e} \overline{\mathbf{c}}_1^\mathrm{T} \right) \widetilde{\mathbf{u}} \mathrm{d} \Gamma + \int_{\partial \Omega_u} \mathbf{w}_{,2}^\mathrm{T} \left( \overline{\mathbf{n}}_{22}^{e \mathrm{T}} + \overline{\mathbf{n}}_3^{e \mathrm{T}} \overline{\mathbf{s}}_4^{e \mathrm{T}} \overline{\mathbf{s}}_4^{e} \overline{\mathbf{c}}_2^\mathrm{T} \right) \widetilde{\mathbf{u}} \mathrm{d} \Gamma,
\end{split} \\ ~\nonumber \\
\begin{split}
\label{eqn:Q_full_mixed}
& \int_\Omega \boldsymbol{\xi}^\mathrm{T} \boldsymbol{\Phi}^{-1} \mathbf{Q} \mathrm{d} \Omega - \int_\Omega  \boldsymbol{\xi}^\mathrm{T} \boldsymbol{\Phi}^{-1} \mathbf{b}^\mathrm{T} \boldsymbol{\varepsilon} \left( \mathbf{u} \right) \mathrm{d} \Omega + \int_\Omega \boldsymbol{\xi}_{,1}^\mathrm{T} \overline{\mathbf{c}}_5^\mathrm{T} \mathbf{E} \left( \phi \right) \mathrm{d} \Omega + \int_\Omega \boldsymbol{\xi}_{,2}^\mathrm{T} \overline{\mathbf{c}}_6^\mathrm{T} \mathbf{E} \left( \phi \right) \mathrm{d} \Omega \\
& \qquad + \int_{\Gamma_h \cup \partial \Omega_u} \left\{ \boldsymbol{\xi}^\mathrm{T} \boldsymbol{\Phi}^{-1} \left( \mathbf{b}^\mathrm{T} \overline{\mathbf{n}}_1^{e \mathrm{T}} - \widehat{\mathbf{Q}}^\mathrm{T} \overline{\mathbf{n}}_4^{e \mathbf{T}} \right) \right\} \left[ \! \left[ \mathbf{u} \right] \! \right] \mathrm{d} \Gamma  - \int_{\Gamma_h} \left[ \! \left[ \boldsymbol{\xi}^\mathrm{T} \right] \! \right] \overline{\mathbf{n}}_6^{e \mathrm{T}} \left\{ \mathbf{E} \left( \phi \right) \right\} \mathrm{d} \Gamma \\
& \qquad + \int_{\Gamma_h \cup \partial \Omega_\phi} \left\{ \boldsymbol{\xi}_{,1}^\mathrm{T} \right\} \overline{\mathbf{n}}_{51}^{e \mathrm{T}} \left[ \! \left[ \phi \right] \! \right] \mathrm{d} \Gamma + \int_{\Gamma_h \cup \partial \Omega_\phi} \left\{ \boldsymbol{\xi}_{,2}^\mathrm{T} \right\} \overline{\mathbf{n}}_{52}^{e \mathrm{T}} \left[ \! \left[ \phi \right] \! \right] \mathrm{d} \Gamma - \int_{\partial \Omega_\omega} \boldsymbol{\xi}_{,1}^\mathrm{T}  \overline{\mathbf{n}}_6^{e \mathrm{T}} \mathbf{s}^e \mathbf{s}^{e \mathrm{T}} \overline{\mathbf{c}}_3^\mathrm{T} \phi \mathrm{d} \Gamma \\
& \qquad  - \int_{\partial \Omega_\omega} \boldsymbol{\xi}_{,2}^\mathrm{T}  \overline{\mathbf{n}}_6^{e \mathrm{T}} \mathbf{s}^e \mathbf{s}^{e \mathrm{T}} \overline{\mathbf{c}}_4^\mathrm{T} \phi \mathrm{d} \Gamma - \int_{\partial \Omega_Z} \boldsymbol{\xi}^\mathrm{T}  \overline{\mathbf{n}}_6^{e \mathrm{T}} \mathbf{n}^e \mathbf{n}^{e \mathrm{T}} \mathbf{E} \left( \phi \right) \mathrm{d} \Gamma + \int_{\partial \Omega_Z} \boldsymbol{\xi}^\mathrm{T} \overline{\mathbf{n}}_6^{e \mathrm{T}} \mathbf{n}^e \lambda_Z \mathrm{d} \Gamma \\
& \qquad + \sum_{C_h} \left[ \! \left[ \boldsymbol{\xi}^\mathrm{T} \right] \! \right] \overline{\mathbf{n}}_6^{e \mathrm{T}} \mathbf{s}^e \left\{ \phi \right\} 
\end{split} \\
\begin{split} \nonumber
& \qquad = - \int_{\partial \Omega_P} \boldsymbol{\xi}^\mathrm{T} \overline{\mathbf{n}}_6^{e \mathrm{T}} \mathbf{n}^e \widetilde{P} \mathrm{d} \Gamma + \int_{\partial \Omega_u} \boldsymbol{\xi}^\mathrm{T} \boldsymbol{\Phi}^{-1} \left( \mathbf{b}^\mathrm{T} \overline{\mathbf{n}}_1^{e \mathrm{T}} - \widehat{\mathbf{Q}}^\mathrm{T} \overline{\mathbf{n}}_4^{e \mathbf{T}} \right) \widetilde{\mathbf{u}} \mathrm{d} \Gamma \\
& \qquad + \int_{\partial \Omega_\phi} \boldsymbol{\xi}_{,1}^\mathrm{T} \left( \overline{\mathbf{n}}_{51}^{e \mathrm{T}} + \overline{\mathbf{n}}_6^{e \mathrm{T}} \mathbf{s}^e \mathbf{s}^{e \mathrm{T}} \overline{\mathbf{c}}_3^\mathrm{T} \right) \widetilde{\phi} \mathrm{d} \Gamma + \int_{\partial \Omega_\phi} \boldsymbol{\xi}_{,2}^\mathrm{T} \left( \overline{\mathbf{n}}_{52}^{e \mathrm{T}} + \overline{\mathbf{n}}_6^{e \mathrm{T}} \mathbf{s}^e \mathbf{s}^{e \mathrm{T}} \overline{\mathbf{c}}_4^\mathrm{T} \right) \widetilde{\phi} \mathrm{d} \Gamma,
\end{split} \\ ~\nonumber \\
\begin{split}
& \int_{\partial \Omega_R} \hat{\boldsymbol{\lambda}}_R^\mathrm{T} \overline{\mathbf{n}}_4^e \overline{\mathbf{n}}_3^e \boldsymbol{\mu} \mathrm{d} \Gamma = \int_{\partial \Omega_R} \hat{\boldsymbol{\lambda}}_R^\mathrm{T} \widetilde{\mathbf{R}} \mathrm{d} \Gamma,
\end{split} \\ ~\nonumber \\
\begin{split} \label{eqn:Z_full_mixed}
& \int_{\partial \Omega_Z} \hat{\lambda}_Z \mathbf{n}^{e \mathrm{T}} \overline{\mathbf{n}}_6^e \mathbf{Q} \mathrm{d} \Gamma = \int_{\partial \Omega_Z} \hat{\lambda}_Z \widetilde{Z} \mathrm{d} \Gamma.
\end{split}
\end{align}
Note that when $\mathbf{C}_{\mu \kappa} = \mathbf{0}$ or $\boldsymbol{\Phi} = \mathbf{0}$, $\mathbf{C}_{\mu \kappa}^{-1}$ or $\boldsymbol{\Phi}^{-1}$ can be considered as an infinite large value, and the above equations can be satisfied by letting the corresponding coefficients of $\mathbf{C}_{\mu \kappa}^{-1}$ or $\boldsymbol{\Phi}^{-1}$ to be zero. For example, when $\mathbf{C}_{\mu \kappa} = \mathbf{0}$, the linear coefficient of $\mathbf{C}_{\mu \kappa}^{-1}$ in Eqn.~\ref{eqn:mu_full_mixed} equals zero:
\begin{align}
\begin{split} \label{eqn:cmk0}
& -\int_\Omega \mathbf{w}^\mathrm{T} \boldsymbol{\mu} \mathrm{d} \Omega - \int_\Omega \mathbf{w}^\mathrm{T} \mathbf{a} \mathbf{E} \left( \phi \right) \mathrm{d} \Omega - \int_{\Gamma_h \cup \partial \Omega_\phi} \left\{ \mathbf{w}^\mathrm{T} \right\} \mathbf{a} \mathbf{n}^e \left[ \! \left[ \phi \right] \! \right] \mathrm{d} \Gamma  = - \int_{\partial \Omega_\phi} \mathbf{w}^\mathrm{T} \mathbf{a} \mathbf{n}^e  \widetilde{\phi} \mathrm{d} \Gamma,
\end{split}
\end{align}
The terms containing $\mathbf{C}_{\mu \kappa}^{-1}$ in Eqn.~\ref{eqn:ele_full_mixed} are balanced as long as Eqn.~\ref{eqn:cmk0} is satisfied. Therefore, they can be eliminated. Note that when $\mathbf{C}_{\mu \kappa} = \mathbf{0}$, the strain gradient $\boldsymbol{\kappa}$ can not be achieved directly from the constitutive relation. Therefore, a post-processing algorithm for $\boldsymbol{\kappa}$ is required based on the constant terms in Eqn.\ref{eqn:mu_full_mixed}:
\begin{align}
\begin{split}
& -\int_\Omega \mathbf{w}^\mathrm{T} \boldsymbol{\kappa} \mathrm{d} \Omega - \int_\Omega \mathbf{w}_{,1}^{\mathrm{T}} \overline{\mathbf{c}}_9 \hat{\boldsymbol{\varepsilon}} \left( \mathbf{u} \right) \mathrm{d} \Omega - \int_\Omega \mathbf{w}_{,2}^{\mathrm{T}} \overline{\mathbf{c}}_{10} \hat{\boldsymbol{\varepsilon}} \left( \mathbf{u} \right) \mathrm{d} \Omega + \int_{\Gamma_h} \left[ \! \left[ \mathbf{w}^\mathrm{T} \right] \! \right] \overline{\mathbf{n}}_3^{e \mathrm{T}} \left\{  \hat{\boldsymbol{\varepsilon}} \left( \mathbf{u} \right)  \right\} \mathrm{d} \Gamma \\
& \qquad + \int_{\Gamma_h \cup \partial \Omega_u} \left\{ \mathbf{w}_{,1}^\mathrm{T} \overline{\mathbf{n}}_{21}^{e \mathrm{T}} + \mathbf{w}_{,2}^\mathrm{T} \overline{\mathbf{n}}_{22}^{e \mathrm{T}} \right\} \left[ \! \left[ \mathbf{u} \right] \! \right] \mathrm{d} \Gamma + \int_{\partial \Omega_R} \mathbf{w}^\mathrm{T} \overline{\mathbf{n}}_3^{e \mathrm{T}} \overline{\mathbf{n}}_4^{e \mathrm{T}} \overline{\mathbf{n}}_4^{e} \hat{\boldsymbol{\varepsilon}} \left( \mathbf{u} \right) \mathrm{d} \Gamma \\
& \qquad - \int_{\partial \Omega_Q} \mathbf{w}_{,1}^\mathrm{T}  \overline{\mathbf{n}}_3^{e \mathrm{T}} \overline{\mathbf{s}}_4^{e \mathrm{T}} \overline{\mathbf{s}}_4^{e} \overline{\mathbf{c}}_1^\mathrm{T} \mathbf{u} \mathrm{d} \Gamma - \int_{\partial \Omega_Q} \mathbf{w}_{,2}^\mathrm{T}  \overline{\mathbf{n}}_3^{e \mathrm{T}} \overline{\mathbf{s}}_4^{e \mathrm{T}} \overline{\mathbf{s}}_4^{e} \overline{\mathbf{c}}_2^\mathrm{T} \mathbf{u} \mathrm{d} \Gamma  + \sum_{C_h} \left[ \! \left[ \mathbf{w}^\mathrm{T} \right] \! \right] \overline{\mathbf{n}}_3^{e \mathrm{T}} \overline{\mathbf{s}}_4^{e \mathrm{T}} \left\{ \mathbf{u} \right\} \\
& \qquad + \int_{\partial \Omega_R} \mathbf{w}^\mathrm{T} \overline{\mathbf{n}}_3^{e \mathrm{T}}  \overline{\mathbf{n}}_4^{e \mathrm{T}} \boldsymbol{\lambda}_R \mathrm{d} \Gamma  = - \int_{\partial \Omega_d} \mathbf{w}^\mathrm{T} \overline{\mathbf{n}}_3^{e \mathrm{T}} \overline{\mathbf{n}}_4^{e \mathrm{T}} \widetilde{\mathbf{d}} \mathrm{d} \Gamma  \\
& \qquad + \int_{\partial \Omega_u} \mathbf{w}_{,1}^\mathrm{T} \left( \overline{\mathbf{n}}_{21}^{e \mathrm{T}} + \overline{\mathbf{n}}_3^{e \mathrm{T}} \overline{\mathbf{s}}_4^{e \mathrm{T}} \overline{\mathbf{s}}_4^{e} \overline{\mathbf{c}}_1^\mathrm{T} \right) \widetilde{\mathbf{u}} \mathrm{d} \Gamma + \int_{\partial \Omega_u} \mathbf{w}_{,2}^\mathrm{T} \left( \overline{\mathbf{n}}_{22}^{e \mathrm{T}} + \overline{\mathbf{n}}_3^{e \mathrm{T}} \overline{\mathbf{s}}_4^{e \mathrm{T}} \overline{\mathbf{s}}_4^{e} \overline{\mathbf{c}}_2^\mathrm{T} \right) \widetilde{\mathbf{u}} \mathrm{d} \Gamma.
\end{split}
\end{align}

Combining Eqn.~\ref{eqn:dis_full_mixed} -- \ref{eqn:Z_full_mixed}, we get the algebraic equations in matrix form in terms of the global nodal vectors ($\overline{\mathbf{u}}$, $\overline{\boldsymbol{\phi}}$, $\overline{\boldsymbol{\mu}}$, $\overline{\mathbf{Q}}$, $\overline{\boldsymbol{\lambda}}_R$ and $\overline{\boldsymbol{\lambda}}_Z$).
\begin{align}
\begin{split} \label{eqn:full_mixed}
\left[ \begin{matrix} \mathbf{K}_{uu}^M & \mathbf{K}_{u \phi}^M & \mathbf{K}_{u \mu}^M & \mathbf{K}_{uQ}^M & \mathbf{0} & \mathbf{0} \\
\left( \mathbf{K}_{u \phi}^M \right)^\mathrm{T} & \mathbf{K}_{\phi \phi}^M & \mathbf{K}_{\phi \mu}^M & \mathbf{K}_{\phi Q}^M & \mathbf{0} & \mathbf{0} \\
\left( \mathbf{K}_{u \mu}^M \right)^\mathrm{T} & \left( \mathbf{K}_{\phi \mu}^M \right)^\mathrm{T} & \mathbf{K}_{\mu \mu}^M & \mathbf{0} & \mathbf{K}_{\mu R}^M & \mathbf{0} \\
\left( \mathbf{K}_{uQ}^M \right)^\mathrm{T} & \left( \mathbf{K}_{\phi Q}^M \right)^\mathrm{T} & \mathbf{0} &  \mathbf{K}_{QQ}^M & \mathbf{0} & \mathbf{K}_{QZ}^M \\
\mathbf{0} & \mathbf{0} & \left( \mathbf{K}_{\mu R}^M \right)^\mathrm{T} & \mathbf{0} & \mathbf{0} & \mathbf{0} \\
\mathbf{0} & \mathbf{0} & \mathbf{0} & \left( \mathbf{K}_{QZ}^M \right)^\mathrm{T} & \mathbf{0} & \mathbf{0} \end{matrix} \right] \left[ \begin{matrix} \overline{\mathbf{u}} \\ \overline{\boldsymbol{\phi}} \\ \overline{\boldsymbol{\mu}} \\  \overline{\mathbf{Q}} \\ \overline{\boldsymbol{\lambda}}_R \\ \overline{\boldsymbol{\lambda}}_Z \end{matrix} \right] = \left[ \begin{matrix} \mathbf{f}_u^M \\ \mathbf{f}_\phi^M \\ \mathbf{f}_\mu^M \\ \mathbf{f}_Q^M \\ \mathbf{f}_R^M \\ \mathbf{f}_Z^M \end{matrix} \right].
\end{split}
\end{align}

\subsubsection{Reduced theory}

In the reduced theory, the high-order electric displacement $\mathbf{Q}$ vanishes. The displacement $\mathbf{u}$, electric potential $\phi$, high-order stress $\boldsymbol{\mu}$ are still employed as independent variables at each point and interpolated by linear trial functions. Thus the basic mixed FPM formula for the reduced flexoelectric theory can be written as: 
\begin{align}
\begin{split}
\nonumber
& \int_\Omega \boldsymbol{\varepsilon}^\mathrm{T} \left( \mathbf{v} \right) \mathbf{C}_{\sigma \varepsilon} \boldsymbol{\varepsilon} \left( \mathbf{u} \right) \mathrm{d} \Omega - \int_\Omega \boldsymbol{\varepsilon}^\mathrm{T} \left( \mathbf{v} \right) \mathbf{e} \mathbf{E} \left( \phi \right) \mathrm{d} \Omega  - \int_\Omega \hat{\boldsymbol{\varepsilon}}^\mathrm{T} \left( \mathbf{v} \right) \overline{\mathbf{c}}_9^\mathrm{T} \boldsymbol{\mu}_{,1} \mathrm{d} \Omega - \int_\Omega \hat{\boldsymbol{\varepsilon}}^\mathrm{T} \left( \mathbf{v} \right) \overline{\mathbf{c}}_{10}^\mathrm{T} \boldsymbol{\mu}_{,2} \mathrm{d} \Omega \\
& \qquad - \int_{\Gamma_h \cup \partial \Omega_u} \left[ \! \left[ \mathbf{v}^\mathrm{T} \right] \! \right] \overline{\mathbf{n}}_1^e \mathbf{C}_{\sigma \varepsilon}  \left\{ \boldsymbol{\varepsilon} \left( \mathbf{u} \right) \right\} \mathrm{d} \Gamma -\int_{\Gamma_h \cup \partial \Omega_u} \left\{ \boldsymbol{\varepsilon}^\mathrm{T} \left( \mathbf{v} \right) \right\} \mathbf{C}_{\sigma \varepsilon} \overline{\mathbf{n}}_1^{e \mathrm{T}} \left[ \! \left[ \mathbf{u} \right] \! \right] \mathrm{d} \Gamma \\
& \qquad + \int_{\Gamma_h \cup \partial \Omega_u} \left[ \! \left[ \mathbf{v}^\mathrm{T} \right] \! \right] \overline{\mathbf{n}}_1^e \mathbf{e}  \left\{ \mathbf{E} \left( \phi \right) \right\}  \mathrm{d} \Gamma  + \int_{\Gamma_h \cup \partial \Omega_u} \left[ \! \left[ \mathbf{v}^\mathrm{T} \right] \! \right] \left\{ \overline{\mathbf{n}}_{21}^e \boldsymbol{\mu}_{,1} + \overline{\mathbf{n}}_{22}^e \boldsymbol{\mu}_{,2} \right\} \mathrm{d} \Gamma \\
& \qquad  - \int_{\Gamma_h \cup \partial \Omega_\phi} \left\{ \boldsymbol{\varepsilon}^\mathrm{T} \left( \mathbf{v} \right) \right\} \mathbf{e} \mathbf{n}^e \left[ \! \left[ \phi \right] \! \right] \mathrm{d} \Gamma + \int_{\Gamma_h} \left\{ \hat{\boldsymbol{\varepsilon}}^\mathrm{T} \left( \mathbf{v} \right) \right\} \overline{\mathbf{n}}_3^e \left[ \! \left[ \boldsymbol{\mu} \right] \! \right] \mathrm{d} \Gamma  \\
& \qquad - \int_{\partial \Omega_Q} \mathbf{v}^\mathrm{T} \left( \overline{\mathbf{c}}_1 \overline{\mathbf{s}}_4^{e \mathrm{T}} \overline{\mathbf{s}}_4^{e} \boldsymbol{\mu}_{,1} + \overline{\mathbf{c}}_2 \overline{\mathbf{s}}_4^{e \mathrm{T}} \overline{\mathbf{s}}_4^{e} \boldsymbol{\mu}_{,2}  \right) \mathrm{d} \Gamma + \int_{\partial \Omega_R} \hat{\boldsymbol{\varepsilon}}^\mathrm{T} \left( \mathbf{v} \right)  \overline{\mathbf{n}}_4^{e \mathrm{T}} \overline{\mathbf{n}}_4^{e} \overline{\mathbf{n}}_3^{e} \boldsymbol{\mu} \mathrm{d} \Gamma \\
& \qquad  + \int_{\partial \Omega_u} \frac{\eta_{11}}{h_e} \mathbf{v}^\mathrm{T} \mathbf{u} \mathrm{d} \Gamma + \int_{\Gamma_h} \frac{\eta_{21}}{h_e} \left[ \! \left[ \mathbf{v}^\mathrm{T} \right] \! \right] \left[ \!\left[ \mathbf{u} \right] \! \right] \mathrm{d} \Gamma  + \sum_{C_h} \left\{ \mathbf{v}^\mathrm{T} \right\} \overline{\mathbf{s}}_4^e \overline{\mathbf{n}}_3^e \left[ \! \left[ \boldsymbol{\mu} \right] \! \right]
\end{split} \\
\begin{split} \label{eqn:dis_red_mixed}
& \qquad  = \int_\Omega \mathbf{v}^\mathrm{T} \mathbf{b} \mathrm{d} \Omega + \int_{\partial \Omega_Q} \mathbf{v}^\mathrm{T} \widetilde{\mathbf{Q}} \mathrm{d} \Gamma  - \int_{\partial \Omega_\phi} \boldsymbol{\varepsilon}^\mathrm{T} \left( \mathbf{v} \right) \mathbf{e} \mathbf{n}^e \widetilde{\phi} \mathrm{d} \Gamma + \int_{\partial \Omega_R} \hat{\boldsymbol{\varepsilon}}^\mathrm{T} \left( \mathbf{v} \right) \overline{\mathbf{n}}_4^{e \mathrm{T}} \widetilde{\mathbf{R}} \mathrm{d} \Gamma  \\
& \qquad - \int_{\partial \Omega_u} \boldsymbol{\varepsilon}^\mathrm{T} \left( \mathbf{v} \right) \mathbf{C}_{\sigma \varepsilon} \overline{\mathbf{n}}_1^{e \mathrm{T}} \widetilde{\mathbf{u}} \mathrm{d} \Gamma + \int_{\partial \Omega_u} \frac{\eta_{11}}{h_e} \mathbf{v}^\mathrm{T} \widetilde{\mathbf{u}} \mathrm{d} \Gamma,
\end{split} \\~\nonumber \\
\begin{split}
\label{eqn:ele_red_mixed}
& - \int_\Omega \mathbf{E}^\mathrm{T} \left( \tau \right) \left( \boldsymbol{\Lambda} + \mathbf{a}^\mathrm{T} \mathbf{C}_{\mu \kappa}^{-1} \mathbf{a} \right) \mathbf{E} \left( \phi \right) \mathrm{d} \Omega - \int_\Omega \mathbf{E}^\mathrm{T} \left( \tau \right) \mathbf{e}^\mathrm{T} \boldsymbol{\varepsilon} \left( \mathbf{u} \right) \mathrm{d} \Omega \\
& \qquad - \int_\Omega \mathbf{E}^\mathrm{T} \left( \tau \right) \mathbf{a}^\mathrm{T} \mathbf{C}_{\mu \kappa}^{-1} \boldsymbol{\mu} \mathrm{d} \Omega \\
& \qquad - \int_{\Gamma_h \cup \partial \Omega_\phi} \left[ \! \left[ \tau \right] \! \right] \mathbf{n}^{e \mathrm{T}} \left( \boldsymbol{\Lambda} + \mathbf{a}^\mathrm{T} \mathbf{C}_{\mu \kappa}^{-1} \mathbf{a} \right) \left\{ \mathbf{E} \left( \phi \right) \right\} \mathrm{d} \Gamma - \int_{\Gamma_h \cup \partial \Omega_\phi} \left[ \! \left[ \tau \right] \! \right] \mathbf{n}^{e \mathrm{T}} \mathbf{e}^\mathrm{T} \left\{ \boldsymbol{\varepsilon} \left( \mathbf{u} \right) \right\} \mathrm{d} \Gamma \\
& \qquad - \int_{\Gamma_h \cup \partial \Omega_\phi} \left\{ \mathbf{E}^\mathrm{T} \left( \tau \right) \right\} \left( \boldsymbol{\Lambda} + \mathbf{a}^\mathrm{T} \mathbf{C}_{\mu \kappa}^{-1} \mathbf{a} \right) \mathbf{n}^e \left[ \! \left[ \phi \right] \! \right] \mathrm{d} \Gamma - \int_{\Gamma_h \cup \partial \Omega_\phi} \left[ \! \left[ \tau \right] \! \right] \mathbf{n}^{e \mathrm{T}} \mathbf{a}^\mathrm{T} \mathbf{C}_{\mu \kappa}^{-1} \left\{ \boldsymbol{\mu} \right\} \mathrm{d} \Gamma \\
& \qquad + \int_{\Gamma_h \cup \partial \Omega_u} \left\{ \mathbf{E}^\mathrm{T} \left( \tau \right) \right\} \mathbf{e}^\mathrm{T} \overline{\mathbf{n}}_1^{e \mathrm{T}} \left[ \! \left[ \mathbf{u} \right] \! \right] \mathrm{d} \Gamma + \int_{\partial \Omega_\phi} \frac{\eta_{13}}{h_e} \tau \phi \mathrm{d} \Gamma + \int_{\Gamma_h} \frac{\eta_{23}}{h_e} \left[ \! \left[ \tau \right] \! \right] \left[ \!\left[ \phi \right] \! \right] \mathrm{d} \Gamma 
\end{split} \\
\begin{split} \nonumber
& \qquad = - \int_\Omega \tau q \mathrm{d} \Omega - \int_{\partial \Omega_\omega} \tau \widetilde{\omega} \mathrm{d} \Gamma - \int_{\partial \Omega_\phi} \mathbf{E}^\mathrm{T} \left( \tau \right) \left( \boldsymbol{\Lambda} + \mathbf{a}^\mathrm{T} \mathbf{C}_{\mu \kappa}^{-1} \mathbf{a} \right) \mathbf{n}^e \widetilde{\phi} \mathrm{d} \Gamma \\
& \qquad  + \int_{\partial \Omega_u} \mathbf{E}^\mathrm{T} \left( \tau \right) \mathbf{e}^\mathrm{T} \overline{\mathbf{n}}_1^{e \mathrm{T}} \widetilde{\mathbf{u}} \mathrm{d} \Gamma + \int_{\partial \Omega_\phi} \frac{\eta_{13}}{h_e} \tau \widetilde{\phi} \mathrm{d} \Gamma.
\end{split} \\~\nonumber \\
\begin{split}
\label{eqn:mu_red_mixed}
& -\int_\Omega \mathbf{w}^\mathrm{T} \mathbf{C}_{\mu \kappa}^{-1} \boldsymbol{\mu} \mathrm{d} \Omega - \int_\Omega \mathbf{w}^\mathrm{T} \mathbf{C}_{\mu \kappa}^{-1} \mathbf{a} \mathbf{E} \left( \phi \right) \mathrm{d} \Omega - \int_\Omega \mathbf{w}_{,1}^{\mathrm{T}} \overline{\mathbf{c}}_9 \hat{\boldsymbol{\varepsilon}} \left( \mathbf{u} \right) \mathrm{d} \Omega - \int_\Omega \mathbf{w}_{,2}^{\mathrm{T}} \overline{\mathbf{c}}_{10} \hat{\boldsymbol{\varepsilon}} \left( \mathbf{u} \right) \mathrm{d} \Omega \\
& \qquad + \int_{\Gamma_h \cup \partial \Omega_u} \left\{ \mathbf{w}_{,1}^\mathrm{T} \overline{\mathbf{n}}_{21}^{e \mathrm{T}} + \mathbf{w}_{,2}^\mathrm{T} \overline{\mathbf{n}}_{22}^{e \mathrm{T}} \right\} \left[ \! \left[ \mathbf{u} \right] \! \right] \mathrm{d} \Gamma  - \int_{\Gamma_h \cup \partial \Omega_\phi} \left\{ \mathbf{w}^\mathrm{T} \right\} \mathbf{C}_{\mu \kappa}^{-1} \mathbf{a} \mathbf{n}^e \left[ \! \left[ \phi \right] \! \right] \mathrm{d} \Gamma  \\
& \qquad + \int_{\Gamma_h} \left[ \! \left[ \mathbf{w}^\mathrm{T} \right] \! \right] \overline{\mathbf{n}}_3^{e \mathrm{T}} \left\{  \hat{\boldsymbol{\varepsilon}} \left( \mathbf{u} \right)  \right\} \mathrm{d} \Gamma + \int_{\partial \Omega_R} \mathbf{w}^\mathrm{T} \overline{\mathbf{n}}_3^{e \mathrm{T}} \overline{\mathbf{n}}_4^{e \mathrm{T}} \overline{\mathbf{n}}_4^{e} \hat{\boldsymbol{\varepsilon}} \left( \mathbf{u} \right) \mathrm{d} \Gamma + \int_{\partial \Omega_R} \mathbf{w}^\mathrm{T} \overline{\mathbf{n}}_3^{e \mathrm{T}}  \overline{\mathbf{n}}_4^{e \mathrm{T}} \boldsymbol{\lambda}_R \mathrm{d} \Gamma \\
& \qquad - \int_{\partial \Omega_Q} \mathbf{w}_{,1}^\mathrm{T}  \overline{\mathbf{n}}_3^{e \mathrm{T}} \overline{\mathbf{s}}_4^{e \mathrm{T}} \overline{\mathbf{s}}_4^{e} \overline{\mathbf{c}}_1^\mathrm{T} \mathbf{u} \mathrm{d} \Gamma - \int_{\partial \Omega_Q} \mathbf{w}_{,2}^\mathrm{T}  \overline{\mathbf{n}}_3^{e \mathrm{T}} \overline{\mathbf{s}}_4^{e \mathrm{T}} \overline{\mathbf{s}}_4^{e} \overline{\mathbf{c}}_2^\mathrm{T} \mathbf{u} \mathrm{d} \Gamma  + \sum_{C_h} \left[ \! \left[ \mathbf{w}^\mathrm{T} \right] \! \right] \overline{\mathbf{n}}_3^{e \mathrm{T}} \overline{\mathbf{s}}_4^{e \mathrm{T}} \left\{ \mathbf{u} \right\}
\end{split} \\
\begin{split} \nonumber
& \qquad = - \int_{\partial \Omega_d} \mathbf{w}^\mathrm{T} \overline{\mathbf{n}}_3^{e \mathrm{T}} \overline{\mathbf{n}}_4^{e \mathrm{T}} \widetilde{\mathbf{d}} \mathrm{d} \Gamma - \int_{\partial \Omega_\phi} \mathbf{w}^\mathrm{T} \mathbf{C}_{\mu \kappa}^{-1} \mathbf{a} \mathbf{n}^e  \widetilde{\phi} \mathrm{d} \Gamma  \\
& \qquad + \int_{\partial \Omega_u} \mathbf{w}_{,1}^\mathrm{T} \left( \overline{\mathbf{n}}_{21}^{e \mathrm{T}} + \overline{\mathbf{n}}_3^{e \mathrm{T}} \overline{\mathbf{s}}_4^{e \mathrm{T}} \overline{\mathbf{s}}_4^{e} \overline{\mathbf{c}}_1^\mathrm{T} \right) \widetilde{\mathbf{u}} \mathrm{d} \Gamma + \int_{\partial \Omega_u} \mathbf{w}_{,2}^\mathrm{T} \left( \overline{\mathbf{n}}_{22}^{e \mathrm{T}} + \overline{\mathbf{n}}_3^{e \mathrm{T}} \overline{\mathbf{s}}_4^{e \mathrm{T}} \overline{\mathbf{s}}_4^{e} \overline{\mathbf{c}}_2^\mathrm{T} \right) \widetilde{\mathbf{u}} \mathrm{d} \Gamma,
\end{split} \\ ~\nonumber \\
\begin{split}
& \int_{\partial \Omega_R} \hat{\boldsymbol{\lambda}}_R^\mathrm{T} \overline{\mathbf{n}}_4^e \overline{\mathbf{n}}_3^e \boldsymbol{\mu} \mathrm{d} \Gamma = \int_{\partial \Omega_R} \hat{\boldsymbol{\lambda}}_R^\mathrm{T} \widetilde{\mathbf{R}} \mathrm{d} \Gamma.
\end{split}
\end{align}
Or in matrices form:
\begin{align}
\begin{split} \label{eqn:red_mixed}
\left[ \begin{matrix} \mathbf{K}_{uu}^M & \mathbf{K}_{u \phi}^M & \mathbf{K}_{u \mu}^M  & \mathbf{0} \\
\left( \mathbf{K}_{u \phi}^M \right)^\mathrm{T} & \mathbf{K}_{\phi \phi}^M & \mathbf{K}_{\phi \mu}^M & \mathbf{0} \\
\left( \mathbf{K}_{u \mu}^M \right)^\mathrm{T} & \left( \mathbf{K}_{\phi \mu}^M \right)^\mathrm{T} & \mathbf{K}_{\mu \mu}^M & \mathbf{K}_{\mu R}^M \\
\mathbf{0} & \mathbf{0} & \left( \mathbf{K}_{\mu R}^M \right)^\mathrm{T}  & \mathbf{0} \end{matrix} \right] \left[ \begin{matrix} \overline{\mathbf{u}} \\ \overline{\boldsymbol{\phi}} \\ \overline{\boldsymbol{\mu}} \\  \overline{\boldsymbol{\lambda}}_R \end{matrix} \right] = \left[ \begin{matrix} \mathbf{f}_u^M \\ \mathbf{f}_\phi^M \\ \mathbf{f}_\mu^M \\  \mathbf{f}_R^M \end{matrix} \right].
\end{split}
\end{align}

\subsection{Eliminating high-order variables}

In Eqn.~\ref{eqn:full_mixed}, there are 13 independent variables at each point (3 primal variables which are  the same as in the primal method, and 10 high-order variables) and $2 n_R + n_Z$ Lagrange multipliers. However, the number of unknown variables can be reduced by transforming the high-order variables back to the primal displacement and electric potential at each Fragile Point at the local level. First, we rearrange Eqn.~\ref{eqn:full_mixed} as:
\begin{align}
\begin{split}
\left[ \begin{matrix} \mathbf{K}_{\mu \mu}^M  &  \mathbf{K}_{\mu R}^M \\ \left( \mathbf{K}_{\mu R}^M \right)^\mathrm{T}  & \mathbf{0} \end{matrix} \right] \left[ \begin{matrix} \overline{\boldsymbol{\mu}} \\ \overline{\boldsymbol{\lambda}}_R \end{matrix} \right] & = \left[ \begin{matrix} \mathbf{f}_\mu^M - \left( \mathbf{K}_{u \mu}^M \right)^\mathrm{T} \overline{\mathbf{u}} - \left( \mathbf{K}_{\phi \mu}^M \right)^\mathrm{T} \overline{\boldsymbol{\phi}} \\  \mathbf{f}_R^M\end{matrix} \right],
\end{split} \\
\begin{split}
\left[ \begin{matrix} \mathbf{K}_{QQ}^M  &  \mathbf{K}_{QZ}^M \\ \left( \mathbf{K}_{QZ}^M \right)^\mathrm{T}  & \mathbf{0} \end{matrix} \right] \left[ \begin{matrix}  \overline{\mathbf{Q}} \\ \overline{\boldsymbol{\lambda}}_Z \end{matrix} \right] & = \left[ \begin{matrix} \mathbf{f}_Q^M - \left( \mathbf{K}_{uQ}^M \right)^\mathrm{T} \overline{\mathbf{u}} - \left( \mathbf{K}_{\phi Q}^M \right)^\mathrm{T} \overline{\boldsymbol{\phi}} \\  \mathbf{f}_Z^M\end{matrix} \right].
\end{split}
\end{align}
Thus we can obtain a linear transformation relation between the two related nodal variable sets :
\begin{align}
\begin{split}
\overline{\boldsymbol{\mu}} = & - \mathbf{G}_{\mu} \left( \mathbf{K}_{u \mu}^M \right)^\mathrm{T} \overline{\mathbf{u}} - \mathbf{G}_{\mu} \left( \mathbf{K}_{\phi \mu}^M \right)^\mathrm{T} \overline{\boldsymbol{\phi}} + \mathbf{g}_\mu,
\end{split} \\
\begin{split}
\overline{\mathbf{Q}} = & - \mathbf{G}_{Q} \left( \mathbf{K}_{uQ}^M \right)^\mathrm{T} \overline{\mathbf{u}} - \mathbf{G}_{Q} \left( \mathbf{K}_{\phi Q}^M \right)^\mathrm{T} \overline{\boldsymbol{\phi}} + \mathbf{g}_Q,
\end{split}
\end{align}
where
\begin{align}
\begin{split} \nonumber
\mathbf{G}_{\mu} = & \left( \mathbf{K}_{\mu \mu}^M \right)^{-1} -  \left( \mathbf{K}_{\mu \mu}^M \right)^{-1}  \mathbf{K}_{\mu R}^M \left( \left( \mathbf{K}_{\mu R}^M \right)^\mathrm{T} \left( \mathbf{K}_{\mu \mu}^M \right)^{-1} \mathbf{K}_{\mu R}^M \right)^{-1} \left( \mathbf{K}_{\mu R}^M \right)^\mathrm{T} \left( \mathbf{K}_{\mu \mu}^M \right)^{-1}, \\
\mathbf{G}_{Q} = &  \left( \mathbf{K}_{QQ}^M \right)^{-1} -  \left( \mathbf{K}_{QQ}^M \right)^{-1}  \mathbf{K}_{QZ}^M \left( \left( \mathbf{K}_{QZ}^M \right)^\mathrm{T} \left( \mathbf{K}_{QQ}^M \right)^{-1} \mathbf{K}_{QZ}^M \right)^{-1} \left( \mathbf{K}_{QZ}^M \right)^\mathrm{T} \left( \mathbf{K}_{QQ}^M \right)^{-1}, \\
 \mathbf{g}_\mu = & \mathbf{G}_{\mu} \mathbf{f}_\mu^M + \left( \mathbf{K}_{\mu \mu}^M \right)^{-1}  \mathbf{K}_{\mu R}^M \left( \left( \mathbf{K}_{\mu R}^M \right)^\mathrm{T} \left( \mathbf{K}_{\mu \mu}^M \right)^{-1} \mathbf{K}_{\mu R}^M \right)^{-1} \mathbf{f}_R^M, \\
 \mathbf{g}_Q = & \mathbf{G}_{Q} \mathbf{f}_Q^M + \left( \mathbf{K}_{QQ}^M \right)^{-1}  \mathbf{K}_{QZ}^M \left( \left( \mathbf{K}_{QZ}^M \right)^\mathrm{T} \left( \mathbf{K}_{QQ}^M \right)^{-1} \mathbf{K}_{QZ}^M \right)^{-1} \mathbf{f}_Z^M.
\end{split}
\end{align}

According to Eqn.~\ref{eqn:mu_full_mixed} and \ref{eqn:Q_full_mixed}, when the Fragile Points are placed at the centroid of each subdomain, and only one integration point is used, $ \left( \mathbf{K}_{\mu \mu}^M \right)^{-1}$ and $\left( \mathbf{K}_{QQ}^M \right)^{-1}$ are known banded matrices:
\begin{align}
\begin{split}
\left( \mathbf{K}_{\mu \mu}^M \right)^{-1} = & -  \text{diag} \left( \left[ \frac{1}{A_1}, \frac{1}{A_2}, \cdots, \frac{1}{A_N} \right] \right) \otimes \mathbf{C}_{\mu \kappa}, \\
\left( \mathbf{K}_{QQ}^M \right)^{-1} = & \quad \; \text{diag} \left( \left[ \frac{1}{A_1}, \frac{1}{A_2}, \cdots, \frac{1}{A_N} \right] \right) \otimes \boldsymbol{\Phi},
\end{split}
\end{align}
where $A_i$ denotes the area of each subdomain $E_i$, and $N$ is the number of all the Fragile Points. Therefore, the transforming matrices $\mathbf{G}_{\mu}$ and $\mathbf{G}_{Q}$ can be easily achieved by inverting a $2 n_R \times 2 n_R$ matrix $ \left[ \left( \mathbf{K}_{\mu R}^M \right)^\mathrm{T} \left( \mathbf{K}_{\mu \mu}^M \right)^{-1} \mathbf{K}_{\mu R}^M \right]$ and a $n_Z \times n_Z$ matrix $\left[ \left( \mathbf{K}_{QZ}^M \right)^\mathrm{T} \left( \mathbf{K}_{QQ}^M \right)^{-1} \mathbf{K}_{QZ}^M \right]$ respectively.

We can prove that $\mathbf{G}_{\mu}$ and $\mathbf{G}_{Q}$ are also banded. All the Fragile Points can be divided into two groups: for $P_i \in \widehat{E}^{in}$, the nodal displacement or electric potential value at $P_i$ has no influence on trial functions on the external boundary; while at $P_i \in \widehat{E}^{ex}$, a non-zero nodal value will result in a non-zero trial function on $\partial \Omega$, i.e., $\widehat{E}^{ex} = \underset{\partial E_i \cup \partial \Omega \not= \varnothing} {\bigcup} \widehat{E}^i$, where $\widehat{E}^i $ is the set of all the supporting points of $P_i$. Thus, we have:
\begin{align}
\begin{split}
\boldsymbol{\mu}^i = & \left\{ \begin{matrix} \mathbf{G}_{\mu u}^i \mathbf{u}_{E_2}^i + \mathbf{G}_{\mu \phi}^i \boldsymbol{\phi}_{E_2}^i, & \text{for} P_i \in \widehat{E}^{in} \\ \mathbf{G}_{\mu u}^i \mathbf{u}_{E_2}^{ex} + \mathbf{G}_{\mu \phi}^i \boldsymbol{\phi}_{E_2}^{ex} + \mathbf{g}_\mu^i, & \text{for} P_i \in \widehat{E}^{ex}  \end{matrix} \right. , \\
\mathbf{Q}^i =  & \left\{ \begin{matrix} \mathbf{G}_{Qu}^i \mathbf{u}_{E_2}^i + \mathbf{G}_{Q \phi}^i \boldsymbol{\phi}_{E_2}^i, & \text{for} P_i \in \widehat{E}^{in} \\ \mathbf{G}_{Qu}^i \mathbf{u}_{E_2}^{ex} + \mathbf{G}_{Q \phi}^i \boldsymbol{\phi}_{E_2}^{ex} + \mathbf{g}_Q^i, & \text{for} P_i \in \widehat{E}^{ex}  \end{matrix} \right. ,
\end{split}
\end{align}
where $\boldsymbol{\mu}^i$ and $\mathbf{Q}^i$ are the nodal values of $\boldsymbol{\mu}$ and $\mathbf{Q}$ at $P_i$. $\mathbf{u}_{E_2}^i $ and $\boldsymbol{\phi}_{E_2}^i$ are the nodal displacement and electric potential vectors for all $P_j \in \widehat{E}^i_2$, where $\widehat{E}^i_2 = \underset{P_k \in \widehat{E}^i} {\bigcup} \widehat{E}^k$ is the ``generalized’’ set of supporting points of $P_i$. $\mathbf{u}_{E_2}^{ex} $ and $\boldsymbol{\phi}_{E_2}^{ex}$ are the nodal displacement and electric potential vectors for all $P_j \in \widehat{E}^{ex}_{E_2}$, where $\widehat{E}^{ex}_{E_2} = \underset{P_k \in \widehat{E}^{ex}} {\bigcup} \widehat{E}^k$. $\mathbf{G}_{\mu u}^i$, $\mathbf{G}_{\mu \phi}^i$, $\mathbf{G}_{Qu}^i$, $\mathbf{G}_{Q \phi}^i$, $\mathbf{g}_\mu^i$ and $\mathbf{g}_Q^i$ are linear algebraic matrices (vectors).

Thus, all the high-order variables are eliminated at the point level. Finally, the mixed FPM formula with only nodal displacement and electric potential is achieved:
\begin{align}
\begin{split} \label{eqn:mixed}
\left[ \begin{matrix} \overline{\mathbf{K}}_{uu}^M &  \overline{\mathbf{K}}_{u \phi}^M \\ \left( \overline{\mathbf{K}}_{u \phi}^M \right)^\mathrm{T} &  \overline{\mathbf{K}}_{\phi \phi}^M \end{matrix} \right] \left[ \begin{matrix} \overline{\mathbf{u}} \\ \overline{\boldsymbol{\phi}} \end{matrix} \right] = \left[ \begin{matrix} \overline{\mathbf{f}}_u^M \\ \overline{\mathbf{f}}_\phi^M \end{matrix} \right]
\end{split},\quad \text{or} \; \mathbf{K}^M \overline{\mathbf{x}} = \mathbf{f}^M.
\end{align}
In full flexoelectric theory:
\begin{align}
\begin{split} \nonumber
\overline{\mathbf{K}}_{uu}^M = & \mathbf{K}_{uu}^M - \mathbf{K}_{u \mu}^M \mathbf{G}_\mu \left(  \mathbf{K}_{u \mu}^M \right)^\mathrm{T} -  \mathbf{K}_{uQ}^M \mathbf{G}_Q \left(  \mathbf{K}_{uQ}^M \right)^\mathrm{T}, \\
 \overline{\mathbf{K}}_{u \phi}^M = & \mathbf{K}_{u \phi}^M -  \mathbf{K}_{u \mu}^M \mathbf{G}_\mu \left(  \mathbf{K}_{\phi \mu}^M \right)^\mathrm{T} -  \mathbf{K}_{uQ}^M \mathbf{G}_Q \left(  \mathbf{K}_{\phi Q}^M \right)^\mathrm{T}, \\
 \overline{\mathbf{K}}_{\phi \phi}^M = & \mathbf{K}_{\phi \phi}^M -  \mathbf{K}_{\phi \mu}^M \mathbf{G}_\mu \left(  \mathbf{K}_{\phi \mu}^M \right)^\mathrm{T} -  \mathbf{K}_{\phi Q}^M \mathbf{G}_Q \left(  \mathbf{K}_{\phi Q}^M \right)^\mathrm{T}, \\
 \overline{\mathbf{f}}_u^M = & \mathbf{f}_u^M - \mathbf{K}_{u \mu}^M \mathbf{g}_\mu  - \mathbf{K}_{u Q}^M \mathbf{g}_Q, \\
  \overline{\mathbf{f}}_\phi^M = & \mathbf{f}_\phi^M - \mathbf{K}_{\phi \mu}^M \mathbf{g}_\mu  - \mathbf{K}_{\phi Q}^M \mathbf{g}_Q,
\end{split}
\end{align}
While in the reduced theory, the last terms relating to $\mathbf{Q}$ in the above equations are absent. The same as in the primal FPM, when the full flexoelectric theory is employed, Eqn.~\ref{eqn:mixed} is a nonlinear equation system and can be solved by Newton-Raphson method. The stiffness matrix $\mathbf{K}^M$ is sparse and symmetric for both the theories. Also, there are only three explicit unknown variables at each Fragile Point.

\section{Simulations of Crack Initiation \& Propagation in FPM} \label{sec:Crack}

As a benefit of the discontinuous trial and test functions, it is much simpler in the FPM to simulate crack initiation and propagation,  as compared with other numerical methods. In the present work, all the cracks are considered as mechanically traction-free. Yet in the electrical field, two different crack-face boundary conditions are considered. The first one is known as an electrically impermeable boundary condition, in which:
\begin{align}
\begin{split}
\widetilde{\omega}^{+} = \widetilde{\omega}^{-} = 0, \quad \Delta \phi = \phi^{+} - \phi^{-} \neq 0.
\end{split}
\end{align}
This condition implies that the gap between the upper ($+$) and lower ($-$) crack-faces are filled with a medium with zero electric permittivity. Alternatively, an electrically permeable crack-face boundary condition is defined as:
\begin{align}
\begin{split}
\phi^{+} = \phi^{-},  \quad \widetilde{\omega}^{+} = \widetilde{\omega}^{-},
\end{split}
\end{align}
that is, an infinite electric permittivity of the medium between the upper and lower crack-faces is assumed.

\begin{figure}[htbp] 
  \centering 
  \includegraphics[width=0.8\textwidth]{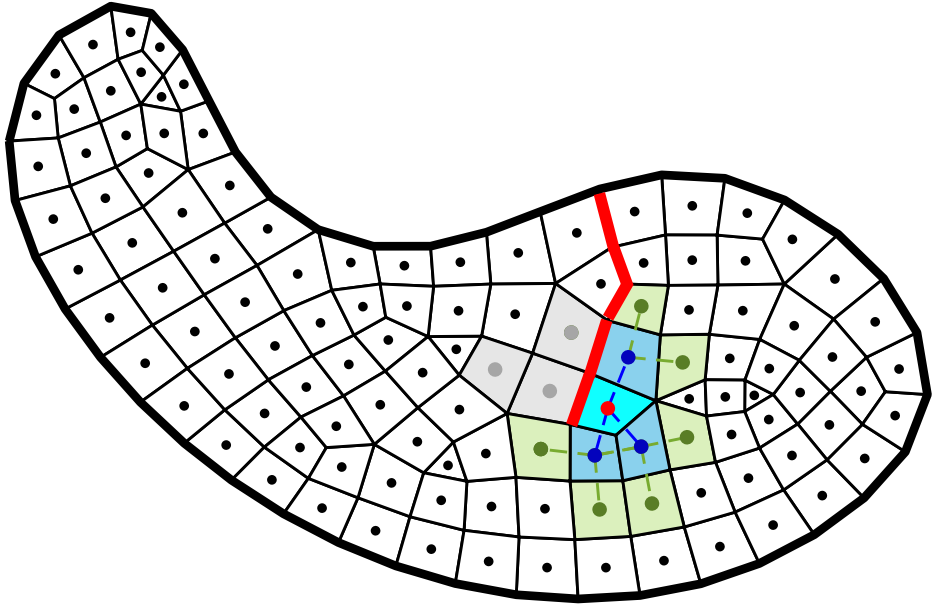}
  \caption{Support of a point by the sides of cracks.} 
  \label{fig:Schem_Crack} 
\end{figure}

When an electrically impermeable crack emerges between two adjacent subdomains, we convert the corresponding internal boundary to two external boundaries with zero mechanical and electrical tractions  ($e \in \partial \Omega_Q \cap \partial \Omega_R \cap \partial \Omega_\omega \cap \partial \Omega_Z$). And we also cut off the interaction between the two Points in the adjacent subdomains, i.e., the Points will be removed from the neighboring point set of each other. Therefore, the support of the point on one side of the crack no longer contain the adjacent points and its neighboring points on the other side of the crack. For example, in Fig.~\ref{fig:Schem_Crack}, as a result of the crack (presented by a red polyline), three points (grey) are removed from the support of the reference point (red).

As for the electrically permeable crack-face boundaries, two different interpolations are employed in the mechanical and electrical field. In the mechanical field, the internal boundaries on the cracks are converted into traction-free boundaries and interactions between the points on either side of the crack are cut off. Whereas in the electrical interpolation, these boundaries remain as internal boundaries with continuous electrical field. As a result, the support of points and definition of $\Gamma_h$ may be different in the mechanical and electrical interpolations.

Therefore, for both the crack-face boundaries conditions, only slight adjustments are required for the stiffness matrix and/or load vectors when a crack occurs. That is, for the two points adjacent to the crack and their closest neighbors, the point stiffness matrices will be regenerated, and the terms related to the numerical fluxes on the cracked boundaries should be deleted. The total number of DoF remains the same. And in the linear reduced theory, the load vector does not need to be adjusted. This is much more convenient than remeshing or deleting elements as in the FEM when simulating the crack propagation.

Various crack initiation and propagation criteria can be applied and incorporated within the FPM. The choice of these criteria should be based on the characteristics of the realistic problem and material. In this paper, for example, the classic Maximum Hoop Stress criterion and an Inter-Subdomain-Boundary Bonding-Energy-Rate(BER)-Based criterion are employed.

In crack initiation analysis, we consider a quasi-static loading process. The entire load history can be divided into several steps. And in each step, for simplicity, the most ``dangerous'' internal boundaries will be cracked. That is, for instance, based on the Maximum Hoop Stress criterion, the internal boundary with the largest normal traction will be cracked. Alternatively, here we introduce an energy-based criterion initiated by the J-integral. Based on the electric Gibbs free energy $G$, the J-integral in domain $\Omega^{int}$ is defined as \cite{Sladek2017}:
\begin{align}\label{eqn:J_1} 
\begin{split}
J = \int_{\partial \Omega^{int}} G {n}_1 \mathrm{d} \Gamma - \int_{\partial \Omega^{int}} \left[ \begin{matrix} \left( \frac{\partial \mathbf{u}}{ \partial \hat{x}_1} \right)^\mathrm{T} \widetilde{\mathbf{Q}} + \left( \frac{\partial \left( \overline{\mathbf{n}}_4^e  \hat{ \boldsymbol{\varepsilon}} \left(\mathbf{u} \right) \right)}{ \partial \hat{x}_1} \right)^\mathrm{T} \widetilde{\mathbf{R}} \\ -  \frac{\partial \phi}{ \partial \hat{x}_1} \widetilde{\omega} - \left( \frac{\partial \left( {\mathbf{n}}^{e \mathrm{T}}  \mathbf{E} \left(\phi \right) \right)}{ \partial \hat{x}_1} \right)^\mathrm{T} \widetilde{Z} \end{matrix} \right] \mathrm{d} \Gamma.
\end{split}
\end{align}
The previous equation is written in a local $\hat{x}_1$-$\hat{x}_2$ coordinate system, in which $\hat{x}_1$ is aligned with the crack and  $\hat{x}_2$ is aligned with the normal vector. $\widetilde{\mathbf{Q}}$, $\widetilde{\mathbf{R}}$, $\widetilde{\omega}$ and $\widetilde{Z}$ are defined in Eqn.~\ref{eqn:BC_full_5} -- \ref{eqn:BC_full_8}, yet they are unknown at here. When $\Omega^{int}$ is divided into several subdomains $\Omega^{int}_I$, we can rewrite Eqn.~\ref{eqn:J_1} as:
\begin{align}\label{eqn:J_2} 
\begin{split}
J = &  \sum_I \int_{\Omega_I^{int}} \left[ \frac{\partial \hat{\boldsymbol{\varepsilon}} \left( \mathbf{u} \right)}{\partial \hat{x}_1} \frac{\partial G}{\partial  \hat{\boldsymbol{\varepsilon}} } + \frac{\partial \boldsymbol{\kappa} \left( \mathbf{u} \right) }{\partial \hat{x}_1} \frac{\partial G}{\partial  \boldsymbol{\kappa} } + \frac{\partial \mathbf{E} \left( \phi \right)}{\partial \hat{x}_1} \frac{\partial G}{\partial \mathbf{E} }+ \frac{\partial \mathbf{V} \left( \phi \right)}{\partial \hat{x}_1} \frac{\partial G}{\partial  \mathbf{V} } \right] \mathrm{d} \Omega \\
& - \int_{\partial \Omega^{int}} \left[ \left( \frac{\partial \mathbf{u}}{ \partial \hat{x}_1} \right)^\mathrm{T} \widetilde{\mathbf{Q}} + \frac{\partial \hat{ \boldsymbol{\varepsilon}}^\mathrm{T} \left( \mathbf{u} \right)}{ \partial \hat{x}_1} \overline{\mathbf{n}}_4^{e \mathrm{T}}  \widetilde{\mathbf{R}} -  \frac{\partial \phi}{ \partial \hat{x}_1} \widetilde{\omega} - \frac{ \partial \mathbf{E}^\mathrm{T} \left( \phi \right)}{\partial \hat{x}_1} {\mathbf{n}}^{e}  \widetilde{Z} \right] \mathrm{d} \Gamma \\
&  -  \sum \int_{\Gamma_h^{int}} \left[ \! \left[ G n_1^e \right] \! \right] \mathrm{d} \Gamma\\
= &  \sum_I \int_{\Omega_I^{int}} \left[ \frac{\partial \boldsymbol{\varepsilon} \left( \mathbf{u} \right)}{\partial \hat{x}_1} \boldsymbol{\sigma}  +  \frac{\partial \hat{\boldsymbol{\varepsilon}} \left( \mathbf{u} \right)}{\partial \hat{x}_1} \boldsymbol{\sigma}^{ES} +  \frac{\partial \boldsymbol{\kappa} \left( \mathbf{u} \right) }{\partial \hat{x}_1} \boldsymbol{\mu} - \frac{\partial \mathbf{E} \left( \phi \right)}{\partial \hat{x}_1} \mathbf{D}+ \frac{\partial \mathbf{V} \left( \phi \right)}{\partial \hat{x}_1} \mathbf{Q}\right] \mathrm{d} \Omega \\
& - \int_{\partial \Omega^{int}} \left[ \left( \frac{\partial \mathbf{u}}{ \partial \hat{x}_1} \right)^\mathrm{T} \widetilde{\mathbf{Q}} + \frac{\partial \hat{ \boldsymbol{\varepsilon}}^\mathrm{T} \left( \mathbf{u} \right)}{ \partial \hat{x}_1} \overline{\mathbf{n}}_4^{e \mathrm{T}}  \widetilde{\mathbf{R}} -  \frac{\partial \phi}{ \partial \hat{x}_1} \widetilde{\omega} - \frac{ \partial \mathbf{E}^\mathrm{T} \left( \phi \right)}{\partial \hat{x}_1} {\mathbf{n}}^{e}  \widetilde{Z} \right] \mathrm{d} \Gamma \\
&  -  \sum \int_{\Gamma_h^{int}} \left[ \! \left[ G n_1^e \right] \! \right] \mathrm{d} \Gamma
\end{split}
\end{align}
where $\Gamma_h^{int}$ is the set of internal boundaries within $\Omega^{int}$. Since $\hat{x}_1$ is aligned with the crack, $n_1^e = 0$ thus the last term in Eqn.~\ref{eqn:J_2} can be eliminated.

The FPM formulations Eqn.~\ref{eqn:dis_full} and \ref{eqn:ele_full} should be satisfied in the present domain $\Omega^{int}$ with any arbitrary test functions $\mathbf{v}$ and $\tau$. Therefore, here we let $\mathbf{v} = \frac{\partial \mathbf{u}}{ \partial \hat{x}_1}$ and $\tau = \frac{\partial \phi}{ \partial \hat{x}_1}$ and comparing the FPM equations with Eqn~\ref{eqn:J_2}. The essential boundaries are omitted here. Finally, we can achieve:
\begin{align}
\begin{split}
J \approx \sum_{e \in \Gamma_h^{int}} BER,
\end{split}
\end{align}
where
\begin{align}
\begin{split}
BER = & - \int_e \frac{\eta_{21}}{h_e} \left[ \! \left[ \frac{ \partial \mathbf{u}^\mathrm{T}}{\partial \hat{x}_1} \right] \! \right] \left[ \!\left[ \mathbf{u} \right] \! \right] \mathrm{d} \Gamma - \int_e \eta_{22} h_e \left[ \! \left[  \frac{ \partial \hat{\boldsymbol{\varepsilon}}^\mathrm{T} \left( \mathbf{u} \right)}{\partial \hat{x}_1} \right] \! \right] \overline{\mathbf{n}}_4^{e \mathrm{T}} \overline{\mathbf{n}}_4^e \left[ \! \left[ \hat{\boldsymbol{\varepsilon}} \left( \mathbf{u} \right) \right] \! \right] \mathrm{d} \Gamma \\
& \ - \int_e \frac{\eta_{23}}{h_e} \left[ \! \left[   \frac{ \partial \phi}{\partial \hat{x}_1}  \right] \! \right] \left[ \!\left[ \phi \right] \! \right] \mathrm{d} \Gamma - \int_e \eta_{24} h_e \left[ \! \left[ \frac{ \partial \mathbf{E}^\mathrm{T} \left( \phi \right)}{\partial \hat{x}_1} \mathbf{n}^{e} \right] \! \right] \overline{\mathbf{n}}_4^{e \mathrm{T}} \overline{\mathbf{n}}_4^e \left[ \! \left[ \mathbf{n}^{e \mathrm{T}}  \mathbf{E} \left( \phi \right)  \right] \! \right] \mathrm{d} \Gamma\\
&   +\int_e \left( \left[ \! \left[ \frac{ \partial \mathbf{u}^\mathrm{T}}{\partial \hat{x}_1} \right] \! \right] \overline{\mathbf{n}}_1^e \mathbf{C}_{\sigma \varepsilon} \left\{ \boldsymbol{\varepsilon} \left( \mathbf{u} \right) \right\} + \left\{ \frac{ \partial \boldsymbol{\varepsilon}^\mathrm{T} \left( \mathbf{u} \right)}{\partial \hat{x}_1} \right\} \mathbf{C}_{\sigma \varepsilon} \overline{\mathbf{n}}_1^{e \mathrm{T}} \left[ \! \left[ \mathbf{u} \right] \! \right]  \right) \mathrm{d} \Gamma \\
& +\int_e \left( \left[ \! \left[ \frac{ \partial \hat{\boldsymbol{\varepsilon}}^\mathrm{T} \left( \mathbf{u} \right)}{\partial \hat{x}_1} \right] \! \right] \overline{\mathbf{n}}_3^e \mathbf{C}_{\mu \kappa} \left\{ \boldsymbol{\kappa} \left( \mathbf{u} \right) \right\} + \left\{  \frac{ \partial {\boldsymbol{\kappa}}^\mathrm{T} \left( \mathbf{u} \right)}{\partial \hat{x}_1} \right\} \mathbf{C}_{\mu \kappa} \overline{\mathbf{n}}_3^{e \mathrm{T}} \left[ \! \left[ \hat{\boldsymbol{\varepsilon}} \left( \mathbf{u} \right) \right] \! \right] \right)  \mathrm{d} \Gamma \\
&   + \int_e \left(  \left[ \! \left[ \frac{ \partial \phi}{\partial \hat{x}_1} \right] \! \right] \mathbf{n}^{e \mathrm{T}} \boldsymbol{\Lambda} \left\{ \mathbf{E} \left( \phi \right) \right\} + \left\{ \frac{ \partial \mathbf{E}^\mathrm{T} \left( \phi \right)}{\partial \hat{x}_1} \right\} \boldsymbol{\Lambda}  \mathbf{n}^{e} \left[ \! \left[ \phi \right] \! \right] \right) \mathrm{d} \Gamma \\
&   - \int_e \left( \left[ \! \left[ \frac{ \partial \mathbf{E}^\mathrm{T} \left( \phi \right)}{\partial \hat{x}_1} \right] \! \right] \overline{\mathbf{n}}_{6}^e \boldsymbol{\Phi} \left\{ \mathbf{V} \left( \phi \right) \right\} + \left\{ \frac{ \partial \mathbf{V}^\mathrm{T} \left( \phi \right)}{\partial \hat{x}_1} \right\} \boldsymbol{\Phi} \overline{\mathbf{n}}_{6}^{e \mathrm{T}} \left[ \! \left[ \mathbf{E} \left( \phi \right) \right] \! \right] \right) \mathrm{d} \Gamma \\
&   - \int_e \left[ \! \left[ \frac{ \partial \mathbf{u}^\mathrm{T}}{\partial \hat{x}_1} \right] \! \right] \overline{\mathbf{n}}_{21}^e \mathbf{C}_{\mu \kappa} \left\{ \boldsymbol{\kappa}_{,1} \left( \mathbf{u} \right) \right\} \mathrm{d} \Gamma - \int_e \left[ \! \left[ \frac{ \partial \mathbf{u}^\mathrm{T}}{\partial \hat{x}_1} \right] \! \right] \overline{\mathbf{n}}_{22}^e \mathbf{C}_{\mu \kappa} \left\{ \boldsymbol{\kappa}_{,2} \left( \mathbf{u} \right\} \right) \mathrm{d} \Gamma \\
&    - \int_e \left[ \! \left[ \frac{ \partial \mathbf{u}^\mathrm{T}}{\partial \hat{x}_1} \right] \! \right] \left\{ \left( \overline{\mathbf{n}}_1^e \mathbf{b} - \overline{\mathbf{n}}_{21}^e \mathbf{a} \overline{\mathbf{c}}_5 - \overline{\mathbf{n}}_{22}^e \mathbf{a} \overline{\mathbf{c}}_6 - \overline{\mathbf{n}}_4^e  \widetilde{\mathbf{Q}} \left( \mathbf{u}, \phi \right) \right) \mathbf{V} \left( \phi \right) \right\} \mathrm{d} \Gamma \\
&  - \int_e  \left\{ \frac{ \partial \mathbf{V}^\mathrm{T} \left( \phi \right)}{\partial \hat{x}_1}  \left( \mathbf{b}^\mathrm{T} \overline{\mathbf{n}}_1^{e \mathrm{T}} - \overline{\mathbf{c}}_5^\mathrm{T} \mathbf{a}^\mathrm{T} \overline{\mathbf{n}}_{21}^{e \mathrm{T}} - \overline{\mathbf{c}}_6^\mathrm{T} \mathbf{a}^\mathrm{T} \overline{\mathbf{n}}_{22}^{e \mathrm{T}} - \widetilde{\mathbf{Q}}^\mathrm{T} \left( \mathbf{u}, \phi \right) \overline{\mathbf{n}}_4^{e \mathrm{T}} \right) \right\} \left[ \! \left[ \mathbf{u} \right] \! \right] \mathrm{d} \Gamma \\
&   - \int_e \left[ \! \left[ \frac{ \partial \mathbf{u}^\mathrm{T}}{\partial \hat{x}_1} \right] \! \right] \left\{ \left( \overline{\mathbf{n}}_1^e \mathbf{e} - \overline{\mathbf{n}}_4^e  \widetilde{\mathbf{D}} \left( \mathbf{u}, \phi \right) \right) \mathbf{E} \left( \phi \right) \right\} \mathrm{d} \Gamma  + \int_e \left\{ \frac{ \partial \boldsymbol{\varepsilon}^\mathrm{T} \left( \mathbf{u} \right)}{\partial \hat{x}_1} \right\}  \mathbf{e} \mathbf{n}^{e}   \left[ \! \left[ \phi \right] \! \right] \mathrm{d} \Gamma \\
&- \int_e  \left\{  \frac{ \partial \mathbf{E}^\mathrm{T} \left( \phi \right)}{\partial \hat{x}_1} \left( \mathbf{e}^\mathrm{T} \overline{\mathbf{n}}_1^{e \mathrm{T}} - \widetilde{\mathbf{D}}^\mathrm{T} \left( \mathbf{u}, \phi \right) \overline{\mathbf{n}}_4^{e \mathrm{T}} \right) \right\} \left[ \! \left[ \mathbf{u} \right] \! \right] \mathrm{d} \Gamma + \int_e \left[ \! \left[ \frac{ \partial \phi}{\partial \hat{x}_1} \right] \! \right] \mathbf{n}^{e \mathrm{T}} \mathbf{e}^\mathrm{T} \left\{ \boldsymbol{\varepsilon} \left( \mathbf{u} \right) \right\} \mathrm{d} \Gamma \\
& + \int_e \left\{ \frac{ \partial \boldsymbol{\kappa}^\mathrm{T} \left( \mathbf{u} \right)}{\partial \hat{x}_1} \right\}  \left( \mathbf{a} \mathbf{n}^{e} -  \overline{\mathbf{c}}_7^\mathrm{T} \mathbf{b} \overline{\mathbf{n}}_{51}^{e \mathrm{T}} -  \overline{\mathbf{c}}_8^\mathrm{T} \mathbf{b} \overline{\mathbf{n}}_{52}^{e \mathrm{T}}  \right)  \left[ \! \left[ \phi \right] \! \right] \mathrm{d} \Gamma \\
&  + \int_e \left[ \! \left[ \frac{ \partial \phi}{\partial \hat{x}_1} \right] \! \right] \left( \mathbf{n}^{e \mathrm{T}} \mathbf{a}^\mathrm{T} - \overline{\mathbf{n}}_{51}^e \mathbf{b}^\mathrm{T} \overline{\mathbf{c}}_7 - \overline{\mathbf{n}}_{52}^e \mathbf{b}^\mathrm{T} \overline{\mathbf{c}}_8 \right) \left\{ \boldsymbol{\kappa} \left( \mathbf{u} \right) \right\} \mathrm{d} \Gamma \\
& -\int_e \left[ \! \left[  \frac{ \partial \hat{\boldsymbol{\varepsilon}}^\mathrm{T} \left( \mathbf{u} \right)}{\partial \hat{x}_1} \right] \! \right] \overline{\mathbf{n}}_3^e \mathbf{a} \left\{ \mathbf{E} \left( \phi \right) \right\}  \mathrm{d} \Gamma - \int_e \left\{  \frac{ \partial {\boldsymbol{\varepsilon}}^\mathrm{T} \left( \mathbf{u} \right)}{\partial \hat{x}_1} \right\} \mathbf{b} \overline{\mathbf{n}}_{6}^{e \mathrm{T}} \left[ \! \left[ \mathbf{E} \left( \phi \right) \right] \! \right] \mathrm{d} \Gamma \\
&  - \int_e \left\{ \frac{ \partial \mathbf{E}^\mathrm{T} \left( \phi \right)}{\partial \hat{x}_1} \right\} \mathbf{a}^\mathrm{T} \overline{\mathbf{n}}_3^{e \mathrm{T}} \left[ \! \left[ \hat{\boldsymbol{\varepsilon}} \left( \mathbf{u} \right) \right] \! \right]   \mathrm{d} \Gamma  - \int_e \left[ \! \left[ \frac{ \partial \mathbf{E}^\mathrm{T} \left( \phi \right)}{\partial \hat{x}_1} \right] \! \right] \overline{\mathbf{n}}_{6}^e \mathbf{b}^\mathrm{T} \left\{\boldsymbol{\varepsilon} \left( \mathbf{u} \right) \right\}  \mathrm{d} \Gamma  \\
&  - \int_e \left(  \left[ \! \left[ \frac{ \partial \phi}{\partial \hat{x}_1} \right] \! \right] \overline{\mathbf{n}}_{51}^e \boldsymbol{\Phi} \left\{ \mathbf{V}_{,1} \left( \phi \right) \right\} \right) \mathrm{d} \Gamma - \int_e \left(  \left[ \! \left[ \frac{ \partial \phi}{\partial \hat{x}_1} \right] \! \right] \overline{\mathbf{n}}_{52}^e \boldsymbol{\Phi} \left\{ \mathbf{V}_{,2} \left( \phi \right) \right\} \right) \mathrm{d} \Gamma .
\end{split}
\end{align}

Thus, the J-integral can be approximated by a sum of integrals over all the internal boundaries. Here we postulate the component on each internal boundary as an estimate of the bonding energy rate ($BER$). In the current work, the $BER$ is employed as an energy-based criterion of crack initiation and development. An internal boundary segment will be cracked if its $BER$ exceeds a prescribed critical value. Note that in the FPM, the $BER$ on each internal boundary is generated following the same process of the internal boundary stiffness matrix $\mathbf{K}_h$, except that the test functions $\mathbf{v}$ and $\tau$ are replaced by $\frac{\partial \mathbf{u}}{ \partial \hat{x}_1}$ and $\frac{\partial \phi}{ \partial \hat{x}_1}$, where $\mathbf{u}$ and $\phi$ are solutions of the FPM analysis.

\section{Conclusion}

A meshless Fragile Points Method (FPM) based on Galerkin weak form is developed to analyze flexoelectric effects in dielectric solids. Both primal and mixed formulations are presented. Local, simple, polynomial and discontinuous trial and test functions, which take the form of Taylor expansions at each internal Fragile Point, are employed. Local RBF-DQ method is applied to approximate the higher derivatives. Numerical Flux Corrections are ultilized to ensure the continuity condition. In the mixed FPM, all the additional high order independent variables are eliminated locally at each Point, and only 3 DoFs are retained explicitly in the final formula. Both full and reduced flexoelectric theories are taken into consideration. The full theory involving the electrostatic stress leads a nonlinear algebraic system while the formulations for the reduced theory are linear. The algorithm for simulating crack initiation and propagation in the FPM is also presented, using the stress-based criterion as well as an Inter-Subdomain-Boundary Bonding-Energy-Rate(BER)-Based criterion for crack initiation and development.

To conclude, the primal and mixed FPM proposed in this work have the following features and advantages as compared to other computational methods for flexoelectricity published in prior literature:
\begin{itemize}
  \item The FPM is a meshless or element-free method with point-based trial and test functions that are  piece-wise continuous.
  \item The integration of the Galerkin weak form is very simple by using Guassian quadrature with only one integration point being sufficient most of the time.
  \item The trial function has the Kronecker-delta property; the essential boundary conditions can be imposed directly or weakly using numerical fluxes.
  \item Crack and rupture initiation and propagation can be easily simulated by FPM and do not involve remeshing or trial function enhancement.
  \item Crack and rupture criteria are based on simple and commonly used continuum physics.
  \item The stiffness matrix is sparse and symmetric.
  \item With mixed FPM, the first and higher strain gradient problems can be handled with ease. Only the lowest-order primitive variables are retained as DoFs at each node.
  \item There is no incompressibility locking or shear locking in FPM.
  \item The FPM does not suffer from mesh distortion.
\end{itemize}

In this first part of the present two-paper series, the theoretical formulation and implementation of the proposed methodology are given in detail. Numerical examples and validation are reported in Part II of the present study.

\section*{References}
\bibliography{Part_1}

\section*{Acknowledgment}

Yue Guan thankfully acknowledges the financial support for her work, provided through the funding for Professor Atluri’s Presidential Chair at TTU.

\appendix
\section{Matrices in numerical implementation} \label{app: Mats}

\begin{align}
\begin{split}
\overline{\mathbf{c}}_1 & = \left[ \begin{matrix} 1 & 0 & 0 & 0 \\ 0 & 0 & 0 & 1 \end{matrix} \right], \qquad \  \  \, \qquad \quad \overline{\mathbf{c}}_2 = \left[ \begin{matrix} 0 & 0 & 1 & 0 \\ 0 & 1 & 0 & 0 \end{matrix} \right] \\
\overline{\mathbf{c}}_3 & = \left[ \begin{matrix} 1 & 0 \end{matrix} \right], \quad \overline{\mathbf{c}}_4 = \left[ \begin{matrix} 0 & 1 \end{matrix} \right] \qquad \ \ \overline{\mathbf{c}}_5 = \overline{\mathbf{c}}_3 \otimes \mathbf{I}_{2 \times 2}, \quad \overline{\mathbf{c}}_6 = \overline{\mathbf{c}}_4 \otimes \mathbf{I}_{2 \times 2}, \\
\overline{\mathbf{c}}_7 & = \left[ \begin{matrix} 1 & 0 & 0 & 0 & 0 & 0 \\ 0 & 0 & 0 & 0 & 0 & 1/2 \\ 0 & 0 & 0 & 1 & 1/2 & 0 \end{matrix} \right], \quad \overline{\mathbf{c}}_8 = \left[ \begin{matrix} 0 & 0 & 0 & 0 & 1/2 & 0 \\ 0 & 1 & 0 & 0 & 0 & 0 \\ 0 & 0 & 1 & 0 & 0 & 1/2 \end{matrix} \right] \\
\overline{\mathbf{c}}_9 & = \left[ \begin{matrix} 1 & 0 & 0 & 0 & 0 & 0 \\ 0 & 0 & 0 & 0 & 0 & 1 \\ 0 & 0 & 0 & 0 & 1 & 0 \\ 0 & 0 & 0 & 1 & 0 & 0 \end{matrix} \right], \qquad \; \; \overline{\mathbf{c}}_{10} = \left[ \begin{matrix} 0 & 0 & 0 & 0 & 1 & 0 \\ 0 & 1 & 0 & 0 & 0 & 0 \\ 0 & 0 & 1 & 0 & 0 & 0 \\ 0 & 0 & 0 & 0 & 0 & 1 \end{matrix} \right] \\
\mathbf{n} & = \left[ \begin{matrix} n_1 &  n_2 \end{matrix} \right]^\mathrm{T}, \qquad \qquad \qquad \quad \ \ \, \mathbf{s} = \left[ \begin{matrix} s_1 & s_2 \end{matrix} \right]^\mathrm{T} = \pm \left[ \begin{matrix} n_2 & -n_1 \end{matrix} \right]^\mathrm{T} \\
\overline{\mathbf{n}}_{1} & = \left[ \begin{matrix} n_1 & 0 & n_2 \\ 0 & n_2 &  n_1 \end{matrix} \right] , \qquad \qquad \quad \ \,  \overline{\mathbf{n}}_{6} = \left[ \begin{matrix} n_1 & 0 & n_2 & 0 \\ 0 & n_1 & 0 & n_2 \end{matrix} \right], \\
\overline{\mathbf{n}}_{21} & = \left[ \begin{matrix} n_1 & 0 & 0 & 0 & n_2 & 0  \\ 0 & 0 & 0 &  n_1 & 0 & n_2  \end{matrix} \right], \ \ \overline{\mathbf{n}}_{22} = \left[ \begin{matrix} 0 & 0 & n_2 & 0 & n_1 & 0 \\ 0 & n_2 & 0 &  0 & 0 & n_1 \end{matrix} \right], \\
\overline{\mathbf{n}}_{3} & = \left[ \begin{matrix} n_1 & 0 & 0 & 0 & n_2 & 0 \\ 0 & n_2 & 0 &  0 & 0 & n_1 \\ 0 & 0 & n_2 & 0 & n_1 & 0 \\  0 & 0 & 0 &  n_1 & 0 & n_2 \end{matrix} \right], \\
\overline{\mathbf{n}}_{4} & = \left[ \begin{matrix} n_1 & 0 & n_2 & 0 \\ 0 & n_2 & 0 &  n_1 \end{matrix} \right],  \qquad \quad \ \ \, \overline{\mathbf{s}}_{4} = \left[ \begin{matrix} s_1 & 0 & s_2 & 0 \\ 0 & s_2 & 0 &  s_1 \end{matrix} \right], \\
\overline{\mathbf{n}}_{51} & = \left[ \begin{matrix} n_1 & n_2 & 0 & 0 \end{matrix} \right], \qquad \  \ \,  \qquad \overline{\mathbf{n}}_{52} = \left[ \begin{matrix} 0 & 0 & n_1 & n_2 \end{matrix} \right], \\ 
\end{split}
\end{align}

\end{document}